\documentclass[12pt]{article}
\usepackage{amsmath}
\usepackage{amssymb}
\usepackage[textures]{graphics}
\usepackage{color}
\usepackage{multicol}
\renewcommand{\baselinestretch}{2}

\newtheorem{Def}{Definition}[section]
\newtheorem{Lem}[Def]{Lemma}
\newtheorem{theorem}[Def]{Theorem}
\newtheorem{Proposition}[Def]{Proposition}

\newtheorem{example}[Def]{Example}

\def\a{{\langle\eta\rangle}}
\def\ai{{\langle\eta\rangle^{-1}}}
\def\aa{{ |\eta|}}
\def\aai{{ |\eta|^{-1}}}
\def\aaa{{ \|\eta\|}}
\def\b{{\langle\xi\rangle}}
\def\bi{{\langle\xi\rangle^{-1}}}
\def\bb{{ |\xi|}}
\def\bbi{{ |\xi|^{-1}}}
\def\bbb{{ \|\xi\|}}
\def\ab{{\xi\otimes\eta }}
\def\achl{{\hat{\triangleleft}}}

\def\acbl{{\bar{\triangleleft}}}
\def\acbr{{\bar{\triangleright}}}
\def\actl{{\tilde{\triangleleft}}}
\def\actr{{\tilde{\triangleright}}}
\def\acl{{\triangleleft}}
\def\acr{{\triangleright}}
\def\al{{\alpha}}

\def\dacr{{\dot \acr}}
\def\aabar{{ |\bar\eta|}}
\def\aap{{ |\eta^{'}|}}
\def\ap{{\langle\eta^{'}\rangle}}
\def\eb{{ \bar\eta}}
\def\epr{{ \eta^{'}}}
\def\T{{ \setlength{\unitlength}{0.09cm}
\begin{picture}(5,2)\thicklines
\put(1,1){\circle{9} }
\put(-1.0,-0.5){\text {$T$}}
\end{picture}}}
\def\thi{{ \setlength{\unitlength}{0.09cm}
\begin{picture}(5,2)\thicklines
\put(1,1){\circle{9} }
\put(-1.9,-0.55){\text {$\theta^{-1}$}}
\end{picture}}}
\def\sthi{{ \setlength{\unitlength}{0.09cm}
\begin{picture}(5,2)\thicklines
\put(1,1){\circle{8} }
\put(-1.9,-01.25){\text {$\theta^{-1}$}}
\end{picture}}}
\def\S{{ \setlength{\unitlength}{0.09cm}
\begin{picture}(5,2)\thicklines
\put(1,1){\circle{7} }
\put(-1.0,-0.5){\text {$S$}}
\end{picture}}}
\def\u{{ \setlength{\unitlength}{0.09cm}
\begin{picture}(5,2)\thicklines
\put(1,1){\circle{7} }
\put(-1.0,-0.5){\text {u}}
\end{picture}}}
\def\uu{{ \setlength{\unitlength}{0.5cm}
\begin{picture}(5,2)\thicklines
\put(5.3,-3.25){\scoev}
\put(6.03,-2.){\seval}
\put(3.45,-2.9){\lbraid}
\put(9.95,-3.3){\line(0,1){3.}}
\put(6.37,0.38){\line(0,1){2.}}
\end{picture}}}

\def\Ad{{ \setlength{\unitlength}{0.09cm}
\begin{picture}(5,2)\thicklines
\put(1,1){\circle{10} }
\put(-2.0,-0.5){\text {$Ad$}}
\end{picture}}}
\def\ch{{ \setlength{\unitlength}{0.09cm}
\begin{picture}(5,2)\thicklines
\put(1,1){\circle{8} }
\put(-1.3,0.5){\text {${\chi_{_V}}$}}
\end{picture}}}
\def\chw{{ \setlength{\unitlength}{0.09cm}
\begin{picture}(5,2)\thicklines
\put(1,1){\circle{8} }
\put(-1.6,0.5){\text {${\chi_{_W}}$}}
\end{picture}}}
\def\chvw{{ \setlength{\unitlength}{0.09cm}
\begin{picture}(5,2)\thicklines
\put(1,1){\circle{12} }
\put(-3.7,0.7){\text {${\chi_{_{V \otimes W}}}$}}
\end{picture}}}
\def\t{{ \setlength{\unitlength}{0.09cm}
\begin{picture}(5,2)\thicklines
\put(1,1){\circle{9} }
\put(-1.0,-0.5){\text {$\theta_V$}}
\end{picture}}}
\def\tt{{ \setlength{\unitlength}{0.09cm}
\begin{picture}(5,2)\thicklines
\put(1,1){\circle{9} }
\put(-1.0,-0.5){\text {$\theta_{V^*}$}}
\end{picture}}}
\def\w{{\setlength{\unitlength}{0.09cm}
\begin{picture}(5,2)\thicklines
\put(1,1){\circle{9} }
\put(-1.0,-0.5){\text {$\theta_W$}}
\end{picture} }}
\def\wi{{ \setlength{\unitlength}{0.09cm}
\begin{picture}(5,2)\thicklines
\put(1,1){\circle{11} }
\put(-3.0,-1.5){\text {${\theta_W}^{\!\!-1}$}}
\end{picture}}}
\def\ti{{ \setlength{\unitlength}{0.09cm}
\begin{picture}(5,2)\thicklines
\put(1,1){\circle{11} }
\put(-3.0,-1.5){\text {${\theta_V}^{\!-1}$}}
\end{picture}}}
\def\tti{{ \setlength{\unitlength}{0.09cm}
\begin{picture}(5,2)\thicklines
\put(1,1){\circle{11} }
\put(-3.4,-1.5){\text {${\theta_{V^*}}^{\!\!\!-1}$}}
\end{picture}}}
\def\thsi{{ \setlength{\unitlength}{0.5cm}
\begin{picture}(5,2)\thicklines
\put(16.8,-7.5){\ti}
\put(15.,-9.08){\seval}
\put(15.05,-8.79){\scoev}
\put(17.53,-6.35){\line(0,1){0.5}}
\put(17.53,-8.8){\line(0,1){0.5}}
\put(15.35,-8.8){\line(0,1){4}}
\put(19.75,-9.85){\line(0,1){4.}}
\end{picture}}}
\def\vw{{ \setlength{\unitlength}{0.09cm}
\begin{picture}(5,2)\thicklines
\put(1,1){\circle{12} }
\put(-4.3,-0.2){\text {$\theta_{V\otimes W}$}}
\end{picture}}}
\def\ep{{ \setlength{\unitlength}{0.09cm}
\begin{picture}(5,2)\thicklines
\put(1,1){\circle{7} }
\put(-.30,0){\text {$\epsilon$}}
\end{picture}}}
\def\et{{ \setlength{\unitlength}{0.09cm}
\begin{picture}(5,2)\thicklines
\put(1,1){\circle{7} }
\put(-.50,0){\text {$\eta$}}
\end{picture}}}
\def\ltri{{ \setlength{\unitlength}{0.4cm}
\begin{picture}(5,2)\thicklines
\put(0.1,0.4){\line(0,1){1}}
\put(7.9,0.4){\line(0,1){1}}
\put(0.1,0.4){\line(1,0){7.8}}
\put(0.1,1.42){\line(1,0){7.8}}
\end{picture}}}
\def\stri{{ \setlength{\unitlength}{0.4cm}
\begin{picture}(5,2)\thicklines
\put(0.1,0.4){\line(0,1){1}}
\put(2.0,0.4){\line(0,1){1}}
\put(0.1,0.4){\line(1,0){1.9}}
\put(0.1,1.42){\line(1,0){1.9}}
\end{picture}}}
\def\lbraid{{ \setlength{\unitlength}{0.6cm}
\begin{picture}(5,2)\thicklines
\put(2,1){\line(1,1){.75}}
\put(3,2){\line(1,1){.75}}
\put(3.75,1){\line(-1,1){1.75}}
\end{picture}}}
\def\braid{{ \setlength{\unitlength}{0.5cm}
\begin{picture}(5,2)\thicklines
\put(2,1){\line(1,1){.75}}
\put(3,2){\line(1,1){.75}}
\put(3.75,1){\line(-1,1){1.75}}
\end{picture}}}
\def\braidi{{\setlength{\unitlength}{0.5cm}
\begin{picture}(5,2)\thicklines
\put(2,1){\line(1,1){1.75}}
\put(3.75,1){\line(-1,1){.75}}
\put(2.75,2){\line(-1,1){.75}}
\end{picture}}}
\def\lbraidi{{\setlength{\unitlength}{0.6cm}
\begin{picture}(5,2)\thicklines
\put(2,1){\line(1,1){1.75}}
\put(3.75,1){\line(-1,1){.75}}
\put(2.75,2){\line(-1,1){.75}}
\end{picture}}}
\def\seval{{\setlength{\unitlength}{0.36cm}
\begin{picture}(5,2)\thicklines
\put(0.1,0.4){\line(1,-1){1}}
\put(1.1,-0.6){\line(1,0){1}}
\put(2.1,-0.6){\line(1,1){1}}
\end{picture} }}
\def\lcoev{{\setlength{\unitlength}{0.7cm}
\begin{picture}(5,2)\thicklines
\put(3.1,4.1){\line(1,1){1}}
\put(4.1,5.1){\line(1,0){1}}
\put(5.1,5.1){\line(1,-1){1}}
\end{picture}}}
\def\leval{{\setlength{\unitlength}{0.7cm}
\begin{picture}(5,2)\thicklines
\put(0.1,0.4){\line(1,-1){1}}
\put(1.1,-0.6){\line(1,0){1}}
\put(2.1,-0.6){\line(1,1){1}}
\end{picture} }}
\def\scoev{{\setlength{\unitlength}{0.36cm}
\begin{picture}(5,2)\thicklines
\put(3.1,4.1){\line(1,1){1}}
\put(4.1,5.1){\line(1,0){1}}
\put(5.1,5.1){\line(1,-1){1}}
\end{picture}}}
\def\brai{{\setlength{\unitlength}{0.5cm}
\begin{picture}(5,2)\thicklines
\put(2,1){\line(6,5){1.7}}
\put(2.7,1.9){\line(-4,5){.57}}
\put(3.6,.9){\line(-4,5){.57}}
\put(3.6,.3){\line(0,1){.6}}
\put(2,.3){\line(0,1){.7}}
\put(2.1,2.6){\line(0,1){.6}}
\put(3.71,2.4){\line(0,1){.77}}
\end{picture} }}
\def\eval{{\setlength{\unitlength}{0.5cm}
\begin{picture}(5,2)\thicklines
\put(0.1,0.4){\line(1,-1){1}}
\put(1.1,-0.6){\line(1,0){1}}
\put(2.1,-0.6){\line(1,1){1}}
\end{picture} }}
\def\coev{{\setlength{\unitlength}{0.5cm}
\begin{picture}(5,2)\thicklines
\put(3.1,4.1){\line(1,1){1}}
\put(4.1,5.1){\line(1,0){1}}
\put(5.1,5.1){\line(1,-1){1}}
\end{picture}}}

\def\tri{{ \setlength{\unitlength}{0.4cm}
\begin{picture}(5,2)\thicklines
\put(0.1,0.4){\line(0,1){1}}
\put(4.0,0.4){\line(0,1){1}}
\put(0.1,0.4){\line(1,0){3.9}}
\put(0.1,1.42){\line(1,0){3.9}}
\end{picture}}}
\begin{document}

{\renewcommand{\baselinestretch}{1}
  \title{\bf MAKING NON-TRIVIALLY
    ASSOCIATED MODULAR CATEGORIES FROM FINITE GROUPS\\ }
    \author{M. M. Al-Shomrani
     \&\ E. J. Beggs \\
\vspace{0.15in}\\
Department of Mathematics \\
University of Wales, Swansea\\
Singleton Park, Swansea SA2 8PP, UK\\}
\date{February \, 2003}

\maketitle {ABSTRACT :} We show that the double ${\cal D}$ of the
non-trivially associated tensor category constructed from left
coset representatives of a subgroup of a finite group $X$ is a
modular category.  Also we give a definition of the character of
an object in this category as an element of a braided Hopf algebra
in the category.  This definition is shown to be adjoint invariant
and multiplicative on tensor products.  A detailed example is
given. Finally we show an equivalence of categories between the
non-trivially associated double ${\cal D}$ and the trivially
associated category of
representations of the Drinfeld double of the group $D(X)$. }\\
\vspace{-1.5cm}
\section{Introduction} This paper will make continual use of
formulae and ideas from \cite{BNT}, and these definitions and
formulae will not be repeated, as they would add very considerably
to the length of the paper.  The paper \cite{BNT} is itself based
on the papers \cite{BGM, BM1}, but is mostly self contained in
terms of notation and definitions.  The book \cite{MajBook} has
been used as a standard reference for Hopf algebras, and
\cite{TWen, Bak} as references for modular categories.

In \cite{BNT} there is a construction of a non-trivially
associated tensor category ${\cal C}$ from data which is a choice
of left coset representatives $M$ for a subgroup $G$ of a finite
group $X$.  This introduces a binary operation $\cdot$ and a
$G$-valued `cocycle' $\tau$ on $M$.  There is also a double
construction where $X$ is viewed as a subgroup of a larger group.
This gives rise to a braided category ${\cal D}$, which is the
category of reps of an algebra $D$, which is itself in the
category, and it is the category that we concentrate on in this
paper.

It is our aim to show that the non-trivially associated algebra
$D$ has reps which have characters in the same way that the reps
of a finite group have characters, and also that the category of
its representations has a modular structure in the same way that
the category of reps of the double of a group has a modular
structure.

We begin by describing the indecomposable objects in ${\cal C}$,
in a similar manner to that used in \cite{BM1}.  A detailed
example is given using the group $D_{6}$.  Then we show how to
find the dual objects in the category, and again illustrate this
with the example.

Next we show that the rigid braided category ${\cal D}$ is a
ribbon category.  The ribbon maps are calculated for the
indecomposable objects in our example category.

In the next section we explicitly evaluate in ${\cal D}$ the standard diagram for
trace in a ribbon category \cite{MajBook}.
Then we define the character of an
object in ${\cal D}$ as an element of the dual of the braided Hopf
algebra $D$.  This element is shown to be right adjoint invariant.
Also we show that the character is multiplicative for the tensor
product of objects.  A formula is found for the character in ${\cal
D}$ in terms of characters of group representations.

The last ingredient needed for a modular category is the trace of
the double braiding, and this is calculated in ${\cal D}$ in terms
of group characters.  Then the matrices $S$, $T$ and $C$
implementing the modular representation are calculated explicitly
in our example.

Finally we show an equivalence of categories between the
non-trivially associated double ${\cal D}$ and the category of
representations of the Drinfeld double of the group $D(X)$.

Throughout the paper we assume that all groups mentioned are
finite, and that all vector spaces are finite dimensional.\!  We
take the base field to be the complex numbers $\Bbb C$.

\section{Indecomposable objects in ${\cal C}$}\label{indob}

The objects of ${\cal C}$ are the right representations of the
algebra $A$ described in \cite{BNT}.  We now look at the
indecomposable objects in ${\cal C}$, or the irreducible
representations of $A$, in a manner similar to that used in
\cite{BM1}.
\begin{theorem}  \hspace{0.2cm}
The indecomposable objects in ${\cal C}$ are of the form
$$
V=\bigoplus_{s\in\cal O} V_s
$$
where ${\cal O}$ is an orbit in $M$ under the $G$ action $\acl$,
and each $V_s$ is an irreducible right representation
of the stabilizer of $s$, ${\rm stab}(s)$.  Every object $T$ in ${\cal
C}$ can be written as a direct sum of indecomposable objects in ${\cal
C}$.
\end{theorem}
\textbf{Proof.}\hspace{0.5cm} For an object $T$ in ${\cal C}$ we
can use the $M$-grading to write
\begin{equation}
T=\bigoplus_{s\in M} T_s\, ,
\end{equation}
but as $M$ is a disjoint union of orbits ${{\cal O}_s}=\{s \acl u
: u \in G \,\}$ for $s\in M$, $T$ can be rewritten as a disjoint
sum over orbits,
\begin{equation}
T=\bigoplus_{\cal O} T_{\cal O} \,,
\end{equation}
where
\begin{equation}
 T_{\cal O}=\bigoplus_{s\in\cal O} T_s.
\end{equation}
Now we will define the stabilizer  of $s \in \cal O$,
 which is a subgroup of $G$, as
$$
{\rm stab}(s)=\{ u \in G:s \acl u =s\}.
$$
As $ \langle \eta \acbl u \rangle = \a \acl  u$ for all $\eta \in
T$, $T_s$ is a representation of the group ${\rm stab}(s)$.  Now
fix a base point $t \in {\cal O}$.  Because ${\rm stab}(t)$ is a
finite group, $T_t$ is a direct sum of irreducible group
representations $ {W_i}$ for $i=1,.. , m$, i.e.,
\begin{equation}
T_t={\bigoplus^{m}_{i=1}} W_{i}.
\end{equation}
Suppose that ${\cal O}=\{ t_{1}, t_{2},...., t_n \}$ where
$t_{1}=t$, and take $u_i \in G$ so that $t_i =t \acl u_i$.  Define
\begin{equation}
U_i={\bigoplus^{n}_{j=1}} W_{i}\acbl u_j  \subset
\bigoplus_{s\in\cal O} T_s.
\end{equation}
We claim that each $U_i$ is an indecomposable object in ${\cal
C}$.  For any $v \in G$ and $\xi \acbl u_{k} \in W_{i} \acbl
u_{k}$,
$$
(\xi \acbl u_{k}) \acbl v = \big(\xi \acbl( u_{k} v {u_{j}}^{-1}
)\big) \acbl u_{j},
$$
where $u_{k} v {u_{j}}^{-1} \in {\rm stab}(t)$ for some $u_{j} \in
G$. This shows that $U_{i}$ is a representation of $G$.  By the
definition of $U_{i}$, any subrepresentation of $U_{i}$ which
contains $W_{i}$ must be all of $U_{i}$.  Thus $U_i$ is an
indecomposable object in $\cal C$, and
\begin{equation}
 T_{\cal O}={\bigoplus^{m}_{i=1}} U_{i}. \qquad \square
\end{equation}

\begin{theorem} {\bf \{Schur's lemma\}} \hspace{0.2cm}
Let $V$  and $W$ be two indecomposable objects in ${\cal C}$, and
let \,$ \al:V  \longrightarrow W $\, be a morphism.  Then $\al$ is
zero or a scalar multiple of the identity.

\end{theorem}
\textbf{Proof.}\hspace{0.5cm} $V$ and $W$ are associated to orbits
${\cal O}$ and ${\cal O'}$ so that $V=\bigoplus_{s\in\cal O} V_s$
and $W=\bigoplus_{s\in\cal O^{'}} W_{s}$. As morphisms preserve
grade, if $\al \not= 0$, then ${\cal O}={\cal O^{'}}$. Now if we
take ${s\in\cal O} $, we find that $\al:V_s \longrightarrow W_s $
is a map of irreps of ${\rm stab}(s)$, so by Schur's lemma for
groups, any non-zero map is a scalar multiple of the identity, and
we have $ V_s =W_s $ as representations of ${\rm stab}(s)$.  Now
we need to check that the multiple of the identity is the same for
each ${s\in\cal O}$. Suppose $\al$ is a multiplication by
$\lambda$ on $V_s$.  Given ${t\in\cal O}$, there is a $u \in G$ so
that $ t \acl u =s$.  Then for $\eta \in V_t$,
$$
 \al (\eta)=\al ( \eta \acbl u )\acbl u^{-1}=
\lambda(\eta \acbl u)\acbl u^{-1}=\lambda \eta.\qquad\square
$$

\begin {Lem}\label{moving}{}\hspace{0.1cm}Let $V$ be an
indecomposable object in ${\cal C}$ associated to the orbit ${\cal
O}$.  Choose $s,t\in {\cal O}$ and $u \in G$ so that $s\acl u=t$.
Then $V_s$ and $V_t$ are irreps of ${\rm stab}(s)$ and ${\rm
stab}(t)$ respectively, and the group characters obey $
\chi_{_{V_{t}}}(v)=\chi_{_{V_{s}}}(u\,v \,u^{-1}) $.
\end{Lem}
\textbf{Proof.}\hspace{0.5cm} Note that $\acbl u$ is an invertible
map from $V_s$ to $V_t$.  Then we have the  commuting diagram

\vspace{.7cm}
$$\qquad\quad
\setlength{\unitlength}{0.5cm}
\begin{picture}(5,2)\thicklines
\put(-1,5){\hbox to 70pt{\rightarrowfill}}
\put(-1,1){\hbox
to70pt{\rightarrowfill}}
\put(-1.8,3){$\Bigg\downarrow$}
\put(4.2,3){$\Bigg\downarrow$}
\put(-1.9,5){\text{$V_s$}}
\put(4.1,5){\text{$V_s$}}
\put(-1.9,1){\text{$V_t$}}
\put(4.1,1){\text{$V_t$}}
\put(.3,5.6){\scriptsize{\text{$\acbl u
v {u}^{-1}$}}}
\put(1.2,1.5){\scriptsize{\text{$\acbl  v $}}}
\put(-1.4,3.25){\scriptsize{\text{$\acbl  u $}}}
\put(4.6,3.25){\scriptsize{\text{$\acbl  u $}}}
\end{picture}
\vspace{-.5cm}
$$
\noindent which implies that $ \text{trace} (\acbl u v
{u}^{-1}:V_s \to V_s) =\text{trace} (\acbl  v :V_t \to
V_t).\qquad\square $

\section{An example of indecomposable objects}\label{exind}
We give an example of indecomposable objects in the categories
discussed in the last section.  As we will later want to have a
category with braiding, we use the double construction in
\cite{BNT}. We also use lemma \ref{moving} to list the group
characters \cite{Grove} for every point in the orbit in terms of
the given base points.

 Take $X$ to be the dihedral group
$D_6=\langle a,b : a^6=b^2=e, ab=ba^5 \rangle$, whose elements we
list as $\{e, a, a^2, a^3, a^4, a^5, b, ba, ba^2, ba^3, ba^4, ba^5 \}$,
and $G$ to be the non-abelian normal subgroup of order 6 generated by
$a^2$ and $b$, i.e.\ $G=\{e, a^2, a^4, b, ba^2, ba^4\}$.  We choose
$M=\{e, a \}$.  The center of $D_6$ is the subgroup $\{e, a^3\}$, and
it has the following conjugacy classes: $\{e\}$, $\{a^3\}$, $\{a^2,
a^4\}$, $\{a, a^5\}$, $\{b, ba^2, ba^4\}$ and $\{ba, ba^3, ba^5\}$.

The category ${\cal{D}}$ consists of right
representations of the group $X=D_6$ which are
graded by $Y=D_6$ (as a set), using the actions $\actl:Y \times X
\rightarrow Y$ and $\actr: Y \times X \rightarrow X$ which are defined
as follows:
$$
y \actl x = x^{-1} y x, \quad \text{and}\quad vt \actr x =v^{-1} x
v^{'}= t x {t^{'}}^{-1}, \,\,
$$
 for $\,\,x \in X, \,y \in Y, \, v,v^{'}
 \in G \,\,$and$\,\,\, t,t^{'} \in M$  where $\,\,
 vt \actl x =v^{'}t^{'}$.

Now let  $V$ be an indecomposable object in ${\cal{D}}$.
   We  get the following cases:

   \noindent { \bf Case (1):} Take the orbit $\{e\}$ with base point
   $e$, whose stabilizer is the whole of $D_{6}$.  There are six
   possible irreducible group representations of the stabilizer, with
   their characters given by  table(1) \cite{TWood}:
$$\quad
\setlength{\unitlength}{0.5cm}
\begin{picture}(5,2)\thicklines
\put(.5,-6){\line(0,1){7.5}}
\put(2,-6){\line(0,1){6.}}
\put(3.5,-6){\line(0,1){7.5}}
\put(5.,-6){\line(0,1){7.5}}
\put(6.7,-6){\line(0,1){7.5}}
\put(11.2,-6){\line(0,1){7.5}}
\put(16.,-6){\line(0,1){7.5}}
\put(18.9,-6){\line(0,1){7.5}}
\put(21.5,-6){\line(0,1){7.5}}
\put(.5,-6){\line(1,0){21}}
\put(.5,-5){\line(1,0){21}}
\put(.5,-4){\line(1,0){21}}
\put(.5,-3){\line(1,0){21}}
\put(.5,-2){\line(1,0){21}}
\put(.5,-1){\line(1,0){21}}
\put(.5,0){\line(1,0){21}}
\put(.5,1.5){\line(1,0){21}}
\put(.8,.5){\text{irreps}}
\put(3.6,.5){\text{$\{e\}$}}
\put(5.05,.5){\text{$\{a^3\}$}}
\put(6.9,.5){\text{$\{b, ba^2,ba^4\}$}}
\put(11.3,.5){\text{$\{ba, ba^3,ba^5\}$}}
\put(16.1,.5){\text{$\{a^2,a^4\}$}}
\put(19,.5){\text{$\{a,a^5\}$}}
\put(.7,-.75){\text{$1_1$}}
\put(.7,-1.75){\text{$1_2$}}
\put(.7,-2.75){\text{$1_3$}}
\put(.7,-3.75){\text{$1_4$}}
\put(.7,-4.75){\text{$1_5$}}
\put(.7,-5.75){\text{$1_6$}}
\put(2.3,-.75){\text{$2_1$}}
\put(2.3,-1.75){\text{$2_2$}}
\put(2.3,-2.75){\text{$2_3$}}
\put(2.3,-3.75){\text{$2_4$}}
\put(2.3,-4.75){\text{$2_5$}}
\put(2.3,-5.75){\text{$2_6$}}
\put(4.,-.75){\text{1}}
\put(4.,-1.75){\text{1}}
\put(4.,-2.75){\text{1}}
\put(4.,-3.75){\text{1}}
\put(4.,-4.75){\text{2}}
\put(4.,-5.75){\text{2}}
\put(5.6,-.75){\text{1}}
\put(5.31,-1.75){\text{-1}}
\put(5.31,-2.75){\text{-1}}
\put(5.6,-3.75){\text{1}}
\put(5.3,-4.75){\text{-2}}
\put(5.6,-5.75){\text{2}}
\put(8.7,-.75){\text{1}}
\put(8.41,-1.75){\text{-1}}
\put(8.7,-2.75){\text{1}}
\put(8.41,-3.75){\text{-1}}
\put(8.7,-4.75){\text{0}}
\put(8.7,-5.75){\text{0}}
\put(13.5,-.75){\text{1}}
\put(13.5,-1.75){\text{1}}
\put(13.21,-2.75){\text{-1}}
\put(13.21,-3.75){\text{-1}}
\put(13.5,-4.75){\text{0}}
\put(13.5,-5.75){\text{0}}
\put(17.4,-.75){\text{1}}
\put(17.4,-1.75){\text{1}}
\put(17.4,-2.75){\text{1}}
\put(17.4,-3.75){\text{1}}
\put(17.11,-4.75){\text{-1}}
\put(17.11,-5.75){\text{-1}}
\put(19.9,-.75){\text{1}}
\put(19.62,-1.75){\text{-1}}
\put(19.62,-2.75){\text{-1}}
\put(19.9,-3.75){\text{1}}
\put(19.9,-4.75){\text{1}}
\put(19.62,-5.75){\text{-1}}
\put(10.6,-7){\text{table (1)}}

\end{picture}\qquad\qquad\qquad\qquad
\qquad\qquad\qquad\qquad\qquad\qquad
\qquad\quad
\vspace{.7cm}
$$
$$
$$
$$
$$

\noindent { \bf Case (2):} Take the orbit $\{a^3\}$ with base
point $a^3$, whose stabilizer is the whole of $D_{6}$.  There are
six possible irreps \{$2_1$, $2_2$, $2_3$, $2_4$, $2_5$, $2_6$\},
with characters given by table(1).

\noindent { \bf Case (3):} Take the orbit $\{a^2, a^4\}$ with base
point $a^2$, whose stabilizer is $\{e, a, a^2, a^3, a^4, a^5 \}$.
There are six irreps \{$3_0$, $3_1$, $3_2$, $3_3$, $3_4$, $3_5$\},
with characters  given by table(2), where $\omega=e^{{i\pi}/3}$.
Applying lemma \ref{moving} gives $
\chi_{_{V_{a^4}}}(v)=\chi_{_{V_{a^2}}}(bvb) $.
$$\quad
\setlength{\unitlength}{0.5cm}
\begin{picture}(5,2)\thicklines
\put(.5,-6){\line(0,1){7.5}}
\put(2,-6){\line(0,1){6.}}
\put(3.5,-6){\line(0,1){7.5}}
\put(6.5,-6){\line(0,1){7.5}}
\put(9.5,-6){\line(0,1){7.5}}
\put(12.5,-6){\line(0,1){7.5}}
\put(15.5,-6){\line(0,1){7.5}}
\put(18.5,-6){\line(0,1){7.5}}
\put(21.5,-6){\line(0,1){7.5}}
\put(.5,-6){\line(1,0){21}}
\put(.5,-5){\line(1,0){21}}
\put(.5,-4){\line(1,0){21}}
\put(.5,-3){\line(1,0){21}}
\put(.5,-2){\line(1,0){21}}
\put(.5,-1){\line(1,0){21}}
\put(.5,0){\line(1,0){21}}
\put(.5,1.5){\line(1,0){21}}
\put(.8,.5){\text{irreps}}
\put(4.6,.5){\text{$e$}}
\put(7.6,.5){\text{$a$}}
\put(10.6,.5){\text{$a^2$}}
\put(13.6,.5){\text{$a^3$}}
\put(16.6,.5){\text{$a^4$}}
\put(19.6,.5){\text{$a^5$}}
\put(.7,-.75){\text{$3_0$}}
\put(.7,-1.75){\text{$3_1$}}
\put(.7,-2.75){\text{$3_2$}}
\put(.7,-3.75){\text{$3_3$}}
\put(.7,-4.75){\text{$3_4$}}
\put(.7,-5.75){\text{$3_5$}}
\put(2.3,-.75){\text{$4_0$}}
\put(2.3,-1.75){\text{$4_1$}}
\put(2.3,-2.75){\text{$4_2$}}
\put(2.3,-3.75){\text{$4_3$}}
\put(2.3,-4.75){\text{$4_4$}}
\put(2.3,-5.75){\text{$4_5$}}
\put(4.7,-.75){\text{1}}
\put(4.7,-1.75){\text{1}}
\put(4.7,-2.75){\text{1}}
\put(4.7,-3.75){\text{1}}
\put(4.7,-4.75){\text{1}}
\put(4.7,-5.75){\text{1}}
\put(7.7,-.75){\text{1}}
\put(7.7,-1.75){\text{$\omega^1$}}
\put(7.7,-2.75){\text{$\omega^2$}}
\put(7.7,-3.75){\text{$\omega^3$}}
\put(7.7,-4.75){\text{$\omega^4$}}
\put(7.7,-5.75){\text{$\omega^5$}}
\put(10.9,-.75){\text{1}}
\put(10.9,-1.75){\text{$\omega^2$}}
\put(10.9,-2.75){\text{$\omega^4$}}
\put(10.9,-3.75){\text{1}}
\put(10.9,-4.75){\text{$\omega^2$}}
\put(10.9,-5.75){\text{$\omega^{4}$}}
\put(13.7,-.75){\text{1}}
\put(13.7,-1.75){\text{$\omega^3$}}
\put(13.7,-2.75){\text{1}}
\put(13.7,-3.75){\text{$\omega^3$}}
\put(13.7,-4.75){\text{1}}
\put(13.7,-5.75){\text{$\omega^{3}$}}
\put(16.7,-.75){\text{1}}
\put(16.7,-1.75){\text{$\omega^4$}}
\put(16.7,-2.75){\text{$\omega^2$}}
\put(16.7,-3.75){\text{1}}
\put(16.7,-4.75){\text{$\omega^{4}$}}
\put(16.7,-5.75){\text{$\omega^{2}$}}
\put(19.7,-.75){\text{1}}
\put(19.7,-1.75){\text{$\omega^5$}}
\put(19.7,-2.75){\text{$\omega^{4}$}}
\put(19.7,-3.75){\text{$\omega^{3}$}}
\put(19.7,-4.75){\text{$\omega^{2}$}}
\put(19.7,-5.75){\text{$\omega^{1}$}}
\put(10.6,-7){\text{table
(2)}}
\end{picture}\qquad\qquad\qquad\qquad\qquad\qquad
\qquad\qquad\qquad\qquad\qquad
\vspace{.7cm}
$$
$$
$$
$$
$$

\noindent { \bf Case (4):} Take the orbit $\{a, a^5\}$ with base
point $a$, whose stabilizer is $\{e, a, a^2, a^3, a^4, a^5 \}$.
There are six irreps \{$4_0$, $4_1$, $4_2$, $4_3$, $4_4$, $4_5$ \}
with characters given in table(2).  Applying lemma \ref{moving}
gives $ \chi_{_{V_{a^5}}}(v)=\chi_{_{V_a}}(ba^2vba^2) $.

 \noindent
{ \bf Case (5):} Take the orbit  $\{b, ba^2, ba^4\}$ with base
point $b$, whose stabilizer is $\{e, a^3, b, b a^3 \}$.  There are
four irreps with characters given by  table(3).  Applying lemma
\ref{moving} gives $ \chi_{_{V_{ba^2}}}(v)=\chi_{_{V_b}}(a^4va^2)$
 and
$\chi_{_{V_{ba^4}}}(v)=\chi_{_{V_b}}(a^2va^4)$.
$$
\setlength{\unitlength}{0.5cm}
\begin{picture}(5,2)\thicklines
\put(2,-4){\line(0,1){5.5}}
\put(3.9,-4){\line(0,1){5.5}}
\put(6.5,-4){\line(0,1){5.5}}
\put(9.5,-4){\line(0,1){5.5}}
\put(12.5,-4){\line(0,1){5.5}}
\put(15.5,-4){\line(0,1){5.5}}
\put(2,-4){\line(1,0){13.5}}
\put(2,-3){\line(1,0){13.5}}
\put(2,-2){\line(1,0){13.5}}
\put(2,-1){\line(1,0){13.5}}
\put(2,0){\line(1,0){13.5}}
\put(2,1.5){\line(1,0){13.5}}
\put(4.9,.5){\text{$e$}}
\put(7.6,.5){\text{$a^3$}}
\put(10.7,.5){\text{$b$}}
\put(13.6,.5){\text{$ba^3$}}
\put(2.3,-.75){\text{$5_{++}$}}
\put(2.3,-1.75){\text{$5_{+-}$}}
\put(2.3,-2.75){\text{$5_{-+}$}}
\put(2.3,-3.75){\text{$5_{--}$}}
\put(4.9,-.75){\text{1}}
\put(4.9,-1.75){\text{1}}
\put(4.9,-2.75){\text{1}}
\put(4.9,-3.75){\text{1}}
\put(7.7,-.75){\text{1}}
\put(7.7,-1.75){\text{1}}
\put(7.4,-2.75){\text{-1}}
\put(7.4,-3.75){\text{-1}}
\put(10.9,-.75){\text{1}}
\put(10.6,-1.75){\text{-1}}
\put(10.9,-2.75){\text{1}}
\put(10.6,-3.75){\text{-1}}
\put(13.7,-.75){\text{1}}
\put(13.4,-1.75){\text{-1}}
\put(13.4,-2.75){\text{-1}}
\put(13.7,-3.75){\text{1}}
\put(7.6,-5){\text{table (3)}}
\end{picture}\qquad\qquad\qquad
\qquad\qquad\qquad\qquad
\vspace{.9cm}
$$
$$
$$
\noindent { \bf Case (6):} Take the orbit  $\{ba, ba^3, ba^5\}$
with base point $ba$, whose stabilizer is $\{e, a^3, ba , b
a^4\}$. There are four irreps with character given by table(4).
Applying lemma \ref{moving} gives
$\chi_{_{V_{ba^3}}}(v)=\chi_{_{V_{ba}}}(a^4va^2)$ and
$\chi_{_{V_{ba^5}}}(v)=\chi_{_{V_{ba}}}(a^2va^4)$.
$$
\setlength{\unitlength}{0.5cm}
\begin{picture}(5,2)\thicklines
\put(2,-4){\line(0,1){5.5}}
\put(3.9,-4){\line(0,1){5.5}}
\put(6.5,-4){\line(0,1){5.5}}
\put(9.5,-4){\line(0,1){5.5}}
\put(12.5,-4){\line(0,1){5.5}}
\put(15.5,-4){\line(0,1){5.5}}
\put(2,-4){\line(1,0){13.5}}
\put(2,-3){\line(1,0){13.5}}
\put(2,-2){\line(1,0){13.5}}
\put(2,-1){\line(1,0){13.5}}
\put(2,0){\line(1,0){13.5}}
\put(2,1.5){\line(1,0){13.5}}
\put(4.9,.5){\text{$e$}}
\put(7.6,.5){\text{$a^3$}}
\put(10.6,.5){\text{$ba$}}
\put(13.6,.5){\text{$ba^4$}}
\put(2.3,-.75){\text{$6_{++}$}}
\put(2.3,-1.75){\text{$6_{-+}$}}
\put(2.3,-2.75){\text{$6_{+-}$}}
\put(2.3,-3.75){\text{$6_{--}$}}
\put(4.9,-.75){\text{1}}
\put(4.9,-1.75){\text{1}}
\put(4.9,-2.75){\text{1}}
\put(4.9,-3.75){\text{1}}
\put(7.7,-.75){\text{1}}
\put(7.4,-1.75){\text{-1}}
\put(7.7,-2.75){\text{1}}
\put(7.4,-3.75){\text{-1}}
\put(10.9,-.75){\text{1}}
\put(10.9,-1.75){\text{1}}
\put(10.6,-2.75){\text{-1}}
\put(10.6,-3.75){\text{-1}}
\put(13.7,-.75){\text{1}}
\put(13.4,-1.75){\text{-1}}
\put(13.4,-2.75){\text{-1}}
\put(13.7,-3.75){\text{1}}
\put(7.6,-5){\text{table (4)}}
\end{picture}\qquad\qquad\qquad\qquad\qquad\qquad
\qquad
$$
\vspace{.5cm}
\section {Duals of indecomposable objects in ${\cal C}$}
\label{dualind} Given an irreducible object $V$ with associated
orbit ${\cal O}$ in  ${\cal C}$, how do we find its dual $V^{*}$?
The dual would be described as in  section \ref{indob}, by an
orbit, a base point in the orbit and a right group representation
of the stabilizer of the base point.  Using the formula $ (s^{L}
\cdot s) \acl u = \big( s^{L} \acl (s \acr u) \big)\cdot (s \acl
u)=e $, we see that the left inverse of a point in the orbit
containing $s$ is in the orbit containing $s^{L}$.  By using the
evaluation map from $V^{*}\otimes V$ to the field we can take
$(V^{*})_{s^{L}}=(V_s)^{*}$ as vector spaces.  We use $\check\acl$
as the action of ${\rm stab}(s)$ on $(V_s)^{*}$, i.e. $ (\alpha
\check\acl z) (\xi \acbl z)= \alpha (\xi) $ for $\alpha \in
(V_s)^{*}$ and $\xi \in V_s$.  The action $ \acbl$ of ${\rm
stab}(s^{L}) $ on $(V^{*})_{s^{L}}$ is given by $ \alpha \acbl (s
\acr z)= \alpha \check \acl z$ for $z\in {\rm stab}(s)$.  In terms
of group characters this gives
$$
\chi_{_{(V^{*})_{s^{L}}}}(s\acr z)\,=\,\chi_{_{(V_{s})^{*}}}(z)\,
,  \quad z\in {\rm stab}(s)\ .
$$
If we take ${\cal O}^{L}\,=\,\{s^{L}:s\in{\cal O}\}$ to have base
point $p$, and choose  $u\in G$ so that $p\acl u=s^{L}$, then
using lemma \ref{moving}  gives
\begin{eqnarray}\label{finddual}
    \chi_{_{(V^{*})_{s^{L}}}}(s\acr
z)\,=\,\chi_{_{(V_{s})^{*}}}(z) \,=\,\chi_{_{(V^{*})_{p}}}(u(s\acr
z)u^{-1})\, ,\quad z\in {\rm stab}(s)\ .
\end{eqnarray}
This formula allows us to find the character of $V^{*}$ at its
base point $p$ as a representation of ${\rm stab}(p)$ in terms of the
character of the dual of $V_{s}$ as a representation of ${\rm
stab}(s)$.

\begin{Lem}\label{tendual} In ${\cal C}$ we can regard the dual
$(V\otimes W)^{*}$ as $W^{*}\otimes V^{*}$ with the evaluation
$$
(\al \otimes \beta )(\ab)\,=\,
\big( \al \acbl \tau (\langle\beta\rangle ,\b \cdot \a) \big) (\eta) \
\big(\beta \acbl \tau (\b , \a)^{-1} \big)(\xi).
$$
Given a basis $\{\xi\}$ of $V$ and a basis $\{\eta\}$ of $W$, the dual basis
$\{\widehat{\xi\otimes\eta}\}$
of $W^{*}\otimes V^{*}$ can be written in terms of the dual basis of $V^{*}$
and $W^{*}$ as
$$
\widehat{\xi\otimes\eta}\,=\, \hat\eta
\acbl \tau (\b^{L}\acl \tau (\b , \a) ,\b \cdot \a)^{-1}
\,\otimes\,\hat\xi \acbl \tau (\b ,
\a)\ .
$$
\end{Lem}
\textbf{Proof.} \,\,\,
Applying the associator
to $(\al \otimes \beta )\otimes(\ab)$ gives
$$
\al \acbl \tau (\langle\beta\rangle ,\b \cdot \a) \otimes
\big(\beta \otimes(\ab)\big)\ ,
$$
and then applying the inverse associator gives
$$
\al \acbl \tau (\langle\beta\rangle ,\b \cdot \a) \otimes \Big(\big(\beta
\acbl \tau (\b , \a)^{-1}\otimes \xi\big) \otimes \eta \Big)\ .
$$
Applying the evaluation map first to $\beta \acbl \tau (\b ,
\a)^{-1}\otimes \xi$, then to $\al \acbl \tau (\langle\beta\rangle
,\b \cdot   \a) \otimes \eta $ gives the first equation. For the
evaluation to be non-zero, we need $\big(\langle \beta \rangle
\acl \tau (\b , \a)^{-1} \big)\cdot\b =e$ which implies $\langle
\beta \rangle \acl \tau (\b , \a)^{-1} =\b^{L}$, or equivalently
$\langle \beta \rangle =\b^{L}\acl \tau (\b , \a)$.  This gives
the second equation.\quad$\square$

\begin{example}
    Using (\ref{finddual}) we calculate the duals of the objects given in
    the last section.

    \noindent { \bf Case (1):} The orbit $\{e\}$ has left inverse $\{e\}$,
    so $\chi_{(V^{*})_{e}}=\chi_{(V_{e})^{*}}$.  By a calculation
    with  group characters, all the listed irreps of ${\rm
    stab}(e)$ are self-dual, so $1_{r}^{*}=1_{r}$ for $r\in\{1,\dots,6\}$.

    \noindent { \bf Case (2):} The orbit $\{a^{3}\}$ has left inverse
    $\{a^{3}\}$, so $\chi_{(V^{*})_{a^{3}}}=\chi_{(V_{a^{3}})^{*}}$.  As
    in the last case the group representations are self-dual, so
    $2_{r}^{*}=2_{r}$ for $r\in\{1,\dots,6\}$.

    \noindent { \bf Case (3):} The left inverse of the base point $a^{2}$
    is $a^{4}$, which is still in the orbit.  As {\it group} representations,
    the dual of $3_{r}$ is $3_{6-r}$ (mod 6).  Applying lemma \ref{moving} to
    move the base point, we see that the dual of $3_{r}$ {\it in the
    category} is $3_{r}$.

    \noindent { \bf Case (4):} The left inverse of the base point $a$
    is $a^{5}$, which is still in the orbit.  As in the last case,
    the dual of $4_{r}$ {\it in the
    category} is $4_{r}$.

    \noindent { \bf Case (5):} The left inverse of the base point is
    itself, and as group representations, all  case 5 irreps
    are self dual.  We deduce that in the category the objects are self dual.

    \noindent { \bf Case (6):} Self dual, as in case 5.

\end{example}

\section {The ribbon map on the category ${\cal D}$ }
\begin{theorem}   \hspace{0.1cm}
The ribbon transformation  $ \theta_V : V \longrightarrow
V $  for any object $ V $ in
${\cal D}$ can be defined by
$\theta_V(\xi)=\xi \,  \hat{\triangleleft} \, \|\xi\|$.
\end{theorem}
\textbf{Proof.}\hspace{0.2cm} In the following  lemmas we
 show that the required properties hold.\quad$\square$

\begin{Lem}{}\,\,\,
 $\theta_V$ is a morphism in the category.
\end{Lem}
\textbf{Proof.}\,\,\,Begin by checking the $X$-grade, for $\xi \in V$
$$
\|\theta_V (\xi) \|=\big\|\xi \, \, \hat{\triangleleft} \, \|\xi\| \big\|
=\bbb \actl \bbb =\bbb.
$$
Now we check the $X$-action, i.e. that \, $\theta_V(\xi \achl
\,x)= \theta_V(\xi) \achl \,x$.
\begin{equation*}
\begin{split}
\theta_V(\xi \achl \,x)&=(\xi \achl \,x) \achl \|\xi \achl \,x \|
=(\xi \achl \,x) \achl (\bbb \actl \,x) \\&=\xi \achl \,x x^{-1 }\bbb
x=(\xi \achl \bbb)\achl\,x=\theta_V(\xi)\achl \,x.  \qquad\square
\end{split}
\end{equation*}

\begin{Lem}{}\,\,\, \,\,\,For any two objects  $V$  and  $W$
in ${\cal{D}}$,
\begin{equation*}
\theta_{V\otimes W} =\Psi^{-1}_{V\otimes W }\,\circ\,
\Psi^{-1}_{W\otimes V }\,\circ \,(\,\theta_V \otimes  \theta_W
\,)=(\,\theta_V \otimes  \theta_W \,)\,\circ\,\Psi^{-1}_{V\otimes
W }\,\circ\, \Psi^{-1}_{W\otimes V }\,
\end{equation*}
This can also be described by figure 1:
$$
\put(3,17){\line(0,1){50}}
\put(20,17){\line(0,1){50}}
\put(3,-60){\line(0,1){50}}
\put(20,-60){\line(0,1){50}}
\put(1,1){\vw} \,\,\,\,\qquad \,\,\,\,=\put(12,-24){\brai}
\put(34,31){\t}
\put(65,30){\w}
\put(10.5,-64.5){\brai}
\put(44,46){\line(0,1){20}}
\put(70,44){\line(0,1){22}}
\put(97,1){\text {=}}
\put(108.5,-21){\brai}
\put(110,20){\brai}
\put(129,-31){\t}
\put(160,-31){\w}
\put(140,-59){\line(0,1){18}}
\put(166,-60){\line(0,1){20}}
\put(-57,70){\text {V}}
\put(40,70){\text {V}}
\put(140,70){\text {V}}
\put(-40,70){\text {W}}
\put(64,70){\text {W}}
\put(160,70){\text {W}}
\put(-57,-71){\text {V}}
\put(40,-71){\text {V}}
\put(136,-71){\text {V}}
\put(-40,-71){\text {W}}
\put(61,-71){\text {W}}
\put(160,-71){\text {W}}
\,\,\,\,\,\,\,\,\,\,\,\,\,\,\,\,\,\,\,\,\,\,\,\,\,\,\,\,\,\,\,
\,\,\,\,\,\,\,\qquad\,\,\,\,\,\,\,\qquad\,\,\,\,\,\,\qquad
$$
$$
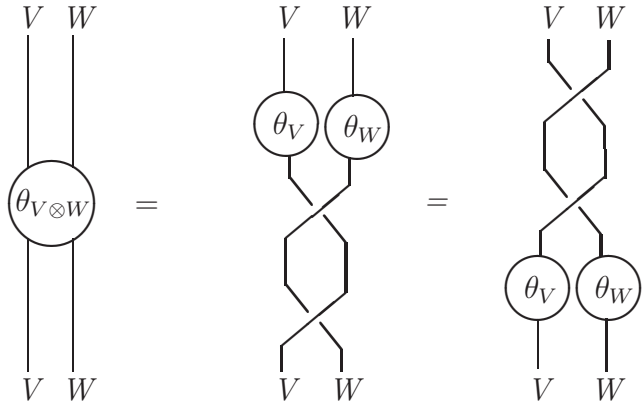
 \centerline{ \rm figure 1}
$$
\end{Lem}
\textbf{Proof.}\,\,\,\,\,
First calculate $ \Psi \, \bigl(\Psi(\ab)\bigr)$ for $ \xi\in V$ and
$\eta\in W$, beginning with
\begin {equation}\label{rrtt} \Psi \,
\bigl(\Psi(\ab)\bigr)\,=\,\Psi\bigl(\eta\,\achl(\b\acl\aa)^{-1}\,\otimes
\xi\achl\aa\,\bigr) .
\end{equation}
To simplify what follows we shall use the substitutions
\begin{equation}
\begin{split}
\eta^ {'}=\xi\achl\aa\,\,\,\,\,\,\,\,\text{and}\,\,\,\,\,\xi^{'}
=\eta\achl(\,\b\acl\aa)^{-1},
\end{split}
\end{equation}
so equation (\ref{rrtt}) can be rewritten as
\begin{equation}
\begin{split}
\Psi \, \bigl(\Psi(\ab)\bigr)\,&=\,\Psi(\xi^{'}\,\otimes\,\eta^ {'})\\
&=\eta^{'}\achl(\langle\xi^{'}\rangle\,\acl |\eta^{'}| )^{-1}\,
\otimes\,\xi^{'}\achl |\eta^{'}|.
\end{split}
\end{equation}
As \,$\eta^ {'}=\xi\,\achl \,\aa\,=\,\xi\,\bar\acl\,\aa $, \,then
\,$|\eta^{'}|=\,\big|\,\xi\,\bar\acl \,\aa\,\big| =\, (\,\b
\triangleright\aa\,)^{-1}\bb\aa$,\,\,so
\begin{equation}\label{gghh}
\begin{split}
\xi^ {'}\achl |\eta^{'}|&=\eta \,\achl(\,\b\acl\aa)^{-1}\,\,
(\,\b \triangleright\aa\,)^{-1}\bb\aa\\
&=\eta \, \achl\big((\,\b \triangleright\aa\,)(\,\b\acl\aa)\big)^{-1}
\bb\aa \\
&=\eta\,\achl\aai\bi\bb\aa.
\end{split}
\end{equation}
Hence if we put $ y=\|\ab\|=\|\xi\| \circ
\|\eta\|=\aai\, \bbi \b\a  $,
\begin{equation}\label{gg3}
    \begin{split} \Psi \bigl( \Psi(\xi\otimes\eta)
\bigr) \hat{\triangleleft}\,\|\,\xi\otimes\eta\,\|= \xi
\,\achl\aa(\langle\xi^{'}\rangle\,\acl |\eta^{'}| )^{-1} (p\,\tilde
\acr \|\ab\|)\otimes\eta\,\achl\aai\a,
\end{split}
\end{equation}
where, using (\ref{gghh}),
\begin{equation*}
\begin{split}
p&=\| \xi'\,\achl|\eta^{'}| \,\|\,=\,
{\big|\xi^{'}\,
\acbl |\eta^{'}|\big|}^{-1} \langle\xi^{'}\,\acbl |\eta^{'}| \rangle
=\|\eta\|\,\actl\|\eta\|\,y^{-1}=\|\eta\|\,\actl \, y^{-1} \\
 p\,\actr\|\ab\|
&=(\|\eta\|\,\actl \, y^{-1} )\, \actr \, y
=(\|\eta\| \, \actr \, y^{-1} )^{-1}, \,\,\,\,\,\,\,
\end{split}
\end{equation*}
As $\big\|\xi^{'}\,
\acbl |\eta^{'}| \big\|=v^{'} t^{'} =
\|\eta\| \, \actl \, y^{-1}$,
 by unique factorization, $t^{'}=\langle \xi^{'}\rangle
\acl|\eta^{'}|$.  Then $\|\eta\| \, \actr \, y^{-1}= \a\,y^{-1}
{t^{'}}^{-1}$, which implies that
\begin{equation}
\aa(\langle\xi^{'}\rangle\,
\acl |\eta^{'}| )^{-1}(\|\eta\|\actr y^{-1} )^{-1}\,=\,
\aa{t^{'}}^{-1}{t^{'}}y \ai\,=\,\bbb\ .
\end{equation}
Substituting this in (\ref{gg3}) gives
$$
\Psi \bigl( \Psi(\xi\otimes\eta) \bigr) \,
\hat{\triangleleft}\,\|\,\xi\otimes\eta\,\|=\xi \,
\achl \bbb \, \otimes\,\eta \,\achl\,\aaa .
\qquad\square
$$

\begin {Lem}{}
For the unit object\, $ \underline{\bf 1}=\Bbb C $ in
 ${\cal {D}}$,
$ \theta_{\underline{\bf 1}}$ is the identity.
\end{Lem}
\textbf{Proof.}\hspace{0.2cm}\,\,For any object $V$ in
 ${\cal {D}}$, $ \theta_V : V \longrightarrow V $
\,is defined by
$$
\theta_V(\xi)=\xi \, \, \hat{\triangleleft} \, \|\xi\| \qquad  \text{for}  \,\, \,\xi
\in V.
$$
If we choose  $V=\underline{\bf 1}=\Bbb C$ then $\theta_{\underline{\bf 1}}(\xi)
=\xi \, \hat{\triangleleft}\,e=\xi$ as $\bbb=e$.\quad$\square$

\begin {Lem}{}\hspace{0.1cm}\,\,\,For any object $V$
in ${\cal D}$, $ (\theta_V )^*=\theta_{V^*} $\,(see figure 2).
$$
$$
$$
$$
$$
\setlength{\unitlength}{0.5cm}
\begin{picture}(5,2)\thicklines
\put(0.1,0.4){\line(0,1){4.9}}
\put(0.1,0.4){\line(1,-1){1}}
\put(1.1,-0.6){\line(1,0){1}}
\put(2.1,-0.6){\line(1,1){1}}
\put(3.1,0.4){\line(0,1){1}}
\put(3.1,3.1){\line(0,1){1}}
\put(2.4,2.1){\t}
\put(3.1,4.1){\line(1,1){1}}
\put(4.1,5.1){\line(1,0){1}}
\put(5.1,5.1){\line(1,-1){1}}
\put(6.1,-0.8){\line(0,1){4.9}}
\put(5.8,-1.6){\text {$V^*$}}
\put(-0.1,5.5){\text {$V^*$}}
\put(8,1.4){\text {$=$}}
\put(11.1,-0.8){\line(0,1){2.2}}
\put(10.4,2.1){\tt}
\put(11.1,3.1){\line(0,1){2.2}}
\put(10.8,-1.6){\text {$V^*$}}
\put(10.8,5.5){\text {$V^*$}}
\end{picture}
\,\,\,\,\,\,\,\,\,\,\,\,\,\,\,\,\,\,\,\,\,\,\,\,\,\,\,\,
\,\,\,\,\,\,\,\,\,\,\qquad \qquad
$$
$$
\,\, \centerline{\rm {\quad\qquad figure 2}}
$$
\end{Lem}
\vspace{-.7cm}
\textbf{Proof.}\hspace{0.5cm}
Begin with
$$
{\rm coev}_V(1)=\sum_{\xi\in \,\text{basis of} \,\,V} \xi\, \achl\,\tilde\tau(
\,\bbb^L ,\bbb \,)^{-1} \otimes \hat\xi=\sum_{\xi\in \,\text{basis of} \,\, V} \xi\, \achl\,\tau(
\,\b^L ,\b \,)^{-1} \otimes \hat\xi
$$
For $\alpha\in V^*$, we follow figure 2  and calculate
\begin{equation}\label{gg4}
\begin{split}
\big(\theta_V \big)^*(\alpha)=({\rm eval}_V\otimes id\,)\sum_{\xi\in \,\text{basis of} \,\, V}
\Phi^{-1}\Big( \al \otimes \Big(\theta_V \big( \xi \achl\tau( \b^L ,\b
)^{-1}\big) \otimes \hat\xi\Big)\Big).
\end{split}
\end{equation}
Now as $\tau(\b^{L},\b) =\b^{L}\b $,
\begin{equation*}
\begin{split}
\big\| \,\xi\, \achl\,\tau(\,\b^L ,\b \,)^{-1} \,\big\|&=\bbb\,\actl
\,(\b^{L}\b)^{-1}\\
&=\b^{L} \b \,\bbi \,\b\bi\b^{L-1}\,\\
&=\b^{L} \b \,\bbi \,\b^{L-1},\,
\\
\theta_V  \,\big( \xi\, \achl\,\tau(\,\b^L ,\b \,)^{-1}\big)
&= \big( \xi\, \achl\,\tau(\,\b^L ,\b \,)^{-1}\big) \, \achl \,
 \big\| \,\xi\, \achl\,\tilde\tau(\,\bbb^L ,\bbb \,)^{-1} \,\big\| \\
&=\xi\, \achl\,\b^{-1} \b^{L-1}\b^{L} \b \,\bbi \,\b^{L-1}\, \\
&=\xi\, \achl\,\bbi \,\b^{L-1}.\,
\end{split}
\end{equation*}
The next step is to find
\begin{equation*}
\begin{split}
\Phi^{-1}\Big( \al \otimes \Big(\big(\xi \achl \,\bbi \b^{L-1} \big)
\otimes \hat\xi\Big)\Big) &= \Big(\al \, \achl \,\tilde{\tau}\big( \|
\,\xi \, \achl \,\bbi \b^{L-1}\|,\|\hat\xi\| \,\big)^{-1} \\& \qquad
\,\,\, \otimes \big(\xi\, \achl \,\bbi \b^{L-1} \big)
\Big)\otimes\hat\xi.
\end{split}
\end{equation*}
As \,\,
\begin{equation*}
\begin{split}
\| \,\xi \, \achl \,\bbi \b^{L-1}\|&=\bbb \,\actl \, \bbi \, \b^{L-1} \\
&=\b^{L} \,\bb \,\bbi \,\b \,\bbi \,\b^{L-1} \\
&=\tau(\b^{L},\b) \, \bbi \,\b^{L-1}\, \\
&=\tau(\b^{L},\b) \, \bbi \,\b \,\tau(\b^{L},\b)^{-1} \\
&=\tau(\b^{L},\b) \, \bbi \,\big( \,\b \,\acr\,\tau(\b^{L},\b)^{-1} \,\big)
  \big( \,\b \,\acl\,\tau(\b^{L},\b)^{-1} \,\big) \,,\,\,\,
\,\,\,\,\,\,\,\,\,\,\,\,\,\,\,\,\,\,\,\,\,\,\,\,\,\,
\end{split}
\end{equation*}
then as $\|\hat\xi\|={\|\xi\|}^{L}=\bb\tau(\b^{L},\b)^{-1}\b^{L}$,
\begin{equation*}
\begin{split}
\Phi^{-1}\Big( \al \otimes \Big(\big(\xi \achl \,\bbi \b^{L-1}
 \big) \otimes \hat\xi\Big)\Big)&=\Big(\al \,\achl \, \tau\,\big(
 \langle \xi\rangle \acl\,\tau(\b^{L},\b)^{-1}, \b^{L} \,\big)^{-1}\\&
 \qquad \,\,\,\otimes \big(\xi \,\achl \,\bbi \,\b^{L-1}\, \big)
 \Big)\otimes\hat\xi.
\end{split}
\end{equation*}
Put $v=\tau(\b^{L},\b)^{-1}=\b^{-1} \b^{L-1}$
and\,$w=\tau(\b\acl v,\b^{L})^{-1} =
\big((\b\acl v)\b^{L} \big)^{-1}$,\,
then substituting in (\ref{gg4}) gives
\begin{equation}\label{gg5}
\begin{split}
\big(\theta_V \big)^*(\alpha)=({\rm eval}_V\otimes {\rm id})\sum_{\xi\in \,\text{basis of} \,\, V}
\Big( (\al \, \achl\,w)\otimes (\xi\, \achl\,\bbi \,\b^{L-1}\,
 )\Big) \otimes \hat\xi.
\end{split}
\end{equation}
For a given term in the sum to be non-zero, we require
\begin{equation}
\|\al \|=\|\hat \xi \|=\bbb^{L}=\bb \bi\ ,
\end{equation}
and we proceed under this assumption. Now calculate
\begin{equation}
\begin{split}
{\rm eval}_V\big( (\al \, \achl\,w)\otimes (\xi\, \achl\,\bbi \,\b^{L-1}\,
 )\big)\,=\,\big(\beta \,\achl \,(\,\bbb \,\actr \, p\,)\big)(\xi \,
 \actl \, p) \,=\,\beta \,( \,\xi \,)
\end{split}
\end{equation}
where  \,\,$p=\bbi \,\b^{L-1}\, \, \text{and} \,\,
\beta=\al \,\achl \, w \,(\,\bbb  \,\actr \, p\,)^{-1}$.
Next we want to find $\bbb \,\actr \, p$.  To do this, we first find
\begin{equation}
\begin{split}
\,\bbb  \,\actl \, p\,&=\b^{L}\, \bb \, \bbi \, \b \,\bbi \,\b^{L-1}\, \\
&=v^{-1} \,\bbi \, \b \,v
=v^{-1} \,\bbi \,  (\b \,\acr \,v \,)(\b \,\acl \,v \,),
\end{split}
\end{equation}
and hence
\begin{equation}
\begin{split}
\,\bbb  \,\actr \, p\,&=\b\,p\,(\b \,\acl \,v \,)^{-1}\\
&=\b \, \bbi\, \b \, v \,(\b \,\acl \,v \,)^{-1} \\
&= \b \, \bbi\,  \,(\b \,\acr \,v \,).
\end{split}
\end{equation}
Thus
\begin{equation}
\begin{split}
\beta &=\al \,\achl \, w \,(\b \,\acr \,v \,)^{-1} \, \bb\, \bi \\
&=\al \,\achl \,\b^{L-1} \,(\b \,\acl \,v \,)^{-1} \,(\b \,\acr \,v \,)^{-1}
\,\bb\, \bi \\
&= \al \,\achl \,\b \,v  \,(\b  \,v \,)^{-1}\,\bb\, \bi
=\al \,\achl \,\bb\, \bi\ .
\end{split}
\end{equation}
Now substituting these last equations in (\ref{gg5}) gives
\begin{equation}
\begin{split}
\big(\theta_V \big)^*(\alpha)=\sum_{ \begin{subarray}\,\xi\in \,\text{basis of} \,\, V \,\,
\text{with} \,\bb \bi=\| \al\, \| \end{subarray}}
\big( \al \, \achl\,\| \al\, \|\, \big)(\xi)  \,. \,\, \hat\xi
\end{split}
\end{equation}
Take a basis $\xi_1  \,,\xi_2 \,, .\,.\,.\,.\,,\xi_n$ with $\big(
\al \, \achl\,\| \al\, \|\, \big)(\xi_i)$ being 1 if $i=1$, and 0
otherwise.  Then
$$
\big(\theta_V \big)^*(\alpha)= \,\hat{\xi}_1\,\,+ \,0=\al \, \achl\,\| \al\, \|
=\theta_{V^*}(\al)\ ,
$$
where $\hat{\xi}_1\,,\,
\hat{\xi}_2\,\,,\,\,.\,\,.\,\,.\,\,,\hat{\xi}_n$ is the dual basis
of $V^*$ defined by $\hat{\xi}_i\,(\,\xi_j\,)=\delta_{i,j}$ \,.
$\qquad\square$

\begin{example}
 We return to the example of section \ref{exind}.
First we calculate the value of the ribbon map on the
indecomposable objects. For an irreducible representation $V$, we
have \,$\theta_V:V \to V$\,defined by
$
\theta_V(\xi)=\xi \, \, \hat{\triangleleft} \, \|\xi\| \,\,
\text{for}   \,\xi \in V.
$
  At the base point ${s\in\cal O}$, we have\,\,$\theta_V(\xi)=\xi \,
\, \bar{\triangleleft} \, s \,\,  \text{for}   \,\xi \in V$\,
and\,$\theta:V_s \to V_s$\, is a multiple $\Theta_V$, say, of the
identity or,  in more explicitly, \rm{trace}\,$(\theta:V_s \to
V_s)= \Theta_V\, \text{\rm{dim}}_{\Bbb{C}}\,(V_s)$,  i.e.,
\begin{equation}
\Theta_V=\frac{\text{\rm {group character}}\,(s)}{\text{\rm
{dim}}_{\Bbb{C}}\,(V_s)}
\end{equation}
And then, for the different cases we will get the following table
:
$$
$$
$$
$$
$$
\vspace{.5cm}
$$\qquad\qquad\qquad\qquad
\setlength{\unitlength}{0.5cm}
\begin{picture}(5,2)\thicklines
\put(2,-8.){\line(0,1){17.52}} \put(4.3,-8){\line(0,1){17.5}}
\put(8,-8){\line(0,1){17.5}} \put(10.3,-8){\line(0,1){17.5}}
\put(14.5,-8){\line(0,1){17.52}}

\put(2,2){\line(1,0){12.5}}
\put(2,3){\line(1,0){12.5}}
\put(2,4){\line(1,0){12.5}}
\put(2,5){\line(1,0){12.5}}
\put(2,6){\line(1,0){12.5}}
\put(2,7){\line(1,0){12.5}}
\put(2,8){\line(1,0){12.5}}
\put(2,9.5){\line(1,0){12.5}}
\put(2,1){\line(1,0){12.5}}
\put(2,3){\line(1,0){12.5}}
\put(2,0){\line(1,0){12.5}}
\put(2,-1){\line(1,0){12.5}}
\put(2,-2){\line(1,0){12.5}}
\put(2,-3){\line(1,0){12.5}}
\put(2,-4){\line(1,0){12.5}}
\put(2,-5){\line(1,0){12.5}}
\put(2,-6){\line(1,0){12.5}}
\put(2,-7){\line(1,0){12.5}}
\put(2,-8){\line(1,0){12.5}}

\put(11.8,8.5){\text{$\Theta_V$}} \put(5.6,8.5){\text{$\Theta_V$}}
\put(2.1,8.5){\text{irreps}} \put(8.1,8.5){\text{irreps}}

\put(2.6,7.25){\text{$1_1$}} \put(2.6,6.25){\text{$1_2$}}
\put(2.6,5.25){\text{$1_3$}} \put(2.6,4.25){\text{$1_4$}}
\put(2.6,3.25){\text{$1_5$}} \put(2.6,2.25){\text{$1_6$}}
\put(2.6,1.25){\text{$2_1$}} \put(2.6,.25){\text{$2_2$}}
\put(2.6,-.75){\text{$2_3$}} \put(2.6,-1.75){\text{$2_4$}}
\put(2.6,-2.75){\text{$2_5$}} \put(2.6,-3.75){\text{$2_6$}}
\put(2.6,-4.75){\text{$3_0$}} \put(2.6,-5.75){\text{$3_1$}}
\put(2.6,-6.75){\text{$3_2$}} \put(2.6,-7.75){\text{$3_3$}}
\put(8.7,7.25){\text{$3_4$}}

\put(8.7,6.25){\text{$3_5$}} \put(8.7,5.25){\text{$4_0$}}
\put(8.7,4.25){\text{$4_1$}} \put(8.7,3.25){\text{$4_2$}}
\put(8.7,2.25){\text{$4_3$}} \put(8.7,1.25){\text{$4_4$}}
\put(8.7,.25){\text{$4_5$}} \put(8.7,-.75){\text{$5_{++}$}}
\put(8.7,-1.75){\text{$5_{+-}$}} \put(8.7,-2.75){\text{$5_{-+}$}}
\put(8.7,-3.75){\text{$5_{--}$}} \put(8.7,-4.75){\text{$6_{++}$}}
\put(8.7,-5.75){\text{$6_{-+}$}} \put(8.7,-6.75){\text{$6_{+-}$}}
\put(8.7,-7.75){\text{$6_{--}$}}

\put(5.7,7.25){1}

\put(5.7,6.25){1}
\put(5.7,5.25){1}
\put(5.7,4.25){1}
\put(5.7,3.25){1}
\put(5.7,2.25){1}
\put(5.7,1.25){1}
\put(5.5,0.25){-1}
\put(5.5,-.75){-1}
\put(5.7,-1.75){1}
\put(5.5,-2.75){-1}
\put(5.7,-3.75){1}
\put(5.7,-4.75){1}
\put(5.6,-5.75){$\omega^{2}$}
\put(5.6,-6.75){$\omega^{4}$}
\put(5.7,-7.75){1}

\put(12.,7.25){$\omega^{2}$}
\put(12.,6.25){$\omega^{4}$}
\put(12.1,5.25){1}
\put(12.,4.25){$\omega^{1}$}
\put(12.,3.25){$\omega^{2}$}
\put(11.8,2.25){-1}
\put(12.,1.25){$\omega^{4}$}
\put(12.,0.25){$\omega^{5}$}
\put(12.1,-.75){1}
\put(11.8,-1.75){-1}
\put(12.1,-2.75){1}
\put(11.8,-3.75){-1}
\put(12.1,-4.75){1}
\put(12.1,-5.75){1}
\put(11.8,-6.75){-1}
\put(11.8,-7.75){-1}
\put(7.,-9){\text{table
(5)}}
\end{picture}
$$
$$
$$
$$
$$
$$
\end{example}

\section{Traces in the category ${\cal D}$}

\begin {Def}\label{trace}\,[13] The trace of a morphism $ T : V
\longrightarrow V $ for any object $ V $ in ${\cal {D}}$ is defined by
\end{Def}
\vspace{.8cm}
$$
$$
$$
$$
$$
$$
$$
$$
$$
\setlength{\unitlength}{0.5cm}
\begin{picture}(5,2)\thicklines
\put(1.4,2.5){\coev}
\put(4,5.05){\T}
\put(4.75,6.09){\line(0,1){.5}}
\put(4.75,3.95){\line(0,1){.4}}
\put(7.8,4){\line(0,1){2.6}}
\put(4.75,3.95){\line(3,-2){3}}
\put(6.5,3.13){\line(3,2){1.3}}
\put(4.78,2){\line(3,2){1.2}}
\put(7.8,-.6){\line(0,1){.5}}
\put(7.78,1.58){\line(0,1){.37}}
\put(4.4,-.99){\eval}
\put(7,.6){\thi}
\put(4.75,-.6){\line(0,1){2.6}}
\put(7.4,7.3){\text {$V^*$}}
\put(4.7,7.3){\text {$V$}}
\end{picture}\,\,\,\,\,\qquad\qquad \,\qquad \,\,\quad
$$
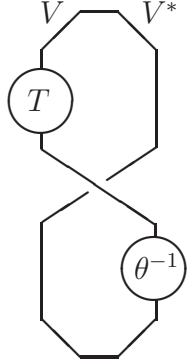
\,\, \centerline{figure 3}

\begin{theorem} \label{trc11} \hspace{0.2cm}
If we evaluate the diagram of definition \ref{trace} in ${\cal
{D}}$, we find
$$
{\rm trace}(T)=\sum_{\xi\in \,\text{basis of} \,\,V}
 \hat\xi\,\,\big(\,T\,(\,\xi\,)\,\big) \ .
$$
\end{theorem}
\textbf{Proof.}\,\,\, Begin with
$$
{\rm coev}_V(1)=\sum_{\xi\in \,\text{basis of} \,\,V} \xi\, \achl\,\tilde\tau(
\,\bbb^L ,\bbb \,)^{-1} \otimes \,\hat\xi=\sum_{\xi\in \,\text{basis of} \,\, V} \xi\, \achl\,\tau(
\,\b^L ,\b \,)^{-1} \otimes \,\hat\xi,
$$
and applying $T \otimes {\rm id}$ to this gives
\begin{equation*}
\begin{split}
\sum_{\xi\in \,\text{basis of} \,\, V}\,T\,\big( \xi\, \achl\,\tau(
\,\b^L ,\b \,)^{-1}\,\big) \otimes \hat\xi \,&=\sum_{\xi\in \,\text{basis of} \,\, V}
T( \,\xi\,)\, \achl\,\tau(
\,\b^L ,\b \,)^{-1} \otimes \hat\xi.
\end{split}
\end{equation*}
Next apply the braiding map to the last equation to get
\begin{equation}\label{gg7}
\begin{split}
\sum_{\xi\in \,\text{basis of} \,\, V}
 \Psi\,\big(\,T( \,\xi\,)\, \achl\,\tau(
\,\b^L ,\b \,)^{-1} \otimes \hat\xi \,\big)=\sum_{\xi\in \,\text{basis of} \,\, V}
\hat\xi \,\achl\,(\langle\xi^{'}\rangle\,\acl\, |\hat\xi\,| \,)^{-1}\,
\otimes\,\xi^{'}\achl \,|\,\hat\xi\,|
\end{split}
\end{equation}
where \,$\,\,\xi^{'}=T( \,\xi\,)\, \achl\,\tau(
\,\b^L ,\b \,)^{-1} $, so
\begin{equation}
\begin{split}
\langle \,\xi^{'}\rangle\,&=\langle \,T( \,\xi\,)\, \achl\,\tau(
\,\b^L ,\b \,)^{-1}\rangle\,=\langle \,T( \,\xi\,)\, \acbl\,\tau(
\,\b^L ,\b \,)^{-1}\rangle\,\\
&=\langle \,T( \,\xi\,)\,\rangle\, \acl\,\tau(
\,\b^L ,\b \,)^{-1}=\langle  \,\xi\,\rangle\, \acl\,\tau(
\,\b^L ,\b \,)^{-1}.
\end{split}
\end{equation}
To calculate $|\,\hat\xi\,| $ we start with
\begin{equation*}
\begin{split}
 \|\hat\xi\|=\bbb^{L}&=(\,\bbi \,\b\,)^L=\bb\,\tau(\,\b^L ,\b \,)^{-1}\,\b^L ,\,
\end{split}
\end{equation*}
which implies that\,\,
$|\hat\xi\,| =\tau(\,\b^L ,\b \,) \,\bbi$.
Then
\begin{equation*}
\begin{split}
\hat\xi \, \achl \,(\langle\xi^{'}\rangle\,\acl\, |\hat\xi\,| \,)^{-1}\,
&=\hat\xi \, \achl \,\big(\langle  \,\xi\,\rangle\, \acl\,\tau(
\,\b^L ,\b \,)^{-1}\,\tau(\,\b^L ,\b \,) \,\bbi \,\big)^{-1}\\
&=\hat\xi \, \achl \,\big(\langle  \,\xi\,\rangle\, \acl\,
\bbi \,\big)^{-1},
\\
\xi^{'}\,\acl\, |\hat\xi|&=\big(T( \xi)\, \achl\,\tau(
\,\b^L ,\b \,)^{-1} \big)\,\achl \big(\tau(\,\b^L ,\b ) \bbi \,\big)
=T( \xi)\, \achl\,\bbi,
\end{split}
\end{equation*}
which gives
\begin{equation}\label{gg8}
\begin{split}
\sum_{\xi\in \,\text{basis of} \,\, V}
\hat\xi \achl(\langle\xi^{'}\rangle\acl\, |\hat\xi| \,)^{-1}\,
\otimes \xi^{'}\achl \,|\hat\xi| =\sum_{\xi\in \,\text{basis of} \,\, V}
\hat\xi  \achl \,\big(\langle  \xi\,\rangle \acl\,\bbi \,\big)^{-1}
\otimes T( \xi) \achl\,\bbi.
\end{split}
\end{equation}
Next
\begin{equation}
\begin{split}
\,\,  \theta^{-1}\big( \,T( \,\xi\,)\, \achl\,\bbi \,\big)
&=\big( \,T( \,\xi\,)\, \achl\,\bbi \,\big)
\achl\,\big\| \,T( \,\xi\,)\, \achl\,\bbi \, \big\|^{-1} \qquad \,\,\,\,\,\,\,\,\,\\
&=\big( \,T( \,\xi\,)\, \achl\,\bbi \,\big)
\achl\,\big(\| \,T( \,\xi\,)\,\| \actl\,\bbi \, \big)^{-1}\\
&=\,T( \,\xi\,)\, \achl\,\bbi \,
\,\big(\|  \,\xi\,\| \actl\,\bbi \, \big)^{-1}\\
&= \,T( \,\xi\,)\, \achl\,\bbi \,
\big(\,\bb \,\bbi \,\b \, \bbi \big)^{-1}\\
&=\,T( \,\xi\,)\, \achl\,\bbi \,
\,\bb  \,\bi =\,T( \,\xi\,)\, \achl \,\bi,
\end{split}
\end{equation}
and finally we need to calculate
\begin{equation}
\begin{split}
{\rm eval} &\big(\hat\xi  \achl \big(\langle  \xi\rangle \acl \bbi \big)^{-1}
\otimes T( \xi) \achl \bi \big)= \big(\hat\xi  \achl
 \big(\langle  \xi\rangle \acl \bbi \big)^{-1}\, \big)
\, (T( \xi) \achl \bi)\ .\\
\end{split}
\end{equation}
We know, from the definition of the action on $V^{*}$, that
\begin{equation}
\Big(  \hat \xi \achl\big(\|T(\xi) \|\actr x \big)\Big)
\big(T(\xi)\achl x \big)=\hat \xi \big(T(\xi) \big).
\end{equation}
If we put $x=\bi$, we want to show that $\|T(\xi) \|\actr x
=\big(\langle \xi\rangle \acl \bbi \big)^{-1}$, so
$$
\|\xi\|\actl x=\bbi \b \actl \bi=\b \bbi=(\b \acr \bbi)
(\b \acl \bbi)=v^{'} t^{'},
$$
which implies that $t^{'}=\b \acl \bbi$,\, and hence
\begin{equation}
\begin{split}
\|T(\xi) \|\actr x&=\|\xi\|\actr x=\bbi \b \actr \bi=
t \bi {t^{'}}^{-1}\\
&=\b \bi(\b \acl \bbi)^{-1}=(\b \acl \bbi)^{-1}.
\end{split} \quad\square
\end{equation}

\section{ Characters in the category ${\cal D}$}

\begin{Def}\label{chardef}
   \cite{MajBook} The right adjoint action in ${\cal{D}}$ of
the algebra $D$ on itself is defined by
\end{Def}
$$
$$
$$
$$
$$
$$
$$
$$
$$
$$\qquad\qquad\qquad\qquad
\setlength{\unitlength}{0.5cm}
\begin{picture}(5,2)\thicklines
\put(10,2.38){\line(0,1){1}}
\put(8.2,4.){\eval}
\put(8.9,5.25){\leval}
\put(8.6,4.4){\line(1,1){.85}}
\put(9.85,5.55){\line(1,1){.74}}
\put(13.5,5.75){\line(0,1){3.45}}
\put(9.3,5.8){\line(0,1){5.4}}
\put(10.6,6.25){\line(0,1){.75}}
\put(10.5,8.15){\line(0,1){1.1}}
\put(12,10.2){\line(0,1){1}}
\put(7.1,5.1){\coev}
\put(9.85,7.4){\S}
\put(9,11.4){\text {$D$}}
\put(11.7,11.4){\text {$D$}}
\put(9.1,1.15){\text\rm{figure 4}}
\end{picture}\qquad\qquad \qquad\qquad\qquad
\vspace{-.8cm}
\begin{Def}
The character $\chi_{V}$ of an object $ V $ in ${\cal{D}}$ is defined by
\end{Def}
\vspace{.5cm}
$$
$$
$$
$$
$$
\setlength{\unitlength}{0.5cm}
\begin{picture}(5,2)\thicklines
\put(1.4,.8){\coev}
\put(4.75,3.95){\line(0,1){.94}}
\put(7.8,4){\line(0,1){.9}}
\put(4.75,3.95){\line(3,-2){3}}
\put(6.5,3.13){\line(3,2){1.3}}
\put(4.78,2){\line(3,2){1.2}}
\put(7.8,-2.11){\line(0,1){.37}}
\put(7.78,1.58){\line(0,1){.37}}
\put(9.4,1.58){\line(0,1){4.9}}
\put(7.78,-.03){\line(0,1){.82}}
\put(4.4,-2.5){\eval}
\put(7,-1.04){\thi}
\put(6.4,.45){\tri}
\put(4.75,-2.11){\line(0,1){4.1}}
\put(9.2,6.7){\text {$D$}}
\put(7.4,5.7){\text {$V^*$}}
\put(4.7,5.7){\text {$V$}}
\end{picture}\,\,\,\,\,\qquad\qquad \,\qquad \,\,\quad
$$
$$
$$
$$
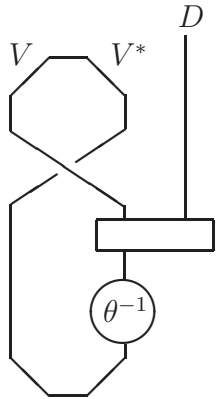
\,\, \centerline{figure 5}
$$
\vspace{-1.4cm}
\begin {Lem}\,\, For an object $V$ in
${\cal{D}}$  we have
\vspace{.5cm}
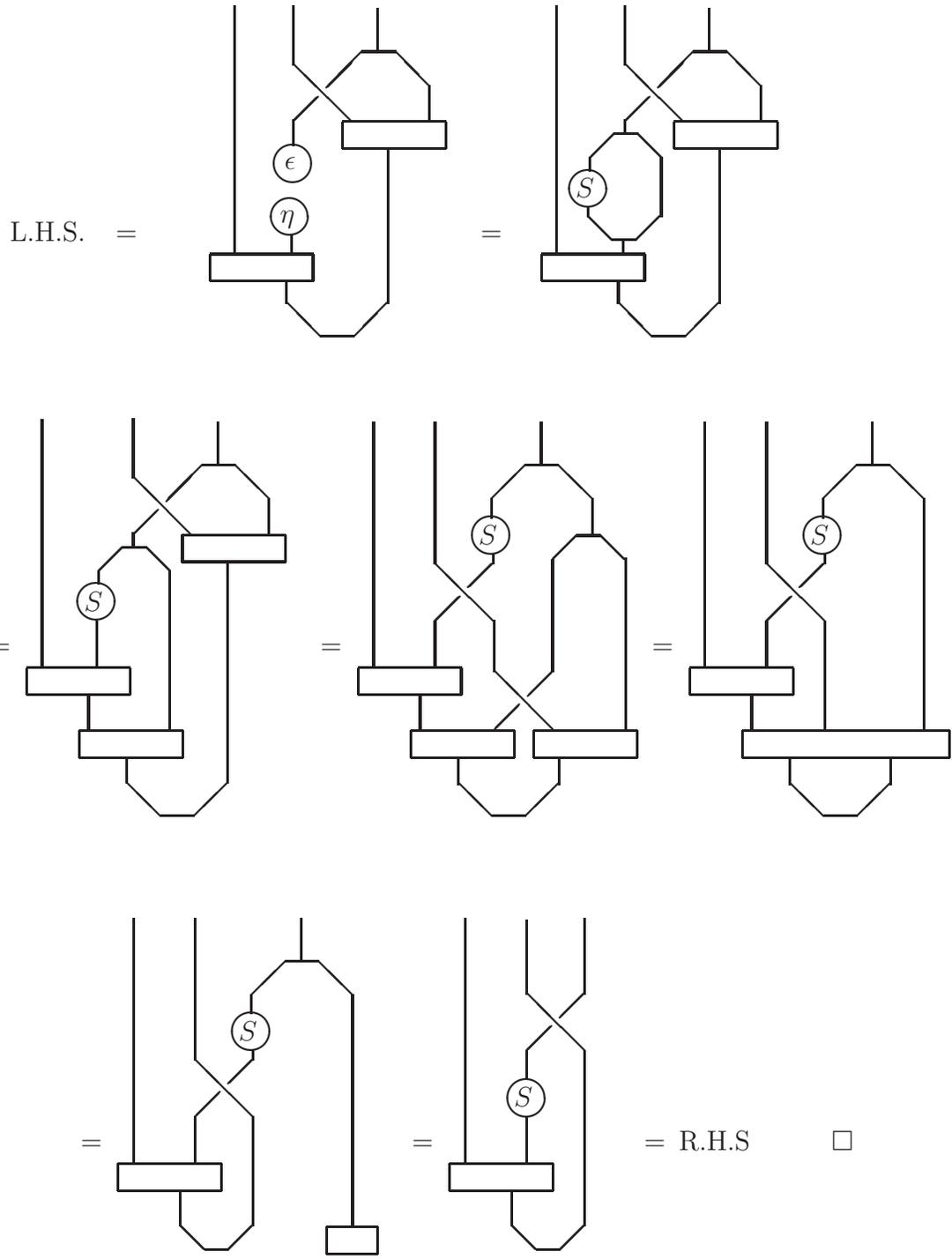
$$
$$
$$
$$
$$
\setlength{\unitlength}{0.5cm}
\begin{picture}(5,2)\thicklines
\put(8.5,.1){\line(0,1){.65}}
\put(7.78,1.58){\line(0,1){4.5}}
\put(9.4,1.58){\line(0,1){4.5}}
\put(5.1,-.29){\eval}
\put(6.4,.45){\tri}
\put(5.45,.1){\line(0,1){6}}
\put(9.2,6.4){\text {$D$}}
\put(7.5,6.4){\text {$V$}}
\put(5.2,6.4){\text {$V^{*}$}}
\put(11.2,2){\text {$=$}}
\put(14.5,6.5){\text {$V^{*}$}}
\put(16.2,6.5){\text {$V$}}
\put(18.1,6.5){\text {$D$}}
\put(15.75,.1){\line(0,1){.65}}
\put(14.78,1.58){\line(0,1){4.6}}
\put(18.25,5.5){\line(0,1){.75}}
\put(16.4,5.5){\line(0,1){.75}}
\put(16.4,1.58){\line(0,1){.6}}
\put(15.4,-.29){\eval}
\put(13.4,.45){\tri}
\put(15.6,2.55){\S}
\put(16.3,3.3){\line(0,1){.45}}
\put(17.5,4.8){\line(1,1){.75}}
\put(16.3,3.75){\line(1,1){.75}}
\put(18.8,.1){\line(0,1){3.}}
\put(16.4,5.5){\line(1,-1){2.4}}
\end{picture}\qquad \qquad\qquad\qquad \qquad
\qquad\qquad\qquad \qquad\qquad
$$
$$
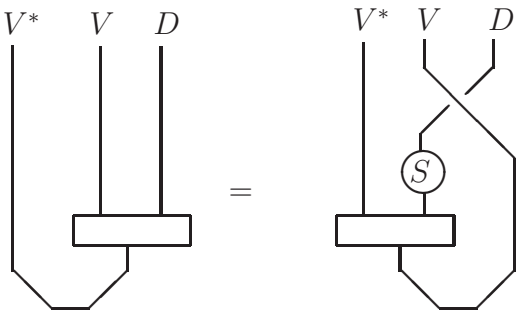
\,\, \centerline{ \rm{figure 6}}
$$
\end{Lem}
$$
$$
$$
$$
\textbf{Proof.}
$$
$$
$$
$$
$$
$$
$$\quad
\setlength{\unitlength}{0.5cm}
\begin{picture}(5,2)\thicklines
\put(0,2){\text {L.H.S.}}
\put(3.2,2){\text {$=$}}
\put(8.35,.1){\line(0,1){.65}}
\put(6.78,1.58){\line(0,1){7.5}}
\put(8.5,1.58){\line(0,1){.6}}
\put(8,-.29){\eval}
\put(5.4,.45){\tri}
\put(9.4,4.45){\tri}
\put(6.,4.6){\braid}
\put(6.25,2.6){\coev}
\put(12.65,5.58){\line(0,1){1.1}}
\put(8.55,4.85){\line(0,1){.8}}
\put(8.55,7.3){\line(0,1){1.8}}
\put(11.1,7.7){\line(0,1){1.3}}
\put(7.7,2.55){\et}
\put(7.8,4.15){\ep}
\put(11.4,.1){\line(0,1){4.7}}
\put(14.2,2){\text {$=$}}
\put(18.35,.1){\line(0,1){.65}}
\put(16.5,1.58){\line(0,1){7.5}}
\put(18.5,1.58){\line(0,1){.4}}
\put(18,-.29){\eval}
\put(15.4,.45){\tri}
\put(19.4,4.45){\tri}
\put(16.,4.6){\braid}
\put(16.25,2.6){\coev}
\put(22.65,5.58){\line(0,1){1.1}}
\put(18.55,5.25){\line(0,1){.4}}
\put(18.55,7.3){\line(0,1){1.8}}
\put(21.1,7.7){\line(0,1){1.3}}
\put(16.7,3.4){\S}
\put(19.65,2.7){\line(0,1){1.8}}
\put(17.45,2.7){\line(0,1){.3}}
\put(17.5,4.1){\line(0,1){.4}}
\put(15,1.55){\scoev}
\put(17.1,2.45){\seval}
\put(21.4,.1){\line(0,1){4.7}}
\end{picture}\qquad \qquad\qquad\qquad \qquad\qquad\qquad
\qquad\qquad\qquad \qquad\qquad\qquad \qquad\qquad\qquad \qquad
\qquad\qquad \qquad\qquad
$$
$$
$$
$$
$$
$$
$$
$$
$$
$$
\setlength{\unitlength}{0.5cm}
\begin{picture}(5,2)\thicklines
\put(.2,2){\text {$=$}}
\put(4.35,-1.9){\line(0,1){.8}}
\put(1.8,1.58){\line(0,1){7.5}}
\put(4,-2.29){\eval}
\put(.7,.45){\tri}
\put(2.3,-1.45){\tri}
\put(5.4,4.45){\tri}
\put(2.,4.6){\braid}
\put(2.25,2.6){\coev}
\put(8.65,5.58){\line(0,1){1.1}}
\put(4.55,5.25){\line(0,1){.4}}
\put(4.55,7.3){\line(0,1){1.8}}
\put(7.1,7.7){\line(0,1){1.3}}
\put(2.7,3.4){\S}
\put(5.65,-.3){\line(0,1){4.8}}
\put(3.2,-.3){\line(0,1){1.1}}
\put(3.45,1.6){\line(0,1){1.4}}
\put(3.5,4.1){\line(0,1){.4}}
\put(1,1.55){\scoev}
\put(7.4,-1.9){\line(0,1){6.7}}
\put(10.2,2){\text {$=$}}
\put(14.35,-1.95){\line(0,1){.8}}
\put(11.8,1.58){\line(0,1){7.4}}
\put(14,-2.29){\eval}
\put(10.7,.45){\tri}
\put(12.3,-1.45){\tri}
\put(16,-1.45){\tri}
\put(11.1,2){\braid}
\put(12.,2.6){\coev}
\put(18.4,5.58){\line(0,1){1.1}}
\put(15.35,6.1){\line(0,1){.6}}
\put(17.2,1.45){\line(0,1){3.4}}
\put(15.45,1.45){\line(0,1){1.55}}
\put(16.85,7.7){\line(0,1){1.3}}
\put(14.6,5.4){\S}
\put(12.9,-1.27){\braid}
\put(13.2,-.3){\line(0,1){1.1}}
\put(19.4,-.3){\line(0,1){5.2}}
\put(13.65,1.6){\line(0,1){1.4}}
\put(15.4,4.7){\line(0,1){.3}}
\put(13.65,4.7){\line(0,1){4.3}}
\put(14.7,1.9){\scoev}
\put(17.4,-1.95){\line(0,1){.8}}
\put(20.2,2){\text {$=$}}
\put(24.35,-1.95){\line(0,1){.8}}
\put(21.8,1.58){\line(0,1){7.4}}
\put(24,-2.29){\eval}
\put(20.7,.45){\tri}
\put(22.3,-1.45){\ltri}
\put(21.1,2){\braid}
\put(22.,2.6){\coev}
\put(25.35,6.1){\line(0,1){.6}}
\put(26.85,7.7){\line(0,1){1.3}}
\put(24.6,5.4){\S}
\put(23.2,-.3){\line(0,1){1.1}}
\put(28.4,-.3){\line(0,1){7}}
\put(25.4,-.3){\line(0,1){3.3}}
\put(23.65,1.6){\line(0,1){1.4}}
\put(25.4,4.7){\line(0,1){.3}}
\put(23.65,4.7){\line(0,1){4.3}}
\put(27.4,-1.95){\line(0,1){.8}}
\end{picture}\qquad \qquad\qquad\qquad \qquad\qquad\qquad
\qquad \qquad\qquad\qquad \qquad\qquad\qquad\qquad\qquad\qquad
\qquad\qquad\qquad \qquad\qquad\qquad
$$
$$
$$
$$
$$
$$
$$
$$
$$
$$
$$
$$
\setlength{\unitlength}{0.5cm}
\begin{picture}(5,2)\thicklines
\put(.2,2){\text {$=$}}
\put(1.8,1.58){\line(0,1){7.4}}
\put(2.87,-.57){\seval}
\put(.7,.45){\tri}
\put(7.,-1.45){\stri}
\put(1.1,2){\braid}
\put(2.,2.6){\coev}
\put(5.35,6.1){\line(0,1){.6}}
\put(6.85,7.7){\line(0,1){1.3}}
\put(4.6,5.4){\S}
\put(3.2,-.3){\line(0,1){1.1}}
\put(8.4,-.3){\line(0,1){7}}
\put(5.4,-.3){\line(0,1){3.3}}
\put(3.65,1.6){\line(0,1){1.4}}
\put(5.4,4.7){\line(0,1){.3}}
\put(3.65,4.7){\line(0,1){4.3}}
\put(10.2,2){\text {$=$}}
\put(11.8,1.58){\line(0,1){7.4}}
\put(12.87,-.57){\seval}
\put(10.7,.45){\tri}
\put(11.1,4){\braid}
\put(15.4,6.7){\line(0,1){2.3}}
\put(12.9,3.4){\S}
\put(13.2,-.3){\line(0,1){1.1}}
\put(15.4,-.3){\line(0,1){5.3}}
\put(13.65,1.6){\line(0,1){1.4}}
\put(13.65,4.15){\line(0,1){.85}}
\put(13.65,6.7){\line(0,1){2.25}}
\put(17.2,2){\text {$=$\,\,\,R.H.S}}
\put(21.2,2){$\qquad\square$}
\end{picture}\qquad\qquad\qquad \qquad\qquad\qquad
\qquad\qquad\qquad \qquad\qquad\qquad\qquad
$$
\begin{Proposition}\,\, The character
is right adjoint invariant, i.e.\ for an object $V$ in
${\cal{D}}$
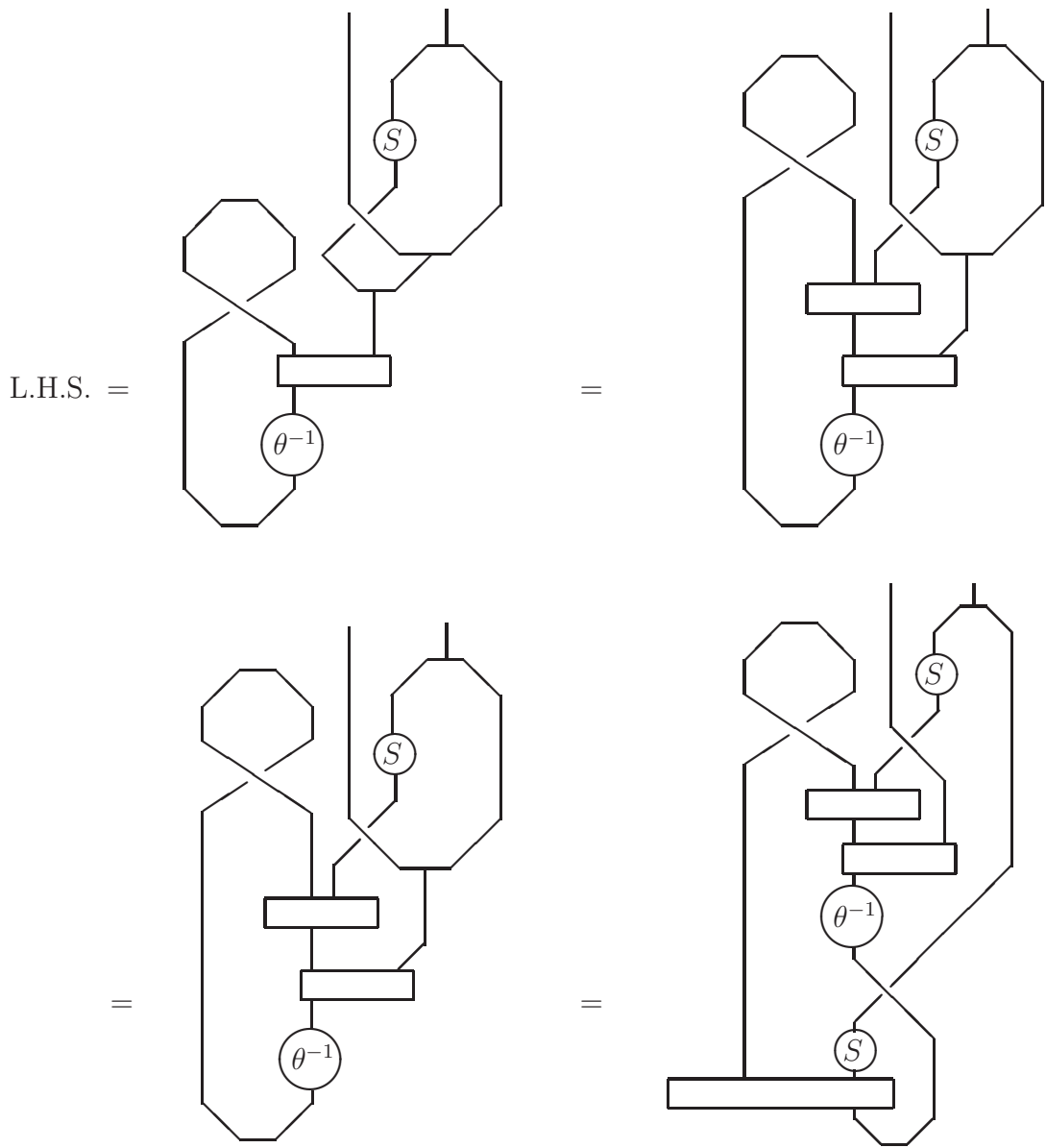
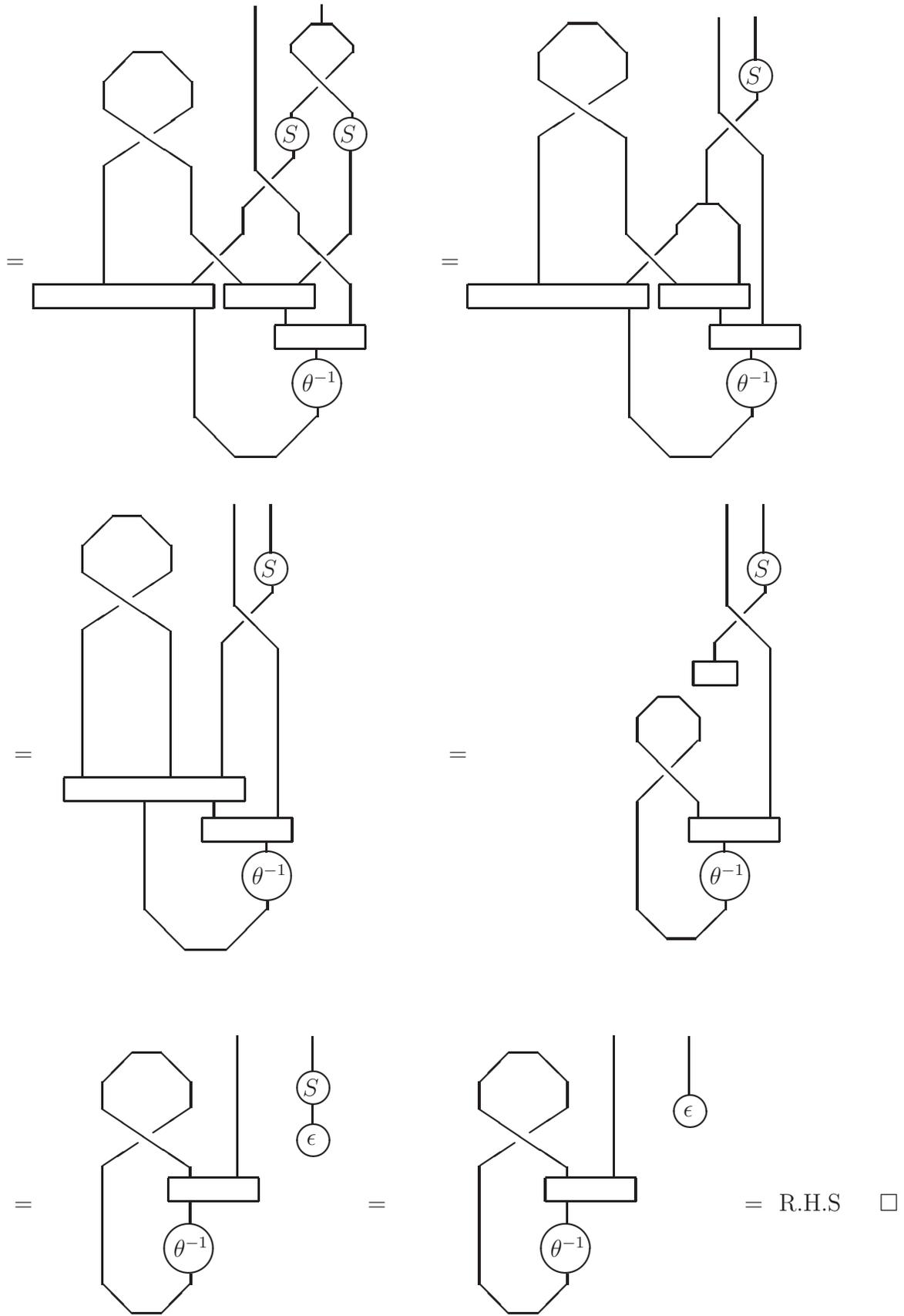
$$
$$
$$
\setlength{\unitlength}{0.5cm}
\begin{picture}(5,2)\thicklines
\put(1.2,1.8){\line(0,1){1.7}}
\put(2.25,1.8){\line(0,1){1.7}}
\put(1,1){\Ad}
\put(1.75,-.65){\line(0,1){1}}
\put(1,-1.5){\ch}
\put(.9,3.65){\text {$D$}}
\put(2.05,3.65){\text {$D$}}
\put(4.05,.65){\text{$=$}}
\put(6.75,1.85){\line(0,1){1.5}}
\put(6,1){\ch}
\put(8.75,1.85){\line(0,1){1.5}}
\put(8.,1.1){\ep}

\put(8.5,3.5){\text {$D$}}
\put(6.5,3.5){\text {$D$}}
\end{picture}\qquad \qquad \qquad\quad
$$
$$
\centerline{ \rm figure 7 \,\,\,\,\,\,\,\,\,}
$$
\end{Proposition}
\textbf{Proof.}
$$
$$
$$
$$
$$
$$
$$
$$
\setlength{\unitlength}{0.5cm}
\begin{picture}(5,2)\thicklines
\put(-.1,.45){\text{L.H.S.\,\,\,$=$}}
\put(1.4,.8){\coev}
\put(4.75,3.95){\line(0,1){.94}}
\put(7.8,4){\line(0,1){.9}}
\put(4.75,3.95){\line(3,-2){3}}
\put(6.5,3.13){\line(3,2){1.3}}
\put(4.78,2){\line(3,2){1.2}}
\put(7.8,-2.11){\line(0,1){.37}}
\put(7.78,1.58){\line(0,1){.37}}
\put(10,1.58){\line(0,1){1.8}}
\put(7.78,-.03){\line(0,1){.82}}
\put(4.4,-2.5){\eval}
\put(7,-1.04){\thi}
\put(6.7,.45){\tri}
\put(4.75,-2.11){\line(0,1){4.1}}
\put(8.2,4.){\eval}
\put(8.9,5.25){\leval}
\put(8.6,4.4){\line(1,1){.85}}
\put(9.85,5.55){\line(1,1){.74}}
\put(13.5,5.75){\line(0,1){3.45}}
\put(9.3,5.8){\line(0,1){5.3}}
\put(10.6,6.25){\line(0,1){.75}}
\put(10.5,8.15){\line(0,1){1.1}}
\put(12,10.2){\line(0,1){1}}
\put(7.1,5.1){\coev}
\put(9.85,7.4){\S}
\put(15.7,.45){\text{$=$}}
\put(16.9,4.8){\coev}
\put(20.25,7.95){\line(0,1){.94}}
\put(23.3,8){\line(0,1){.9}}
\put(20.25,7.95){\line(3,-2){3}}
\put(22.,7.13){\line(3,2){1.3}}
\put(20.28,6){\line(3,2){1.2}}
\put(23.3,-2.11){\line(0,1){.37}}
\put(23.28,1.58){\line(0,1){1.17}}
\put(23.28,3.58){\line(0,1){2.37}}
\put(23.88,3.58){\line(0,1){.95}}
\put(26.4,2.28){\line(0,1){2.1}}
\put(23.28,-.03){\line(0,1){.82}}
\put(19.9,-2.5){\eval}
\put(22.5,-1.04){\thi}
\put(22.35,.45){\tri}
\put(21.35,2.45){\tri}
\put(20.25,-2.11){\line(0,1){8.1}}
\put(23.9,5.25){\leval}
\put(23.9,4.5){\line(1,1){.71}}
\put(24.85,5.55){\line(1,1){.74}}
\put(25.65,1.6){\line(1,1){.74}}
\put(28.5,5.75){\line(0,1){3.45}}
\put(24.3,5.8){\line(0,1){5.3}}
\put(25.6,6.25){\line(0,1){.75}}
\put(25.5,8.15){\line(0,1){1.1}}
\put(27,10.2){\line(0,1){1}}
\put(22.1,5.1){\coev}
\put(24.85,7.4){\S}
\end{picture}
$$
$$
$$
$$
$$
$$
$$
$$
$$
$$
$$
$$
\setlength{\unitlength}{0.5cm}
\begin{picture}(5,2)\thicklines
\put(2.7,.45){\text{$=$}}
\put(1.9,4.8){\coev}
\put(5.25,7.95){\line(0,1){.94}}
\put(8.3,8){\line(0,1){.9}}
\put(5.25,7.95){\line(3,-2){3}}
\put(7.,7.13){\line(3,2){1.3}}
\put(5.28,6){\line(3,2){1.2}}
\put(8.3,-2.11){\line(0,1){.37}}
\put(8.28,1.58){\line(0,1){1.17}}
\put(8.28,3.58){\line(0,1){2.37}}
\put(8.88,3.58){\line(0,1){.95}}
\put(11.4,2.28){\line(0,1){2.1}}
\put(8.28,-.03){\line(0,1){.82}}
\put(4.9,-2.5){\eval}
\put(7.5,-1.04){\thi}
\put(7.35,.45){\tri}
\put(6.35,2.45){\tri}
\put(5.25,-2.11){\line(0,1){8.1}}
\put(8.9,5.25){\leval}
\put(8.9,4.5){\line(1,1){.71}}
\put(9.85,5.55){\line(1,1){.74}}
\put(10.65,1.6){\line(1,1){.74}}
\put(13.5,5.75){\line(0,1){3.45}}
\put(9.3,5.8){\line(0,1){5.3}}
\put(10.6,6.25){\line(0,1){.75}}
\put(10.5,8.15){\line(0,1){1.1}}
\put(12,10.2){\line(0,1){1}}
\put(7.1,5.1){\coev}
\put(9.85,7.4){\S}
\put(15.7,.45){\text{$=$}}
\put(16.9,6.1){\coev}
\put(20.25,9.25){\line(0,1){.94}}
\put(23.3,9.3){\line(0,1){.9}}
\put(20.25,9.25){\line(3,-2){3}}
\put(22.,8.47){\line(3,2){1.3}}
\put(20.28,7.3){\line(3,2){1.2}}
\put(23.3,-.11){\line(0,1){.24}}
\put(22.6,-.81){\S}
\put(20.75,-.85){\braid}
\put(25.05,1.9){\line(1,1){2.6}}
\put(17.5,-2.55){\ltri}
\put(23.3,-1.41){\line(0,1){.24}}
\put(23.28,5.08){\line(0,1){.67}}
\put(23.28,6.58){\line(0,1){.67}}
\put(23.88,6.58){\line(0,1){.45}}
\put(25.8,5.1){\line(0,1){1.75}}
\put(23.28,3.950){\line(0,1){.32}}
\put(23.28,1.9){\line(0,1){.32}}
\put(22.97,-2.8){\seval}
\put(23.3,-2.54){\line(0,1){.28}}
\put(25.5,-2.54){\line(0,1){2.24}}
\put(25.51,-.31){\line(-1,1){.8}}
\put(22.5,2.9){\thi}
\put(22.35,3.95){\tri}
\put(21.35,5.45){\tri}
\put(20.25,-1.41){\line(0,1){8.7}}
\put(25.8,6.85){\line(-1,1){1.5}}
\put(23.9,7.05){\line(1,1){.71}}
\put(24.9,8.05){\line(1,1){.74}}
\put(27.65,4.45){\line(0,1){6.5}}
\put(24.3,8.3){\line(0,1){4.}}
\put(25.6,8.75){\line(0,1){.45}}
\put(25.5,10.35){\line(0,1){.6}}
\put(26.6,11.7){\line(0,1){.6}}
\put(23.,8.){\scoev}
\put(24.85,9.6){\S}
\end{picture}
$$
$$
$$
$$
$$
$$
$$
$$
$$
$$
\setlength{\unitlength}{0.5cm}
\begin{picture}(5,2)\thicklines
\put(.7,.45){\text{$=$}}
\put(11.58,8.9){\line(0,1){.65}}
\put(.72,2.8){\coev}
\put(4.07,5.95){\line(0,1){.94}}
\put(7.12,6){\line(0,1){.9}}
\put(4.07,5.95){\line(3,-2){3}}
\put(5.82,5.13){\line(3,2){1.3}}
\put(4.1,4){\line(3,2){1.2}}
\put(11.4,-2.7){\line(0,1){.37}}
\put(11.45,-4.68){\line(0,1){.28}}
\put(1.,-1.25){\ltri}
\put(7.1,1.6){\line(0,1){2.37}}
\put(8.88,1.6){\line(0,1){.95}}
\put(10.8,1.6){\line(0,1){.75}}
\put(7.2,-4.7){\line(0,1){3.75}}
\put(10.35,-1.5){\line(0,1){.55}}
\put(6.8,-5.23){\leval}
\put(10.7,-3.7){\thi}
\put(9.35,-2.65){\tri}
\put(7.6,-1.25){\tri}
\put(4.07,-.11){\line(0,1){4.1}}
\put(10.8,2.35){\line(-1,1){1.5}}
\put(8.9,2.55){\line(1,1){.71}}
\put(9.9,3.55){\line(1,1){.71}}
\put(9.3,3.8){\line(0,1){5.7}}
\put(12.68,7.88){\line(0,1){.25}}
\put(10.52,7.88){\line(0,1){.27}}
\put(7.6,4.6){\lbraid}
\put(8.,5.2){\scoev}
\put(9.85,4.9){\S}
\put(11.85,4.9){\S}
\put(10.62,4.23){\line(0,1){.25}}
\put(10.55,5.62){\line(0,1){.2}}
\put(12.66,5.62){\line(0,1){.2}}
\put(12.58,1.62){\line(0,1){2.9}}
\put(12.58,-1.5){\line(0,1){1.4}}
\put(4.55,-1.1){\braid}
\put(8.25,-1.1){\braid}
\put(15.7,.45){\text{$=$}}
\put(15.72,3.8){\coev}
\put(19.07,6.95){\line(0,1){.94}}
\put(22.12,7){\line(0,1){.9}}
\put(19.07,6.95){\line(3,-2){3}}
\put(20.82,6.13){\line(3,2){1.3}}
\put(19.1,5){\line(3,2){1.2}}
\put(26.4,-2.7){\line(0,1){.37}}
\put(26.45,-4.68){\line(0,1){.28}}
\put(16.,-1.25){\ltri}
\put(22.1,1.6){\line(0,1){3.37}}
\put(24.88,2.65){\line(0,1){1.9}}
\put(26,-.1){\line(0,1){2.05}}
\put(23.85,1.6){\line(0,1){.35}}
\put(22.2,-4.7){\line(0,1){3.75}}
\put(25.35,-1.5){\line(0,1){.55}}
\put(21.8,-5.23){\leval}
\put(25.7,-3.7){\thi}
\put(24.35,-2.65){\tri}
\put(22.6,-1.25){\tri}
\put(19.07,-.11){\line(0,1){5.1}}
\put(26.8,4.35){\line(-1,1){1.5}}
\put(24.9,4.55){\line(1,1){.71}}
\put(25.9,5.55){\line(1,1){.71}}
\put(21.33,-1.){\scoev}
\put(25.85,6.9){\S}
\put(25.3,5.8){\line(0,1){3.3}}
\put(26.62,6.23){\line(0,1){.25}}
\put(26.55,7.62){\line(0,1){1.5}}
\put(26.8,-1.5){\line(0,1){5.85}}
\put(19.55,-1.1){\braid}
\end{picture}
$$
$$
$$
$$
$$
$$
$$
$$
$$
$$
$$
$$
\vspace{.5cm}
\setlength{\unitlength}{0.5cm}
\begin{picture}(5,2)\thicklines
\put(.7,.45){\text{$=$}}
\put(-.28,3.8){\coev}
\put(3.07,6.95){\line(0,1){.94}}
\put(6.12,7){\line(0,1){.9}}
\put(3.07,6.95){\line(3,-2){3}}
\put(4.82,6.13){\line(3,2){1.3}}
\put(3.1,5){\line(3,2){1.2}}
\put(9.4,-2.7){\line(0,1){.37}}
\put(9.45,-4.68){\line(0,1){.28}}
\put(1.8,-1.25){\ltri}
\put(6.1,-.1){\line(0,1){5.03}}
\put(7.88,-.1){\line(0,1){4.64}}
\put(5.2,-4.7){\line(0,1){3.75}}
\put(7.6,-1.5){\line(0,1){.55}}
\put(4.8,-5.23){\leval}
\put(8.7,-3.7){\thi}
\put(6.55,-2.65){\tri}
\put(3.07,-.11){\line(0,1){5.1}}
\put(9.8,4.35){\line(-1,1){1.5}}
\put(7.9,4.55){\line(1,1){.71}}
\put(8.9,5.55){\line(1,1){.71}}
\put(8.3,5.8){\line(0,1){3.5}}
\put(8.85,6.9){\S}
\put(9.62,6.23){\line(0,1){.25}}
\put(9.55,7.62){\line(0,1){1.7}}
\put(9.8,-1.5){\line(0,1){5.83}}
\put(15.7,.45){\text{$=$}}
\put(19.7,-1.){\scoev}
\put(19.25,-2.15){\lbraid}
\put(22.2,1.1){\line(0,1){.85}}
\put(24.35,1.13){\line(0,1){.82}}
\put(25.2,-2.7){\line(0,1){.37}}
\put(25.25,-4.7){\line(0,1){.3}}
\put(23.5,2.75){\stri}
\put(24.88,3.9){\line(0,1){.64}}
\put(22.2,-4.7){\line(0,1){3.75}}
\put(24.3,-1.5){\line(0,1){.55}}
\put(21.85,-5.08){\eval}
\put(24.5,-3.7){\thi}
\put(23.35,-2.65){\tri}
\put(26.8,4.35){\line(-1,1){1.5}}
\put(24.9,4.55){\line(1,1){.71}}
\put(25.9,5.55){\line(1,1){.71}}
\put(25.3,5.8){\line(0,1){3.5}}
\put(25.85,6.9){\S}
\put(26.62,6.23){\line(0,1){.25}}
\put(26.55,7.62){\line(0,1){1.7}}
\put(26.8,-1.5){\line(0,1){5.83}}
\end{picture}
$$
$$
$$
$$
$$
$$
$$
$$
$$
$$
\vspace{.7cm}
\setlength{\unitlength}{0.5cm}
\begin{picture}(5,2)\thicklines
\put(.7,.45){\text{$=$}}
\put(.4,.8){\coev}
\put(3.75,3.95){\line(0,1){.94}}
\put(6.8,4){\line(0,1){.9}}
\put(3.75,3.95){\line(3,-2){3}}
\put(5.5,3.13){\line(3,2){1.3}}
\put(3.78,2){\line(3,2){1.2}}
\put(6.8,-2.11){\line(0,1){.37}}
\put(6.78,1.58){\line(0,1){.37}}
\put(8.4,1.58){\line(0,1){4.9}}
\put(6.78,-.03){\line(0,1){.82}}
\put(3.4,-2.5){\eval}
\put(6,-1.04){\thi}
\put(5.4,.45){\tri}
\put(3.75,-2.11){\line(0,1){4.1}}
\put(11,5.25){\line(0,1){1.2}}
\put(10.3,4.5){\S}
\put(10.3,2.7){\ep}
\put(11,3.45){\line(0,1){.68}}
\put(12.9,.45){\text{$=$}}
\put(13.4,.8){\coev}
\put(16.75,3.95){\line(0,1){.94}}
\put(19.8,4){\line(0,1){.9}}
\put(16.75,3.95){\line(3,-2){3}}
\put(18.5,3.13){\line(3,2){1.3}}
\put(16.78,2){\line(3,2){1.2}}
\put(19.8,-2.11){\line(0,1){.37}}
\put(19.78,1.58){\line(0,1){.37}}
\put(21.4,1.58){\line(0,1){4.9}}
\put(19.78,-.03){\line(0,1){.82}}
\put(16.4,-2.5){\eval}
\put(19,-1.04){\thi}
\put(18.4,.45){\tri}
\put(16.75,-2.11){\line(0,1){4.1}}
\put(23.3,3.7){\ep}
\put(24,4.45){\line(0,1){2.}}
\put(25.9,.45){\text{$=$\,\,\,\,R.H.S}}
\put(28.9,.45){$\qquad\square$}
\end{picture}
$$
$$

\begin{Proposition}\,\, The character of a tensor product
 of representations is the product of the
 characters, i.e.\ for two objects $V$ and $W$ in
${\cal{D}}$
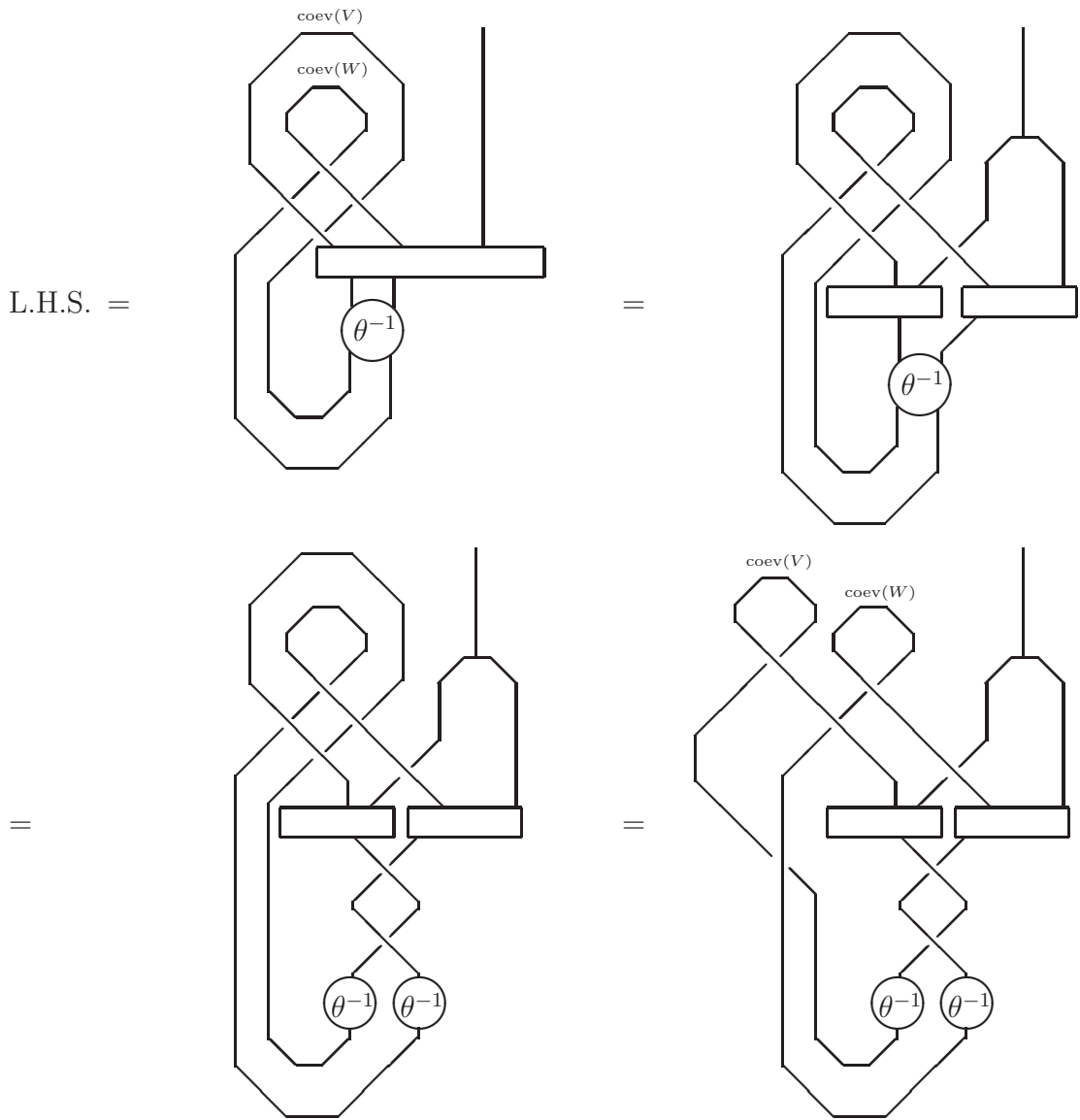
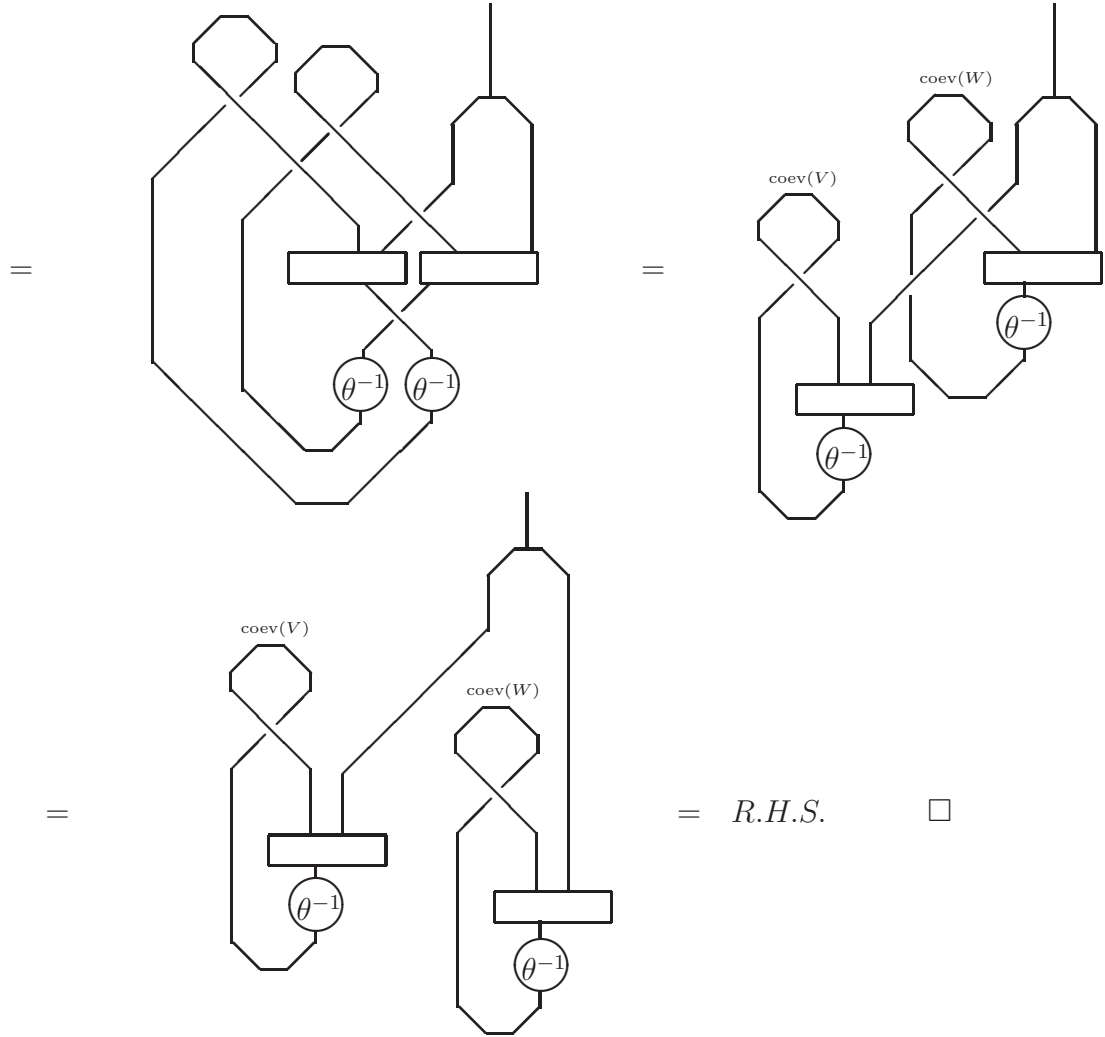
$$
$$
$$\qquad
\setlength{\unitlength}{0.5cm}
\begin{picture}(5,2)\thicklines
\put(1.1,1.2){\line(0,1){2}}
\put(.4,-.1){\chvw}
\put(.8,3.4){\text {$D$}}
\put(7.3,3.5){\text {$D$}}
\put(3.5,1.4){\text {$=$}}
\put(7.65,2.3){\line(0,1){1}}
\put(4.1,-1.4){\scoev}
\put(6.6,.75){\line(0,1){.8}}
\put(8.77,.75){\line(0,1){.8}}
\put(5.9,-.1){\ch}
\put(8.,-.1){\chw}
\end{picture}\qquad \qquad \qquad\quad
$$
$$
\centerline{ \qquad\rm figure 8 \,\,\,\qquad\,\,\,}
$$
\end{Proposition}
\textbf{Proof.}
$$
$$
$$
$$
$$
$$
\setlength{\unitlength}{0.5cm}
\begin{picture}(5,2)\thicklines
\put(.1,.45){\text{L.H.S.\,\,\,$=$}}
\put(8.,8.45){\text{\tiny{coev$(V)$}}}
\put(8.,7.){\text{\tiny{coev$(W)$}}}
\put(13.1,2.25){\line(0,1){6}}
\put(2.1,.95){\lcoev}
\put(5.2,2.95){\scoev}
\put(7.7,5.4){\line(0,1){.5}}
\put(9.9,5.4){\line(0,1){.5}}
\put(6.7,4.5){\line(0,1){2.2}}
\put(10.9,4.6){\line(0,1){2.1}}
\put(9.,2.25){\line(-1,1){2.3}}
\put(10.9,2.25){\line(-1,1){3.2}}
\put(9.8,3.55){\line(1,1){1.1}}
\put(8.92,4.5){\line(1,1){.95}}
\put(7.88,3.58){\line(1,1){.71}}
\put(8.78,2.68){\line(1,1){.71}}
\put(7.88,3.58){\line(1,1){.71}}
\put(7.2,1.24){\line(1,1){1.21}}
\put(6.3,2){\line(1,1){1.31}}
\put(7.9,1.1){\ltri}
\put(9.5,.6){\line(0,1){.8}}
\put(10.65,.55){\line(0,1){.85}}
\put(9.33,-.22){\thi}
\put(9.44,-1.74){\line(0,1){1.1}}
\put(10.55,-2.47){\line(0,1){1.75}}
\put(7.2,-1.74){\line(0,1){3.}}
\put(6.3,-2.47){\line(0,1){4.5}}
\put(6.9,-2.){\seval}
\put(5.9,-3.){\leval}
\put(16.9,.45){\text{$=$}}
\put(27.9,5.25){\line(0,1){3}}
\put(24.35,1.58){\scoev}
\put(26.9,2.95){\line(0,1){1.6}}
\put(29.,1.15){\line(0,1){3.4}}
\put(25.,1.15){\line(1,1){.8}}
\put(26.15,2.25){\line(1,1){.74}}
\put(17.1,.95){\lcoev}
\put(20.2,2.95){\scoev}
\put(22.7,5.4){\line(0,1){.5}}
\put(24.9,5.4){\line(0,1){.5}}
\put(21.7,4.5){\line(0,1){2.2}}
\put(25.9,4.6){\line(0,1){2.1}}
\put(24.4,1.85){\line(-1,1){2.7}}
\put(24.4,1.15){\line(0,1){.7}}
\put(27,1.15){\line(-1,1){4.3}}
\put(24.8,3.55){\line(1,1){1.1}}
\put(23.92,4.5){\line(1,1){.95}}
\put(22.88,3.58){\line(1,1){.71}}
\put(23.78,2.68){\line(1,1){.71}}
\put(22.88,3.58){\line(1,1){.71}}
\put(22.2,1.24){\line(1,1){1.21}}
\put(21.3,2){\line(1,1){1.31}}
\put(21.9,-.){\tri}
\put(25.6,0.){\tri}
\put(24.5,-.9){\line(0,1){1.2}}
\put(25.65,-.95){\line(0,1){.3}}
\put(25.65,-.65){\line(1,1){.96}}
\put(24.33,-1.72){\thi}
\put(24.44,-3.24){\line(0,1){1.1}}
\put(25.55,-3.97){\line(0,1){1.75}}
\put(22.2,-3.24){\line(0,1){4.5}}
\put(21.3,-3.97){\line(0,1){6.}}
\put(21.9,-3.5){\seval}
\put(20.9,-4.5){\leval}
\end{picture}
$$
$$
$$
$$
$$
$$
$$
$$
$$
$$
\setlength{\unitlength}{0.5cm}
\begin{picture}(5,2)\thicklines
\put(.1,.45){\text{$=$}}
\put(12.9,5.25){\line(0,1){3}}
\put(9.35,1.58){\scoev}
\put(11.9,2.95){\line(0,1){1.6}}
\put(14.,1.15){\line(0,1){3.4}}
\put(10.,1.15){\line(1,1){.8}}
\put(11.15,2.25){\line(1,1){.74}}
\put(2.1,.95){\lcoev}
\put(5.2,2.95){\scoev}
\put(7.7,5.4){\line(0,1){.5}}
\put(9.9,5.4){\line(0,1){.5}}
\put(6.7,4.5){\line(0,1){2.2}}
\put(10.9,4.6){\line(0,1){2.1}}
\put(9.4,1.85){\line(-1,1){2.7}}
\put(9.4,1.15){\line(0,1){.7}}
\put(12,1.15){\line(-1,1){4.3}}
\put(9.8,3.55){\line(1,1){1.1}}
\put(8.92,4.5){\line(1,1){.95}}
\put(7.88,3.58){\line(1,1){.71}}
\put(8.78,2.68){\line(1,1){.71}}
\put(7.88,3.58){\line(1,1){.71}}
\put(7.2,1.24){\line(1,1){1.21}}
\put(6.3,2){\line(1,1){1.31}}
\put(6.9,-.){\tri}
\put(10.4,0.){\tri}
\put(7.,-2.43){\braid}
\put(9.52,-1.68){\line(0,1){.24}}
\put(11.32,-1.68){\line(0,1){.24}}
\put(7.,-4.4){\braid}
\put(9.52,-3.5){\line(0,1){.1}}
\put(11.32,-3.5){\line(0,1){.1}}
\put(10.62,-4.4){\sthi}
\put(8.72,-4.4){\sthi}
\put(9.44,-5.24){\line(0,1){.3}}
\put(11.32,-5.2){\line(0,1){.3}}
\put(7.2,-5.24){\line(0,1){6.5}}
\put(6.3,-5.97){\line(0,1){8.}}
\put(6.9,-5.5){\seval}
\put(10.5,-5.97){\line(1,1){.8}}
\put(5.9,-6.5){\leval}
\put(16.9,.45){\text{$=$}}
\put(20.3,7.75){\text{\tiny{coev$(V)$}}}
\put(23.,6.9){\text{\tiny{coev$(W)$}}}
\put(27.9,5.25){\line(0,1){3}}
\put(24.35,1.58){\scoev}
\put(26.9,2.95){\line(0,1){1.6}}
\put(29.,1.15){\line(0,1){3.4}}
\put(25.,1.15){\line(1,1){.8}}
\put(26.15,2.25){\line(1,1){.74}}

\put(20.2,2.95){\scoev}
\put(22.7,5.4){\line(0,1){.5}}
\put(24.9,5.4){\line(0,1){.5}}
\put(17.5,3.75){\scoev}
\put(20,6.2){\line(0,1){.5}}
\put(22.2,6.2){\line(0,1){.5}}
\put(21.25,5.3){\line(1,1){.95}}
\put(18.9,3.1){\line(1,1){1.95}}
\put(22.2,-1.23){\line(-1,1){.71}}
\put(21.,-.2){\line(-1,1){2.1}}
\put(18.9,1.85){\line(0,1){1.25}}

\put(24.4,1.85){\line(-1,1){4.4}}
\put(24.4,1.15){\line(0,1){.7}}
\put(27,1.15){\line(-1,1){4.3}}

\put(23.92,4.5){\line(1,1){.95}}
\put(22.88,3.58){\line(1,1){.71}}

\put(22.88,3.58){\line(1,1){.71}}

\put(21.3,2){\line(1,1){1.31}}
\put(21.9,-.){\tri}
\put(25.4,0.){\tri}
\put(22.,-2.43){\braid}
\put(24.52,-1.68){\line(0,1){.24}}
\put(26.32,-1.68){\line(0,1){.24}}
\put(22.,-4.4){\braid}
\put(24.52,-3.5){\line(0,1){.1}}
\put(26.32,-3.5){\line(0,1){.1}}
\put(25.62,-4.4){\sthi}
\put(23.72,-4.4){\sthi}
\put(24.44,-5.24){\line(0,1){.3}}
\put(26.32,-5.2){\line(0,1){.3}}
\put(22.2,-5.24){\line(0,1){4.}}
\put(21.3,-5.97){\line(0,1){8.}}
\put(21.9,-5.5){\seval}
\put(25.5,-5.97){\line(1,1){.8}}
\put(20.9,-6.5){\leval}
\end{picture}
$$
$$
$$
$$
$$
$$
$$
$$
$$
$$
$$
$$
$$
$$
\setlength{\unitlength}{0.5cm}
\begin{picture}(5,2)\thicklines
\put(.1,.45){\text{$=$}}
\put(12.9,5.25){\line(0,1){2.5}}
\put(9.35,1.58){\scoev}
\put(11.9,2.95){\line(0,1){1.6}}
\put(14.,1.15){\line(0,1){3.4}}
\put(10.,1.15){\line(1,1){.8}}
\put(11.15,2.25){\line(1,1){.74}}
\put(5.2,2.95){\scoev}
\put(7.7,5.4){\line(0,1){.5}}
\put(9.9,5.4){\line(0,1){.5}}
\put(2.5,3.75){\scoev}
\put(5,6.2){\line(0,1){.5}}
\put(7.2,6.2){\line(0,1){.5}}
\put(6.25,5.3){\line(1,1){.95}}
\put(3.9,3.1){\line(1,1){1.95}}
\put(7.23,-3.4){\line(-1,1){.91}}
\put(6.35,-4.18){\line(-1,1){2.45}}
\put(9.4,1.85){\line(-1,1){4.4}}
\put(9.4,1.15){\line(0,1){.7}}
\put(12,1.15){\line(-1,1){4.3}}
\put(8.92,4.5){\line(1,1){.95}}
\put(7.88,3.58){\line(1,1){.71}}
\put(7.88,3.58){\line(1,1){.71}}
\put(6.3,2){\line(1,1){1.31}}
\put(6.9,-.){\tri}
\put(10.4,0.){\tri}
\put(7.,-2.43){\braid}
\put(9.52,-1.68){\line(0,1){.24}}
\put(11.32,-1.68){\line(0,1){.24}}
\put(10.62,-2.55){\sthi}
\put(8.72,-2.55){\sthi}
\put(9.44,-3.44){\line(0,1){.3}}
\put(11.32,-3.4){\line(0,1){.3}}
\put(6.3,-2.55){\line(0,1){4.55}}
\put(3.9,-1.8){\line(0,1){4.9}}
\put(6.9,-3.7){\seval}
\put(10.5,-4.17){\line(1,1){.8}}
\put(5.9,-4.7){\leval}
\put(16.9,.45){\text{$=$}}
\put(20.3,3.){\text{\tiny{coev$(V)$}}}
\put(24.3,5.65){\text{\tiny{coev$(W)$}}}
\put(27.9,5.25){\line(0,1){2.5}}
\put(24.35,1.58){\scoev}
\put(26.9,2.95){\line(0,1){1.6}}
\put(29.,1.15){\line(0,1){3.4}}
\put(23.,-.75){\line(1,1){2.8}}
\put(26.15,2.25){\line(1,1){.74}}
\put(21.5,1.65){\scoev}
\put(24,4.1){\line(0,1){.5}}
\put(26.2,4.1){\line(0,1){.5}}
\put(17.5,-1.){\scoev}
\put(20.03,1.45){\line(0,1){.5}}
\put(22.15,1.45){\line(0,1){.5}}
\put(17.1,-1.8){\lbraid}
\put(27,1.15){\line(-1,1){3.}}
\put(25.22,3.2){\line(1,1){.95}}
\put(24.12,2.15){\line(1,1){.8}}
\put(24.07,-1.75){\line(0,1){1.75}}
\put(24.09,.55){\line(0,1){1.6}}
\put(20.4,-3.5){\tri}
\put(25.4,0.){\tri}
\put(26.35,-.9){\sthi}
\put(22.29,-3.5){\line(0,1){.3}}
\put(22.15,-2.34){\line(0,1){1.75}}
\put(23.,-2.34){\line(0,1){1.6}}
\put(21.57,-4.4){\sthi}
\put(22.29,-5.24){\line(0,1){.3}}
\put(27.1,-1.75){\line(0,1){.3}}
\put(27.1,-.05){\line(0,1){.4}}
\put(20.05,-5.24){\line(0,1){4.65}}
\put(19.75,-5.51){\seval}
\put(23.7,-2.12){\eval}
\end{picture}
$$
$$
$$
$$
$$
$$
$$
$$

\setlength{\unitlength}{0.5cm}
\begin{picture}(5,2)\thicklines
\put(.1,-2){\text{$=$}}
 \put(5.3,3.){\text{\tiny{coev$(V)$}}}
\put(11.3,1.35){\text{\tiny{coev$(W)$}}}
\put(16.9,-2){\text{$=
\,\,\,\,R.H.S.$}}
\put(12.9,5.25){\line(0,1){1.5}}
\put(9.35,1.58){\scoev}
\put(11.88,3.05){\line(0,1){1.5}}
\put(14.,-3.85){\line(0,1){8.4}}
\put(8.,-.75){\line(1,1){3.88}}
\put(8.5,-2.65){\scoev}
\put(11,-.2){\line(0,1){.5}}
\put(13.2,-.2){\line(0,1){.5}}
\put(8.1,-3.5){\lbraid}
\put(10.75,-7.21){\seval}
\put(13.15,-3.85){\line(0,1){1.55}}
\put(11.05,-6.95){\line(0,1){4.65}}
\put(13.25,-6.98){\line(0,1){.45}}
\put(13.25,-5.1){\line(0,1){.45}}
\put(2.5,-1.){\scoev}
\put(5.03,1.45){\line(0,1){.5}}
\put(7.15,1.45){\line(0,1){.5}}
\put(2.1,-1.8){\lbraid}
\put(5.4,-3.5){\tri}
\put(11.4,-5.){\tri}
\put(12.55,-6.){\sthi}
\put(7.29,-3.5){\line(0,1){.3}}
\put(7.15,-2.34){\line(0,1){1.75}}
\put(8.,-2.34){\line(0,1){1.6}}
\put(6.57,-4.4){\sthi}
\put(7.29,-5.24){\line(0,1){.3}}
\put(5.05,-5.24){\line(0,1){4.65}}
\put(4.75,-5.51){\seval}
\put(21.9,-2){$\qquad\square$}
\end{picture}
$$
$$
$$
$$
$$
$$

\begin{theorem}\label{plpl} \hspace{0.2cm}
We have the following formula for the character;
$$
\chi_{_V}(\delta_y\otimes x)= \sum_{ \begin{subarray}\,\xi\in \,\text{basis of} \,\, V
\,\text{with} \,y=\b \bbi \end{subarray}}  \hat\xi\,(\,\xi\,\achl\,\bi \,x \, \b \,),\,
 \,\,\,\qquad
$$
for $xy =yx$, otherwise $\chi_{_V}(\delta_y\otimes x)=0$.
\end{theorem}
\textbf{Proof.}
Set $a=\delta_y\otimes x$.  To have $\chi_{_V}(a)\neq 0$
we must have $\|a\|=e$, i.e.\
$y=y \,\actl\,x$
which implies that $x$ and $y$ commute.  Assuming this, we continue
with the diagrammatic definition of the character, starting with
\begin{equation*}
\Big(\sum_{\xi\in \,\text{basis of} \,\,V} \xi \achl \tilde\tau(
\bbb^L ,\bbb )^{-1} \otimes \,\hat\xi\Big)\otimes
 a=\sum_{\xi\in \,\text{basis of} \,\, V}  \Big(\xi \achl \tau(
\b^L ,\b )^{-1} \otimes \,\hat\xi \Big)\otimes a \,.
\end{equation*}
Next we calculate
\begin{equation}
\begin{split}
\Psi\,\big(\,\xi\, \achl\,\tau(
\,\b^L ,\b \,)^{-1} \otimes \hat\xi \,\big)=
\hat\xi \,\achl\,(\langle\xi^{'}\rangle\,\acl\, |\hat\xi\,| \,)^{-1}\,
\otimes\,\xi^{'}\achl \,|\,\hat\xi\,|
\end{split}
\end{equation}
where \,\,$\,\,\xi^{'}= \,\xi\, \achl\,\tau(
\,\b^L ,\b \,)^{-1} \,\,$,\,\,\,\,so
\begin{equation}
\begin{split}
\langle \xi^{'}\rangle &=\langle \xi \achl\,\tau(
\,\b^L ,\b \,)^{-1}\rangle\,=\langle \xi \acbl\,\tau(
\,\b^L ,\b \,)^{-1}\rangle\,
=\langle  \xi \rangle \acl\,\tau(
\b^L ,\b )^{-1}.
\end{split}
\end{equation}
From a previous calculation we know that $|\,\hat\xi\,|
=\tau(\,\b^L ,\b \,) \,\bbi $, so
\begin{equation*}
\begin{split}
\hat\xi \, \achl \,(\langle\xi^{'}\rangle\,\acl\, |\hat\xi\,| \,)^{-1}\,
&=\hat\xi \, \achl \,\big(\langle  \,\xi\,\rangle\, \acl\,\tau(
\,\b^L ,\b \,)^{-1}\,\tau(\,\b^L ,\b \,) \,\bbi \,\big)^{-1}\\
&=\hat\xi \, \achl \,\big(\langle  \,\xi\,\rangle\, \acl\,\bbi \,\big)^{-1}\\
\xi^{'}\,\acl\, |\hat\xi\,|&= \big(\,\xi\, \achl\,\tau(
\,\b^L ,\b \,)^{-1} \,\big)\achl \,\big(\,\tau(\,\b^L ,\b \,) \,\bbi \,\big)
=\,\xi\, \achl\,\bbi\ ,
\end{split}
\end{equation*}
which gives the next stage in the evaluation of the diagram:
\begin{equation}
\begin{split}
\sum_{\xi\in \,\text{basis of} \,\, V}
\Psi\,\big(\xi \achl\,\tau(
\,\b^L ,\b \,)^{-1} &\otimes \hat\xi \,\big)\otimes a \\
=&\sum_{\xi\in \,\text{basis of} \,\, V}
\Big(\hat\xi  \achl \,\big(\langle  \xi\rangle \acl\,\bbi \big)^{-1}
\otimes \,\xi\, \achl\,\bbi\Big) \otimes a.
\end{split}
\end{equation}
Now we apply the associator to the last equation to get
\begin{equation*}
\begin{split}
\sum_{\xi\in \,\text{basis of} \,\, V}
\Phi \,&\Big( \,\big(\hat\xi  \achl \,(\langle  \xi\rangle \acl\,\bbi )^{-1}
\otimes \,\xi\, \achl\,\bbi \big) \otimes a \Big)
\\=&\sum_{\xi\in \,\text{basis of} \,\, V}
\hat\xi  \achl \,(\langle  \xi\rangle \acl\,\bbi )^{-1}
 \, \tilde\tau(\,\|\,\xi\, \achl\,\bbi \,\|^L ,\,\|a\| \,)\otimes \,
\big(\xi\, \achl\,\bbi  \otimes a \big)
\\=&\sum_{\xi\in \,\text{basis of} \,\, V}
\hat\xi \, \achl \,(\langle  \xi\rangle \, \acl\,\bbi )^{-1}
 \, \tau(\,\langle\,\xi\, \achl\,\bbi \,\rangle ,\,e \,)\otimes \,
\big( \,\xi\, \achl\,\bbi  \otimes a  \,\big)
\\=&\sum_{\xi\in \,\text{basis of} \,\, V}
\hat\xi \, \achl \,(\langle  \xi\rangle \, \acl\,\bbi )^{-1}
 \,\otimes \,
\big( \,\xi\, \achl\,\bbi  \otimes (\delta_y\otimes x)  \,\big)
\end{split}
\end{equation*}
as $\tau(\langle\,\xi \achl\,\bbi \,\rangle ,e)
=e$.  Now apply the action $\achl$ to $ \xi\, \achl\,\bbi \otimes
(\delta_y\otimes x) $ to get
\begin{equation}
\begin{split}
(\xi\, \achl\,\bbi ) \achl (\delta_y\otimes x)=
\delta_{y\,,\|\,\xi\, \achl\,\bbi \,\|}\,
(\xi\, \achl\,\bbi \,) \achl \,x =
\delta_{y\,,\|\,\xi\,\| \actl\,\bbi }\,\,
\xi\, \achl\,\bbi x,
\end{split}
\end{equation}
and to get a non-zero answer we must have
\begin{equation}
\begin{split}
y=\|\,\xi\,\| \actl\,\bbi=\bbi \,\b\, \actl\,\bbi
=\bb\,\bbi\,\b\,\bbi=\b \,\bbi.
\end{split}
\end{equation}
Thus the character of $V$ is given by
\begin{equation*}
\begin{split}
\chi_{_V}(\delta_y\otimes x)= \sum_{ \begin{subarray}\,\xi\in \,
\text{basis of} \,\, V \, \text{with} \,y=\b \bbi \end{subarray}}
{\rm eval}\Big(\hat\xi  \achl \,(\langle  \xi\rangle \, \acl\,\bbi
)^{-1} \, \otimes \theta^{-1}\,( \xi \achl\,\bbi \,x ) \Big)\, .
\end{split}
\end{equation*}
Next
\begin{equation*}
\begin{split}
\theta^{-1}\,(\, \xi\, \achl\,\bbi \,x \,)\,&=
(\, \xi\, \achl\,\bbi \,x \,)\,\achl\,
\big\|\, \xi\, \achl\,\bbi \,\,x \,\big\|^{-1}\\
&=(\, \xi\, \achl\,\bbi \,x \,)\,\achl\,
\big(\|\, \xi\,\| \actl\,\bbi \,\,x \,\big)^{-1}\\
&=(\, \xi\, \achl\,\bbi \,x \,)\,\achl\,
\big(x^{-1}\,\bb \,\bbi \b \,\bbi \,\,x \big)^{-1}\\
&=\, \xi\, \achl\,\bbi \,\,x \,x^{-1}\, \bb \, \bi \,x
=\, \xi\, \achl\, \, \bi \,x.
\end{split}
\end{equation*}
Now we need to calculate
${\rm eval} \big(\,\hat\xi \achl (\langle \xi\rangle \acl \bbi )^{-1}
\otimes \xi \achl \bi x \big)$.
Start with $\bbb \actl \bi x=\b \bbi \actl x =\b \bbi$, as we only
have nonzero summands for $y=\b \bbi$.  Then
\begin{eqnarray*}
 &&   {\rm eval} \big(\,\hat\xi \achl (\langle \xi\rangle \acl \bbi )^{-1}
\otimes \xi \achl \bi x \big) \\ & & =\,
{\rm eval}\big( \big(\,\hat\xi \achl (\langle \xi\rangle \acl \bbi )^{-1}
\otimes \xi \achl \bi x \big)\achl\langle \xi\rangle\big) \\
 & & =\,
{\rm eval} \big(\,\hat\xi \achl (\langle \xi\rangle \acl \bbi )^{-1}
(\b \bbi\actr \langle \xi\rangle) \otimes \xi \achl \bi x \langle \xi\rangle
\big)\ .
\end{eqnarray*}
To find $ \b \bbi \actr \b$, first find $\b \bbi \actl \b =\bbi \b$, so
\begin{equation*}
\begin{split}
\b \bbi \actr \b =(\b \acr \bbi)(\b \acl \bbi) \actr \b
=(\b \acl \bbi) \b \bi =\b \acl \bbi\ .\quad\square
\end{split}
\end{equation*}

\begin {Lem}{}\hspace{0.1cm}Let $V$ be an  object
in ${\cal D}$.  For $\delta_y\otimes x \in D$ the character of $V$
is given by the following formula, where  $y=su^{-1}$ with $s \in
M $ and $u \in G $:
$$
\chi_{_V}(\delta_y\otimes x)=
 \sum_{ \begin{subarray}\,\xi\in \,\text{basis of} \,\, V_{u^{-1}s}
\, \end{subarray}}  \hat\xi\,(\,\xi\,\achl\,s^{-1} \,x \, s \,)
=\chi_{_{V_{u^{-1}s}}}(\,s^{-1} \,x \, s \,)\,
 \,\,\,\qquad
$$
where  $xy=yx$, otherwise $\chi_{_V}(\delta_y\otimes x)=0$.
  Here $ \chi_{_{V_{u^{-1}s}}}$ is the group representation character
of the representation $V_{u^{-1}s}$ of the group stab$(u^{-1}s)$.
\end{Lem}
\textbf{Proof.}\hspace{0.5cm}
From theorem \ref{plpl}, we know that
$$
\chi_{_V}(\delta_y\otimes x)= \sum_{ \begin{subarray}\,\xi\in \,\text{basis of} \,\, V
\,\text{with} \,y=\b \bbi \end{subarray}}  \hat\xi\,(\,\xi\,\achl\,\bi \,x \, \b \,),\,
 \,\,\,\qquad
$$
for with $xy=yx$.  Set $s=\b$ and $u=\bb$, so $ y=su^{-1}$.  We note
that $s^{-1} x s$ is in $stab(u^{-1}s)$, because
$$
u^{-1}s \actl s^{-1} x s =s^{-1} x^{-1} s u^{-1}ss^{-1} x s=s^{-1} x^{-1}
xs u^{-1} s=u^{-1} s.
$$
It just remains to note that $\|\xi\|=|\xi|^{-1}\langle \xi\rangle=
u^{-1}s$.\quad$\square$

\section{Modular Categories}
  Let
 ${\cal{M}}$ be a semisimple ribbon category.  For objects $V$ and $W$
 in ${\cal{M}}$ define ${\tilde S}_{VW} \in
 \underline{\bf 1}$ by
$$
$$
$$
$$
$$
$$
$$
\vspace{.7cm}
\setlength{\unitlength}{0.5cm}
\begin{picture}(5,2)\thicklines
\put(1.42,3.38){\coev}
\put(4.76,3.95){\line(0,1){3.53}}

\put(4.75,3.95){\line(3,-2){3}}
\put(6.5,3.13){\line(3,2){1.25}}
\put(4.78,2){\line(3,2){1.2}}
\put(7.78,.9){\line(0,1){1.03}}
\put(7.8,-2.11){\line(0,1){1.03}}
\put(4.4,-2.5){\eval}
\put(7,-.24){\ti}
\put(4.75,-2.11){\line(0,1){4.1}}
\put(1.1,2.3){\text {$\tilde S_{VW}$\,=}}

\put(5.2,8.7){\text{ \small{coev($V$)}}}

\put(5.5,4.76){\braidi}
\put(5.5,2.97){\braidi}
\put(6.1,3.36){\coev}
\put(12.5,4){\line(0,1){3.46}}
\put(9.48,4.){\line(3,-2){3.}}
\put(11.2,3.13){\line(3,2){1.3}}
\put(9.48,2){\line(3,2){1.2}}
\put(12.48,.9){\line(0,1){1.1}}
\put(12.5,-2.11){\line(0,1){1.03}}
\put(9.1,-2.5){\eval}
\put(11.7,-.24){\wi}
\put(9.45,-2.11){\line(0,1){4.1}}

\put(10.,8.7){\text {\small{coev($W$)}}}

\end{picture}\qquad\qquad\qquad\qquad\qquad \,\qquad \,\,\quad
$$
\vspace{-1.cm}

\centerline{\quad figure 9}

There are standard results \cite{Bak,TWen}:
$$
{\tilde S}_{VW}={\tilde S}_{WV}= {\tilde S}_{V^*W^*}={\tilde
S}_{W^*V^*},\,\,\, {\tilde S}_{V\underline 1}={\rm dim} (V)\ .
$$
Here ${\rm dim} (V)$ is the trace in ${\cal{M}}$ of the identity
map on $V$.

\begin {Def}\,  We call an object $U$ in an abelian category ${\cal
{M}}$ simple if, for any $V$ in ${\cal {M}}$, any injection $V
\hookrightarrow U$ is either $0$ or an isomorphism \cite{Bak}.  A
semisimple category
 is an abelian category whose objects split as a direct sums
of simple objects  \cite{TWen}.
\end {Def}

\begin {Def}\,\cite{Bak} A modular category is
a semisimple ribbon
category ${\cal {M}}$ satisfying the following properties:\\
1- There are only a finite number of isomorphism
classes of simple objects in ${\cal {M}}$.\\
2- Schur's lemma holds, i.e. the morphisms between simple objects
are zero unless they are isomorphic, in which case the morphisms
are a multiple of the identity.\\
3- The matrix ${\tilde S}_{VW}$ with indices in isomorphism
classes of simple objects is invertible.
\end {Def}

\begin {Def}\,\cite{Bak} For a simple object $V$, the ribbon map on $V$ is a
multiple of the identity, and we use ${\Theta _V}$ for the scalar
multiple. The numbers $P^{\pm}$ are defined as the following sums
over simple isomorphism classes:
$$
P^{\pm}=\sum_{V} {\Theta _V}^{\pm 1} ({\rm dim}(V))^2\ ,
$$
and the matrices $T$ and $C$ are defined using the Kronecker delta
function by
$$
T_{VW}=\delta_{VW}\, {\Theta _V}\ ,\quad C_{VW}=\delta_{V W^*}\ .
$$
\end {Def}

\begin {theorem}\,\cite{Bak} In a modular category, if we define the matrix $S$ by
$$
S=\frac {\tilde S } {\sqrt {P^{+}P^{-}} }\ ,
$$
then we have the following matrix equations:
$$
( S T)^{3}=\sqrt {\frac {P^{+}}{P^{-}}}\  S^{2}\ ,\qquad S^{2}=C\ ,\qquad
CT=TC\ , \qquad C^{2}=1\ .
$$
\end {theorem}
We now give some results which allow us to calculate the matrix
$\tilde{S}$ in ${\cal{D}}$ .

\begin{Lem}
\end{Lem}
\vspace{.5cm}
$$
$$
$$
$$
\qquad\qquad\qquad
\setlength{\unitlength}{0.5cm}
\begin{picture}(5,2)\thicklines
\put(-4.,.6){\lcoev}
\put(1,8.2){\text {coev($V^{*}$)}}
\put(4.83,4.25){\line(0,1){2.1}}
\put(.58,4.25){\line(0,1){2.1}}
\put(6,5.2){\text {$=$}}
\put(8.5,9.){\text {coev($V$)}}
\put(6.,5.){\scoev}
\put(10.68,6.95){\line(0,1){1.}}
\put(8.48,6.95){\line(0,1){1.}}
\put(6.45,4){\braidi}
\put(9.98,6.25){\line(1,1){.71}}
\put(9.2,6.28){\line(-1,1){.73}}
\put(10.5,4.4){\line(0,1){.6}}
\put(8.7,2.3){\line(0,1){2.7}}
\put(9.75,3.7){\u}
\put(10.5,2.3){\line(0,1){1.}}
\put(12,5.){\text {where}}
\put(15.4,5){\text {u\,=}}
\put(16.3,2.75){\scoev}
\put(17.03,4.){\seval}
\put(14.45,3.1){\lbraid}
\put(20.95,2.7){\line(0,1){3.}}
\put(17.37,6.38){\line(0,1){2.}}
\end{picture}

\vspace{-1.5cm} \centerline{ figure 10 \qquad }\
\hspace{-0.7cm}\textbf{Proof.}\qquad\qquad
$$
$$
$$
$$
$$
$$
$$
\setlength{\unitlength}{0.5cm}
\begin{picture}(5,2)\thicklines
\put(-6,3.2){\text {R. H. S.$=$}}
\put(-4,5.){\scoev}
\put(.68,6.95){\line(0,1){1.}}
\put(-1.52,6.95){\line(0,1){1.}}
\put(-3.55,4){\braidi}
\put(-.02,6.25){\line(1,1){.71}}
\put(-.8,6.28){\line(-1,1){.73}}
\put(-1.3,-.7){\line(0,1){5.7}}
\put(-6.42,2.6){\uu}
\put(-1.6,-1.4){\text{$V^{*}$}}
\put(3.8,-1.4){\text{$V^{**}$}}
\put(5,3.2){\text {$=$}}
\put(.42,4.6){\uu}
\put(2.75,1.25){\lcoev}
\put(8.65,-1.43){\braidi}
\put(12.7,1.3){\line(0,1){4.55}}
\put(12.7,5.85){\line(-1,1){1.1}}
\put(13.5,3.2){\text {$=$}}
\put(12.8,2.){\scoev}
\put(14.92,3.95){\seval}
\put(14.65,2.73){\braidi}
\put(15.28,4.2){\line(0,1){.75}}
\put(14.08,-.25){\lcoev}
\put(18.85,-1.48){\braidi}
\put(21.7,.7){\line(-1,1){3.1}}
\put(22.9,1.25){\line(0,1){4.22}}
\end{picture}\qquad\qquad\qquad \qquad\qquad\qquad
\qquad\qquad\qquad
$$
$$
$$
$$
$$
$$\qquad\qquad\qquad\qquad\qquad
\setlength{\unitlength}{0.5cm}
\begin{picture}(5,2)\thicklines
\put(-0,1.2){\text {$=$}}
\put(-1.95,-1.){\lcoev}
\put(2.64,3.22){\line(0,1){1.5}}
\put(8.15,3.45){\line(-1,1){1.8}}
\put(8.15,-1.35){\line(0,1){4.8}}
\put(2.3,2.95){\seval}

\put(2.2,.2){\scoev}

\put(6.87,-1.32){\line(0,1){4.47}}
\put(9.8,1.2){\text {$=$}}

\put(8.43,-1.5){\coev}
\put(11.8,1.08){\line(0,1){1.5}}
\put(14.82,1.08){\line(0,1){1.5}}
\put(16.5,1.2){\text {$=\,L. H. S.$}}
\put(20.2,){$\qquad\square$}
\end{picture}\qquad\qquad\qquad \qquad\qquad\qquad
\qquad\qquad\qquad
$$

\begin{Lem}
\end{Lem}
\qquad
$$\qquad\qquad\qquad
\setlength{\unitlength}{0.5cm}
\begin{picture}(5,2)\thicklines
\put(4,2.5){\text {$V$}}
\put(7,2.5){\text {$V^*$}}
\put(7.3,.95){\line(0,1){1.4}}
\put(4.25,.45){\line(0,1){1.9}}
\put(3.5,-.3){\u}
\put(6.55,-.2){\tti}
\put(7.3,-1.65){\line(0,1){0.6}}
\put(4.25,-1.62){\line(0,1){.9}}
\put(3.9,-2.02){\eval}
\put(9.5,-.2){\text {$=$}}
\put(11.2,3.1){\text {$V$}}
\put(13.3,3.1){\text {$V^*$}}
\put(8.45,-1.){\lbraid}
\put(12.8,-1.5){\ti}
\put(11.05,-3.07){\seval}
\put(11.37,-2.8){\line(0,1){3.}}
\put(11.37,2.3){\line(0,1){0.6}}
\put(13.5,2.3){\line(0,1){0.6}}
\put(13.55,-2.8){\line(0,1){0.5}}
\put(13.5,-.35){\line(0,1){0.55}}

\put(-3.5,-6.8){\text {where  \,\,\,\quad\,\,u=}}
\put(-4.,-6.2){\uu} \put(8.5,-6.8){\text {and  \,\,\,\quad\,
\,${\theta_{V^*}}^{\!\!\!-1}=$}} \put(16.8,-6.5){\ti}
\put(15.,-8.08){\seval} \put(15.05,-7.78){\scoev}
\put(17.53,-5.35){\line(0,1){0.5}}
\put(17.53,-7.8){\line(0,1){0.5}} \put(15.35,-7.8){\line(0,1){4}}
\put(19.75,-8.85){\line(0,1){4.}}
\end{picture}\qquad\qquad\qquad
\qquad\qquad\qquad\qquad \qquad\qquad
\qquad\qquad \qquad
$$
$$
\qquad\qquad\qquad\qquad\qquad\qquad \centerline{figure 11}
$$
$$
$$
$$
$$
\textbf{Proof.}\hspace{0.5cm}\qquad\qquad
$$
$$
$$
$$
$$\qquad\qquad\qquad\qquad
\setlength{\unitlength}{0.5cm}
\begin{picture}(5,2)\thicklines
\put(-4.5,-.2){\text {$L. H. S. =$}}
\put(-.08,2.83){\line(0,1){3}}
\put(4.05,2.83){\line(0,1){3}}
\put(-7.,0.47){\uu}
\put(3.1,-3.38){\leval}
\put(-11.85,7.73){\thsi}
\put(7.7,-2.83){\line(1,1){.72}}
\put(9.5,-.2){\text {$ =$}}
\put(17.33,6.18){\line(0,1){.4}}
\put(15.22,6.18){\line(0,1){.4}}
\put(-1.85,7.73){\thsi}
\put(8.8,-1.78){\scoev}
\put(13.1,-3.38){\leval}
\put(17.7,-2.83){\line(1,1){.72}}
\put(13.5,-2.83){\line(0,1){4}}
\put(13.1,-4.38){\leval}
\put(17.7,-3.83){\line(1,1){1.72}}
\put(13.5,-3.8){\line(-1,1){2.2}}
\put(11.27,-1.63){\line(0,1){2.8}}
\put(19.45,-2.13){\line(0,1){4.1}}
\put(12.28,2.9){\lbraid}
\put(14.08,2.95){\line(1,1){1.72}}
\put(19.45,1.98){\line(-1,1){2.2}}
\end{picture}\qquad\qquad\qquad
\qquad\qquad\qquad\qquad \qquad\qquad
\qquad\qquad \qquad
$$
$$
$$
$$
$$
$$
$$
$$
$$
$$\qquad\qquad\qquad\qquad\qquad
\setlength{\unitlength}{0.5cm}
\begin{picture}(5,2)\thicklines
\put(-2.5,-.2){\text {$ =$}}
\put(5.33,6.18){\line(0,1){.4}}
\put(3.22,6.18){\line(0,1){.4}}
\put(-13.85,7.73){\thsi}
\put(8.63,-2.13){\line(0,1){2.95}}
\put(0.28,2.9){\lbraid}
\put(2.08,2.95){\line(1,1){1.72}}
\put(8.6,0.83){\line(-1,1){3.8}}
\put(6.1,-2.4){\seval}
\put(9.5,-.2){\text {$ =$}}
\put(8.45,1.){\lbraid}
\put(12.8,.5){\ti}
\put(11.05,-1.07){\seval}
\put(11.37,-0.8){\line(0,1){3.}}
\put(11.37,4.3){\line(0,1){0.6}}
\put(13.5,4.3){\line(0,1){0.6}}
\put(13.55,-0.8){\line(0,1){0.5}}
\put(13.5,1.65){\line(0,1){0.55}}
\put(15.5,-.2){\text {$ =R. H. S.$}}
\put(20.,-.2){$\qquad\square$}
\end{picture}\qquad\qquad\qquad
\qquad\qquad\qquad\qquad \qquad\qquad
\qquad\qquad \qquad
$$
$$
$$

\begin{Lem} For $V,W$ indecomposable objects in ${\cal D}$,
    $
    {\rm trace}(\Psi_{V^{*}W}\circ\Psi_{WV^{*}})\,=\,
    \tilde S_{VW}
    $.
\end{Lem}
\textbf{Proof.}\hspace{0.5cm}\qquad\qquad
$$
$$
$$
$$
$$
$$
$$
$$\qquad
\setlength{\unitlength}{0.5cm}
\begin{picture}(5,2)\thicklines
\put(.1,3.45){\text{L.H.S.\,\,\,$=$}}
\put(8,9.85){\text{\tiny{coev($V^*$)}}}
\put(8,11.4){\text{\tiny{coev($W$)}}}

\put(2.1,3.95){\lcoev}
\put(6.,8.96){\line(1,1){1.21}}
\put(3.45,4.45){\braid}
\put(3.45,6.2){\braid}
\put(5.2,5.95){\scoev}

\put(9.9,5.4){\line(0,1){3.5}}

\put(10.9,4.6){\line(0,1){5.1}}
\put(9.3,2.15){\line(-1,1){3.3}}
\put(10.7,2.5){\line(-1,1){3.2}}
\put(9.8,3.55){\line(1,1){1.1}}
\put(8.92,4.5){\line(1,1){.95}}
\put(7.98,3.68){\line(1,1){.703}}
\put(8.88,2.78){\line(1,1){.703}}

\put(7.2,1.24){\line(1,1){1.21}}
\put(6.3,2){\line(1,1){1.31}}

\put(9.33,.4){\line(0,1){1.75}}
\put(10.7,.5){\line(0,1){2.}}
\put(9.33,-.22){\thi}
\put(9.44,-1.74){\line(0,1){1.1}}
\put(10.55,-2.47){\line(0,1){1.75}}
\put(7.2,-1.74){\line(0,1){3.}}
\put(6.3,-2.47){\line(0,1){4.5}}
\put(6.9,-2.){\seval}
\put(5.9,-3.){\leval}
\put(15.9,3.45){\text{$=$}}

\put(17.1,3.95){\lcoev}
\put(20.2,5.95){\scoev}
\put(23.,11.4){\text{\tiny{coev($W$)}}}
\put(23.1,9.8){\text{\tiny{coev($V^*$)}}}
\put(22.7,8.4){\line(0,1){.5}}
\put(24.9,8.4){\line(0,1){.5}}
\put(21.7,7.5){\line(0,1){2.2}}
\put(25.9,7.6){\line(0,1){2.1}}
\put(24.4,4.85){\line(-1,1){2.7}}
\put(24.4,2.45){\line(0,1){2.4}}
\put(27,4.15){\line(-1,1){4.3}}
\put(24.8,6.55){\line(1,1){1.1}}
\put(23.92,7.5){\line(1,1){.95}}
\put(22.88,6.58){\line(1,1){.71}}
\put(23.78,5.68){\line(1,1){.71}}
\put(22.88,6.58){\line(1,1){.71}}
\put(22.2,4.24){\line(1,1){1.21}}
\put(21.3,5){\line(1,1){1.31}}
\put(25.2,-1.25){\line(1,1){1.8}}
\put(23.65,1.28){\wi}
\put(26.25,1.78){\tti}
\put(24.44,-0.24){\line(0,1){.7}}
\put(27.0,.5){\line(0,1){.45}}
\put(27.0,2.95){\line(0,1){1.2}}
\put(22.2,-0.24){\line(0,1){4.5}}
\put(21.3,-0.97){\line(0,1){6.}}
\put(21.9,-0.5){\seval}
\put(20.9,-1.5){\leval}
\end{picture}
$$
$$
$$
$$
$$
$$
$$
$$
$$
$$
$$
$$
$$
$$
$$
$$
\setlength{\unitlength}{0.5cm}
\begin{picture}(5,2)\thicklines
\put(.1,4.45){\text{\qquad\,\,\,$=$}}
\put(8,14.15){\text{\tiny{coev($V$)}}}
\put(8,15.25){\text{\tiny{coev($W$)}}}
\put(2.1,7.85){\lcoev}
\put(5.2,10.25){\scoev}
\put(5.07,9.01){\lbraidi}
\put(9.15,9.25){\u}
\put(7.7,8.4){\line(0,1){1.8}}
\put(9.88,8.4){\line(0,1){.45}}
\put(9.88,9.98){\line(0,1){.23}}
\put(7.7,12.31){\line(0,1){0.9}}
\put(9.88,12.31){\line(0,1){0.9}}
\put(6.7,7.5){\line(0,1){6.1}}
\put(10.9,7.6){\line(0,1){6}}
\put(9.4,4.85){\line(-1,1){2.7}}
\put(9.4,2.45){\line(0,1){2.4}}
\put(12,4.15){\line(-1,1){4.3}}
\put(9.8,6.55){\line(1,1){1.1}}
\put(8.92,7.5){\line(1,1){.95}}
\put(7.88,6.58){\line(1,1){.71}}
\put(8.78,5.68){\line(1,1){.71}}
\put(7.88,6.58){\line(1,1){.71}}
\put(7.2,4.24){\line(1,1){1.21}}
\put(6.3,5){\line(1,1){1.31}}
\put(10.2,-1.25){\line(1,1){1.8}}
\put(8.65,1.28){\wi}
\put(11.25,1.78){\tti}
\put(9.44,-0.24){\line(0,1){.7}}
\put(12.0,.5){\line(0,1){.45}}
\put(12.0,2.95){\line(0,1){1.2}}
\put(7.2,-0.24){\line(0,1){4.5}}
\put(6.3,-0.97){\line(0,1){6.}}
\put(6.9,-0.5){\seval}
\put(5.9,-1.5){\leval}
\put(15.9,4.45){\text{$=$}}

\put(17.1,6.15){\lcoev}
\put(20.9,7.95){\scoev}
\put(23.4,10.4){\line(0,1){.5}}
\put(25.6,10.4){\line(0,1){.5}}
\put(21.7,8.5){\line(0,1){3.4}}

\put(24.4,5.85){\line(-1,1){2.7}}
\put(24.4,3.45){\line(0,1){2.4}}
\put(27.35,10.45){\line(-1,1){2.4}}
\put(23.05,7.67){\braid}
\put(22.8,9.7){\u}

\put(23.6,8.25){\line(0,1){1.1}}
\put(23.8,6.88){\line(1,1){2.5}}
\put(22.88,7.58){\line(1,1){.71}}
\put(22.2,5.24){\line(1,1){1.21}}
\put(21.3,6){\line(1,1){1.31}}
\put(25.2,-0.25){\line(1,1){2.0}}
\put(23.65,2.28){\wi}
\put(26.55,2.78){\tti}

\put(24.44,.76){\line(0,1){.7}}
\put(27.2,1.7){\line(0,1){.25}}
\put(27.35,3.95){\line(0,1){4.7}}
\put(22.2,.76){\line(0,1){4.5}}
\put(21.3,.03){\line(0,1){6.}}
\put(21.9,0.5){\seval}
\put(20.9,-0.5){\leval}
\end{picture}
\vspace{0.5cm}
$$
$$
$$
$$
$$
$$
$$
$$\qquad\qquad\qquad
\setlength{\unitlength}{0.5cm}
\begin{picture}(5,2)\thicklines
\put(.1,3.45){\text{$=$}}
\put(2.1,5.15){\lcoev}
\put(5.9,6.95){\scoev}
\put(8.45,7.5){\line(0,1){2.4}}
\put(10.6,9.4){\line(0,1){.5}}
\put(6.7,7.5){\line(0,1){3.4}}

\put(9.4,4.85){\line(-1,1){1.1}}
\put(9.4,2.45){\line(0,1){2.4}}
\put(12.35,9.45){\line(-1,1){2.4}}
\put(8.05,6.67){\braid}
\put(5.57,1.4){\u}

\put(8.8,5.88){\line(1,1){2.5}}

\put(7.2,4.24){\line(1,1){1.21}}
\put(4.18,4.77){\braid}
\put(10.2,-1.25){\line(1,1){2.0}}
\put(8.65,1.28){\wi}
\put(11.55,1.78){\tti}
\put(9.44,-0.24){\line(0,1){.7}}
\put(12.2,.7){\line(0,1){.25}}
\put(12.35,2.95){\line(0,1){4.7}}
\put(7.2,-0.24){\line(0,1){4.5}}
\put(6.3,-0.97){\line(0,1){2.}}
\put(6.32,2.15){\line(0,1){3.25}}
\put(6.35,5.4){\line(1,1){1.}}
\put(6.9,-0.5){\seval}
\put(5.9,-1.5){\leval}

\put(15.9,3.45){\text{$=$}}
\put(16.1,5.15){\lcoev}
\put(19.9,6.95){\scoev}
\put(22.45,7.5){\line(0,1){2.4}}
\put(24.6,7.67){\line(0,1){2.23}}
\put(24.9,5.93){\line(0,1){4.97}}
\put(20.7,7.5){\line(0,1){3.4}}

\put(24.4,3.85){\line(-1,1){1.1}}
\put(24.4,1.45){\line(0,1){2.4}}

\put(17.63,0.4){\u}
\put(23.8,4.88){\line(1,1){1.1}}
\put(22.8,5.93){\line(1,1){1.8}}
\put(22.2,3.24){\line(1,1){1.21}}
\put(20.5,3.65){\line(1,1){1.21}}
\put(18.38,3.42){\line(1,1){2.51}}
\put(18.18,4.77){\braid}
\put(19.4,3.82){\braidi}

\put(23.65,0.28){\wi}
\put(19.75,0.58){\tti}
\put(24.44,-1.24){\line(0,1){.7}}

\put(22.2,-1.24){\line(0,1){4.5}}
\put(18.35,-1.97){\line(0,1){2.}}
\put(20.5,-1.95){\line(0,1){1.7}}
\put(20.5,1.7){\line(0,1){1.95}}
\put(18.35,1.15){\line(0,1){2.27}}

\put(21.9,-1.5){\seval}
\put(18,-2.2){\seval}

\end{picture}
$$
$$
$$
$$
$$
$$
$$
$$
\setlength{\unitlength}{0.5cm}
\begin{picture}(5,2)\thicklines
\put(1.,2.3){\text {=}}
\put(1.42,3.38){\coev}
\put(4.76,.48){\line(0,1){7.01}}

\put(7.75,.9){\line(0,1){3.07}}
\put(7.8,-2.11){\line(0,1){1.03}}
\put(4.4,-2.5){\eval}
\put(7,-.24){\tti}
\put(4,-.24){\u}
\put(4.75,-2.11){\line(0,1){1.5}}

\put(5.5,4.76){\braidi}
\put(5.5,2.97){\braidi}
\put(6.1,3.36){\coev}
\put(12.5,4){\line(0,1){3.46}}
\put(9.48,4.){\line(3,-2){3.}}
\put(11.2,3.13){\line(3,2){1.3}}
\put(9.48,2){\line(3,2){1.2}}
\put(12.48,.9){\line(0,1){1.1}}
\put(12.5,-2.11){\line(0,1){1.03}}
\put(9.1,-2.5){\eval}
\put(11.7,-.24){\wi}
\put(9.45,-2.11){\line(0,1){4.1}}
\put(16.,2.3){\text {=}}
\put(25.05,8.7){\text{\tiny{coev($W$)}}}
\put(20.4,8.7){\text{\tiny{coev($V$)}}}
\put(16.42,3.38){\coev}
\put(19.76,3.95){\line(0,1){3.53}}

\put(19.75,3.95){\line(3,-2){3}}
\put(21.5,3.13){\line(3,2){1.25}}
\put(19.78,2){\line(3,2){1.2}}
\put(22.78,.9){\line(0,1){1.03}}
\put(22.8,-2.11){\line(0,1){1.03}}
\put(19.4,-2.5){\eval}
\put(22,-.24){\ti}
\put(19.75,-2.11){\line(0,1){4.1}}
\put(20.5,4.76){\braidi}
\put(20.5,2.97){\braidi}
\put(21.1,3.36){\coev}
\put(27.5,4){\line(0,1){3.46}}
\put(24.48,4.){\line(3,-2){3.}}
\put(26.2,3.13){\line(3,2){1.3}}
\put(24.48,2){\line(3,2){1.2}}
\put(27.48,.9){\line(0,1){1.1}}
\put(27.5,-2.11){\line(0,1){1.03}}
\put(24.1,-2.5){\eval}
\put(26.7,-.24){\wi}
\put(24.45,-2.11){\line(0,1){4.1}}
\put(1.,-5.3){\text {=\,R. H. S.}}
\put(5.8,-5.3){\text {$\square$}}

\end{picture}\,\,\,\,\,\qquad \,\qquad \,\,\quad
$$
$$
$$
$$

\begin{Lem}\label{mod88} \hspace{0.2cm}
For two objects $V$  and $W$ in  ${\cal D}$,
$$
{\rm trace}\bigl( \Psi_{W \otimes V}\circ\Psi_{V \otimes W}
\bigr)= \sum_{ \begin{subarray}\,\ab\in \,\text{basis of} \,\,
V\otimes W \text{and}\\ \bbi \b \text{commutes with }\aa \ai
\end{subarray}} \hat\eta (\eta \achl \aai\bi \bb \aa ) \
\hat\xi(\xi \achl \aa\ai )
$$
\end{Lem}
\textbf{Proof.} From theorem \ref{trc11}, we know that
\begin{equation}
{\rm trace}\bigl(  \Psi_{W \otimes V}\circ\Psi_{V \otimes W}
\bigr)= \sum_{ \begin{subarray}\,(\ab)\in \,\text{basis of} \,
V\otimes W
\end{subarray}}
 \widehat{(\ab)} \big( \Psi^{2} (\ab)\big)
 \,\,\,.
\end{equation}
From the definition of the ribbon map, we know that $\Psi \,
\bigl(\Psi(\ab)\bigr) \achl \| \ab\|=\xi \achl \bbb\otimes \eta
\achl \aaa$, so
\begin{equation*}
\begin{split}
\Psi \bigl(\Psi(\ab)&\bigr) =\big(\xi \achl \bbb\otimes \eta \achl
\aaa
\big) \achl \| \ab\|^{-1}\\
&=(\xi \achl \bbi \b \otimes \eta \achl \aai \a) \achl
\ai \bi \bb \aa \\
&=(\xi \achl \bbi \b ) \achl \big( \big\|\eta \achl \aaa \big\|
\actr \ai \bi \bb \aa \big) \otimes \eta \achl \aai \a
\ai \bi \bb \aa\\
&=\xi \achl \bbi \b \big(  \aaa \actr \ai \bi \bb \aa \big)
\otimes \eta \achl \aai
 \bi \bb \aa.
\end{split}
\end{equation*}
Put $\Psi \bigl(\Psi(\ab)\bigr) =\xi'\otimes\eta'$ and
$\widehat{\xi\otimes\eta}=\alpha\otimes\beta$, and then from lemma
\ref{tendual} we get
$$
(\widehat{\xi\otimes\eta})(\xi'\otimes\eta')\,=\, \big( \al \acbl
\tau (\langle\beta\rangle ,\langle\xi'\rangle \cdot
\langle\eta'\rangle) \big) (\eta') \ \big(\beta \acbl \tau
(\langle\xi'\rangle , \langle\eta'\rangle)^{-1} \big)(\xi')\,.
$$
As $\widehat{\xi\otimes\eta}$ is part of a dual basis, the last
expression can only be non-zero if $\|\xi'\|=\|\xi\|$ and
$\|\eta'\|=\|\eta\|$. A simple calculation shows that
$\|\eta'\|=\|\eta\|$ if and only if $\bbi \b$ commutes with $\aa
\ai$.  We use this to find
\begin{equation*}
\begin{split}
\aaa \actl \ai \bi \bb \aa &=\aai \bbi \b \a \aai \a
\ai \bi \bb \aa \\
&=\aai \a \aai \bbi \b \bi \bb \aa =\aai \a \ ,
\end{split}
\end{equation*}
and then
$$
\aaa \actr \ai \bi \bb \aa= \a\ai \bi \bb \aa \ai= \bi \bb \aa
\ai\ .
$$
Now using the formula for
$\widehat{\xi\otimes\eta}=\alpha\otimes\beta$ from lemma
\ref{tendual} gives the result.\quad$\square$

\begin {Lem}{}\hspace{0.1cm}Let $V$
 and $W$ be   objects
in ${\cal D}$. Then
  in terms of group characters:
$$
{\rm trace}\bigl( \Psi^{2}{_{V \otimes W}} \bigr)= \sum_{
\begin{subarray}\,u, v \in G,\,\, s, t \in M \,\,\text{and}\,\\
su\,\, \text{commutes with }vt \end{subarray}}
 \chi_{_{W_{u s}}}(s^{-1} t^{-1} v^{-1} s)\
\chi_{_{V_{v t}}}(u^{-1} s^{-1})\ .
$$
\end{Lem}
\textbf{Proof.}\hspace{0.5cm} This is more or less immediate from
 lemma \ref{mod88}.  Put $\aaa =u^{-1}s$ and $\bbb=v^{-1}t$
 and sum over basis elements of constant degree first.  \quad$\square$

\section{An example of a modular category}

Using the order of the indecomposable objects in  table (5), we
get $T$ to be a diagonal $32\times 32$ matrix whose diagonal
entries are taken from the table. As every indecomposable object
in our example is self-dual, the matrix $C$ is the $32\times 32$
identity matrix.

To find $S$, we calculate the trace of the double braiding ${\rm
trace}(\Psi_{VW}\circ\Psi_{WV})$.  We do this using the result
from \ref{mod88}, and split into different cases for the objects
$V$ and $W$, and move the points the characters are evaluated at
to the base
points for each orbit using \ref{moving}.  For example:\\
Case (1) $\otimes$ Case (1): (i.e. the orbit of
$W$ is $\{e\}$ and the orbit of $V$ is $\{e\}$)
$$
{\rm trace} (\Psi^{2})=\chi_{_{W_e}}(e)\chi_{_{V_e}}(e).
$$
Case (2) $\otimes$ Case (5): (i.e. the orbit of
$W$ is $\{a^3\}$ and the orbit of $V$ is $\{b, ba^2, ba^4\}$)
\begin{equation*}
\begin{split}
{\rm trace} (\Psi^{2}) &=\Big(\chi_{_{W_{a^3}}}(ba^2)+
\chi_{_{W_{a^3}}}(ba^4)+
\chi_{_{W_{a^3}}}(b)\Big)\chi_{_{V_b}}(a^3).
\end{split}
\end{equation*}
Case (5) $\otimes$ Case (3): (i.e. the orbit of
$W$ is $\{b, ba^2, ba^4\}$ and the orbit of $V$ is $\{a^2, a^4\}$)
\begin{equation*}
\begin{split}
{\rm trace} (\Psi^{2})=0.
\end{split}
\end{equation*}
Case (6) $\otimes$ Case (5): (i.e. the orbit of
$W$ is $\{ba, ba^3, ba^5\}$ and the orbit of $V$ is $\{b, ba^2, ba^4\}$)
\begin{equation*}
\begin{split}
{\rm trace} (\Psi^{2})
&=3\Big(\chi_{_{W_{ba}}}(ba^4)\chi_{_{V_b}}(ba^3)\Big).
\end{split}
\end{equation*}
From these cases we get $S$ to be the following $32\times 32$
symmetric matrix:

$$
\renewcommand{\baselinestretch}{1.0}{
\!\! S=\tiny{\left(\begin{array}{cccccccccccccccccccccccccccccccc}
   1 & 1 & 1 & 1 & 2 & 2 & 1 & 1 & 1 & 1 & 2 & 2 & 2 & 2 &
   2 & 2 \\
    1 & 1 & 1 & 1 & 2 & 2 & -1 & -1 &
    -1 & -1 & -2 & -2 & 2 & 2 & 2 & 2 \\
     1 & 1 & 1 & 1 & 2 & 2 & -1 &
    -1 & -1 & -1 & -2 & -2 & 2 & 2 & 2 & 2 \\
     1 & 1 & 1 & 1 & 2 & 2 & 1 & 1 & 1 & 1 & 2 & 2 &
   2 & 2 & 2 & 2 \\
    2 & 2 & 2 & 2 & 4 & 4 & -2 & -2 & -2 & -2 &
    -4 & -4 & -2 & -2 & -2 & -2  \\
     2 & 2 & 2 & 2 & 4 & 4 & 2 & 2 & 2 & 2 & 4 & 4 &
    -2 & -2 & -2 & -2 \\ 1 & -1 &
    -1 & 1 & -2 & 2 & 1 & -1 & -1 & 1 & -2 & 2 & 2 &
    -2 & 2 & -2  \\ 1 & -1 &
    -1 & 1 & -2 & 2 & -1 & 1 & 1 & -1 & 2 & -2 & 2 &
    -2 & 2 & -2 \\ 1 & -1 & -1 & 1 &
    -2 & 2 & -1 & 1 & 1 & -1 & 2 & -2 & 2 & -2 & 2 &
    -2  \\
     1 & -1 & -1 & 1 &
    -2 & 2 & 1 & -1 & -1 & 1 & -2 & 2 & 2 & -2 & 2 &
    -2 \\
     2 & -2 & -2 & 2 &
    -4 & 4 & -2 & 2 & 2 & -2 & 4 & -4 & -2 & 2 & -2 & 2\\
     2 & -2 &
    -2 & 2 & -4 & 4 & 2 & -2 & -2 & 2 & -4 & 4 & -2 & 2 &
    -2 & 2\\
     2 & 2 & 2 &
   2 & -2 & -2 & 2 & 2 & 2 & 2 & -2 & -2 & 4 & -2 &
    -2 & 4\\
     2 & 2 & 2 & 2 &
    -2 & -2 & -2 & -2 & -2 & -2 & 2 & 2 & -2 & -2 & 4 &
    -2\\
     2 & 2 & 2 &
   2 & -2 & -2 & 2 & 2 & 2 & 2 & -2 & -2 & -2 & 4 & -2 &
    -2\\
     2 & 2 & 2 &
   2 & -2 & -2 & -2 & -2 & -2 & -2 & 2 & 2 & 4 & -2 &
    -2 & 4\\
     2 & 2 &
   2 & 2 & -2 & -2 & 2 & 2 & 2 & 2 & -2 & -2 & -2 &
    -2 & 4 & -2\\
     2 & 2 & 2 & 2 &
    -2 & -2 & -2 & -2 & -2 & -2 & 2 & 2 & -2 & 4 & -2 &
    -2\\
     2 & -2 &
    -2 & 2 & 2 & -2 & 2 & -2 & -2 & 2 & 2 & -2 & 4 & 2 &
    -2 & -4\\
     2 & -2 &
    -2 & 2 & 2 & -2 & -2 & 2 & 2 & -2 & -2 & 2 & -2 & -4 &
    -2 & 2\\
     2 &
    -2 & -2 & 2 & 2 & -2 & 2 & -2 & -2 & 2 & 2 & -2 &
    -2 & 2 & 4 & 2\\
     2 &
    -2 & -2 & 2 & 2 & -2 & -2 & 2 & 2 & -2 &
    -2 & 2 & 4 & 2 & -2 & -4\\
     2 & -2 &
    -2 & 2 & 2 & -2 & 2 & -2 & -2 & 2 & 2 & -2 & -2 & -4 &
    -2 & 2\\
     2 & -2 &
    -2 & 2 & 2 & -2 & -2 & 2 & 2 & -2 & -2 & 2 &
    -2 & 2 & 4 & 2 \\
     3 & -3 & 3 &
    -3 & 0 & 0 & 3 & -3 & 3 &
    -3 & 0 & 0 & 0 & 0 & 0 & 0 \\
     3 & -3 & 3 &
    -3 & 0 & 0 & 3 & -3 & 3 &
    -3 & 0 & 0 & 0 & 0 & 0 & 0\\
     3 & -3 & 3 &
    -3 & 0 & 0 & -3 & 3 &
    -3 & 3 & 0 & 0 & 0 & 0 & 0 & 0\\
     3 &
    -3 & 3 & -3 & 0 & 0 & -3 & 3 &
    -3 & 3 & 0 & 0 & 0 & 0 & 0 & 0 \\
     3 & 3 &
    -3 & -3 & 0 & 0 & 3 & 3 & -3 &
    -3 & 0 & 0 & 0 & 0 & 0 & 0\\
     3 & 3 & -3 &
    -3 & 0 & 0 & -3 &
    -3 & 3 & 3 & 0 & 0 & 0 & 0 & 0 & 0\\
     3 & 3 & -3 & -3 & 0 & 0 & 3 & 3 & -3 &
    -3 & 0 & 0 & 0 & 0 & 0 & 0\\
     3 & 3 & -3 &
    -3 & 0 & 0 & -3 &
    -3 & 3 & 3 & 0 & 0 & 0 & 0 & 0 & 0\\
 \end{array}\right.}}
$$
$$
\,\,\,\,\, \qquad \renewcommand{\baselinestretch}{1.0}{
\tiny{\left.
\begin{array}{cccccccccccccccccccccccccccccccc}
   2 & 2 & 2 & 2 & 2 & 2 & 2 & 2 & 3 & 3 & 3 & 3 &
   3 & 3 & 3 & 3 \\
     2 & 2 & -2 & -2 &
    -2 & -2 & -2 & -2 & -3 & -3 & -3 &
    -3 & 3 & 3 & 3 & 3 \\
      2 & 2 & -2 &
    -2 & -2 & -2 & -2 & -2 & 3 & 3 & 3 & 3 & -3 & -3 &
    -3 &
    -3 \\
      2 & 2 & 2 & 2 & 2 & 2 & 2 & 2 & -3 &
    -3 & -3 & -3 & -3 & -3 & -3 &
    -3 \\
     -2 &
    -2 & 2 & 2 & 2 & 2 & 2 & 2 & 0 & 0 & 0 & 0 & 0 & 0 & 0 &
   0 \\
    -2 & -2 & -2 & -2 & -2 & -2 & -2 &
    -2 & 0 & 0 & 0 & 0 & 0 & 0 & 0 & 0 \\
      2 & -2 & 2 & -2 & 2 & -2 & 2 &
    -2 & 3 & 3 & -3 & -3 & 3 & -3 & 3 & -3 \\
      2 & -2 & -2 & 2 & -2 & 2 & -2 & 2 & -3 &
    -3 & 3 & 3 & 3 & -3 & 3 & -3 \\
      2 & -2 & -2 & 2 & -2 & 2 & -2 & 2 & 3 & 3 & -3 &
    -3 & -3 & 3 & -3 & 3 \\
    2 & -2 & 2 & -2 & 2 & -2 & 2 & -2 & -3 &
    -3 & 3 & 3 & -3 & 3 & -3 & 3 \\
    -2 & 2 & 2 & -2 & 2 & -2 & 2 &
    -2 & 0 & 0 & 0 & 0 & 0 & 0 & 0 & 0 \\
     -2 & 2 & -2 & 2 & -2 & 2 &
    -2 & 2 & 0 & 0 & 0 & 0 & 0 & 0 & 0 & 0 \\
      -2 & -2 & 4 & -2 & -2 & 4 & -2 &
    -2 & 0 & 0 & 0 & 0 & 0 & 0 & 0 & 0 \\
     -2 & 4 & 2 & -4 & 2 & 2 &
    -4 & 2 & 0 & 0 & 0 & 0 & 0 & 0 & 0 & 0 \\
     4 & -2 & -2 & -2 & 4 & -2 &
    -2 & 4 & 0 & 0 & 0 & 0 & 0 & 0 & 0 & 0 \\
    -2 & -2 & -4 & 2 & 2 &
    -4 & 2 & 2 & 0 & 0 & 0 & 0 & 0 & 0 & 0 & 0 \\
      -2 & 4 & -2 & 4 & -2 & -2 & 4 &
    -2 & 0 & 0 & 0 & 0 & 0 & 0 & 0 & 0 \\
     4 & -2 & 2 & 2 & -4 & 2 & 2 &
    -4 & 0 & 0 & 0 & 0 & 0 & 0 & 0 & 0 \\
     -2 & 2 & 4 & 2 & -2 & -4 &
    -2 & 2 & 0 & 0 & 0 & 0 & 0 & 0 & 0 & 0 \\
     4 & 2 & 2 & -2 & -4 &
    -2 & 2 & 4 & 0 & 0 & 0 & 0 & 0 & 0 & 0 & 0 \\
     -2 & -4 & -2 & -4 &
    -2 & 2 & 4 & 2 & 0 & 0 & 0 & 0 & 0 & 0 & 0 & 0 \\
     -2 & 2 & -4 &
    -2 & 2 & 4 & 2 &
    -2 & 0 & 0 & 0 & 0 & 0 & 0 & 0 & 0 \\
     4 & 2 & -2 & 2 & 4 & 2 & -2 &
    -4 & 0 & 0 & 0 & 0 & 0 & 0 & 0 & 0 \\
     -2 & -4 & 2 & 4 & 2 & -2 & -4 &
    -2 & 0 & 0 & 0 & 0 & 0 & 0 & 0 & 0 \\
      0 & 0 & 0 & 0 & 0 & 0 & 0 &
   0 & 3 & -3 & 3 & -3 & 3 & -3 & -3 & 3 \\
    0 & 0 & 0 & 0 & 0 & 0 & 0 &
   0 & -3 & 3 & -3 & 3 & -3 & 3 & 3 & -3 \\
     0 & 0 & 0 & 0 & 0 & 0 &
   0 & 0 & 3 & -3 & 3 & -3 & -3 & 3 & 3 & -3 \\
    0 & 0 & 0 & 0 & 0 & 0 &
   0 & 0 & -3 & 3 & -3 & 3 & 3 & -3 & -3 & 3 \\
     0 & 0 & 0 & 0 & 0 & 0 & 0 &
   0 & 3 & -3 & -3 & 3 & 3 & 3 & -3 & -3 \\
     0 & 0 & 0 & 0 & 0 &
   0 & 0 & 0 & -3 & 3 & 3 & -3 & 3 & 3 & -3 &
    -3 \\
      0 & 0 & 0 & 0 & 0 & 0 & 0 &
   0 & -3 & 3 & 3 & -3 & -3 & -3 & 3 & 3 \\
     0 & 0 & 0 & 0 & 0 &
   0 & 0 & 0 & 3 & -3 & -3 & 3 & -3 & -3 & 3 & 3
 \end{array}\right)}}
$$
$$
$$
 Now it is possible to check that the matrices $S$, $T$ and
$C$ satisfy the following relations:
\[
S^2 \,=\,\frac{1}{12}\,(ST)^3, \qquad CS \,=\, SC, \qquad CT \,=\,
TC\ .
\]
\section {An equivalence of tensor categories}
In this section we will generalize some results of \cite{BGM}
which considered group doublecross products, i.e. a group $X$
factoring into two subgroups $G$ and $M$.
\begin{Def}\cite{BGM}
 For the doublecross product group $X=GM$ there is a quantum double
$D(X)=k(X) >\!\!\!\acl\, kX$  which has the following operations
\begin{equation*}
\begin{split}
(\delta_y\otimes x)(\delta_{y^{'}}\otimes x^{'})=\delta_{{x^{-1}}yx,y^{'}}(\delta_y\otimes x{x^{'}}),
\,\,\,\Delta(\delta_y\otimes x)=\sum_{ \begin{subarray}\,
ab=y \end{subarray}}\delta_a\otimes x\otimes\delta_b\otimes x\\
1=\sum_{ \begin{subarray}\,
y \end{subarray}}\delta_y\otimes e,\,\,\,\epsilon(\delta_y\otimes x)
=\delta_{y,e},\,\,\,S(\delta_y\otimes x)=
\delta_{x^{-1}y^{-1}x}\otimes x^{-1},\,\,\,\\
 (\delta_y\otimes x)^{*}
=\delta_{x^{-1}yx}\otimes x^{-1},\,\,\,R=\sum_{ \begin{subarray}\,
x,z \end{subarray}}\delta_x\otimes e \otimes\delta_z\otimes x.\,\,\,\,\,\,\,\,
\end{split}
\end{equation*}
The representations of $D(X)$ are given by $X$-graded left $kX$-modules.
  The $kX$ action will be denoted by $\dot \acr$
 and the grading by $|\!|\!|.|\!|\!|$.  The grading and $X$ action are
 related by
\begin{equation}
|\!|\!|x \dot \acr \xi|\!|\!|=x|\!|\!| \xi |\!|\!| x^{-1}\,\,\,
\,\,\,x\in X,\,\,\,\xi \in V,
\end{equation}
and the action of\, $(\delta_y\otimes x) \in D(X)$\, is given by
 \begin{equation}
(\delta_y\otimes x)\dot \acr \xi=\delta_{y,|\!|\!|x \dot \acr \xi|\!|\!|}x \dot \acr \xi.
\end{equation}
\end{Def}

\begin{Proposition}
There is a functor $\chi$ from ${\cal D}$ to the category of representations of
$D(X)$ given by the following: As vector spaces, $\chi (V)$ is the same
as $V$, and $\chi$ is the identity map.  The $X$-grading
$|\!|\!|.|\!|\!|$ on $\chi (V)$ and the action of $us\in kX$ are
defined by
\begin{eqnarray*}
|\!|\!|\chi(\eta)|\!|\!| &=& \a^{-1} \aa \,\,\,\,\,\,\,\,{\text{for}}\,\,\, \eta \in V,
\\
us \dot \acr \chi(\eta) &=& \chi \Big( \big( (s\acl\aa^{-1})
\acbr \eta \big) \acbl u^{-1}\Big) , \,\,\,\,\,\, s \in M \, \qquad u
\in G\ .  \end{eqnarray*} A morphism $\phi:V\to W$ in ${\cal D}$ is
sent to the morphism $\chi(\phi):\chi(V)\to \chi(W)$ defined by
$\chi(\phi)(\chi(\xi))=\chi(\phi(\xi))$.
\end{Proposition}
\textbf{Proof.} First we show that $\dot \acr$ is an action, i.e.\
$
vt \dot \acr\big(us \dot \acr \chi(\eta)\big)=vtus \dot \acr \chi(\eta)
$
for all $s,t \in M$ and $u,v\in G$.  Note that
\begin{equation*}
\begin{split}
vt \dot \acr\big(us \dot \acr \chi(\eta)\big)&=vt \dot \acr \chi \big(
\big( (s\acl\aa^{-1}) \acbr \eta \big) \acbl u^{-1}\big)\\
&= \chi\big( \big( (t\acl\aabar^{-1}) \acbr \bar\eta \big) \acbl
v^{-1}\big),
\end{split}
\end{equation*}
where $\eb= \big( (s\acl\aa^{-1})
\acbr \eta \big) \acbl u^{-1}$.  On the other hand we have
 $$vtus = v(t \acr u) \tau(t \acl u ,s)\big( (t \acl u)\cdot s\big), $$
 where $v(t \acr u) \tau(t \acl u ,s) \in G$ and
$ (t \acl u)\cdot s \in M $, so
\begin{equation*}
\begin{split}
vtus \dot \acr \chi(\eta)&=\chi\Big(\Big(\big(( (t \acl u)\cdot s)\acl
\aa^{-1}\big)\acbr \eta\Big) \acbl \tau(t \acl u ,s)^{-1}(t \acr
u)^{-1} v^{-1} \Big)\ .
\end{split}
\end{equation*}
We need to show that
\begin{equation}\label{pupi}
\begin{split}
(t\acl\aabar^{-1}) \acbr \bar\eta &=\Big(\big(( (t \acl u)\cdot s)\acl
\aa^{-1}\big)\acbr \eta\Big) \acbl \tau(t \acl u ,s)^{-1}(t \acr
u)^{-1}\\
&=\Big(\big(( t \acl u(s\acr \aa^{-1}))\cdot (s\acl \aa^{-1})\big)\acbr
\eta\Big) \acbl \tau(t \acl u ,s)^{-1}(t \acr u)^{-1}
\end{split}
\end{equation}
  Put $\bar s =s\acl\aa^{-1}$
and $\epr=\bar s\, \acbr \eta$ which give $\eb=\epr \acbl u^{-1}$.
Then using the connections between the gradings
and actions, 
$$
\aabar =|\epr \acbl u^{-1} |=(\ap \acr u^{-1})^{-1} \aap u^{-1}.
$$
Putting $\bar t =t \acl u \aap^{-1}$, the left hand side of (\ref{pupi}) will
become
\begin{equation*}
\begin{split}
(t\acl\aabar^{-1}) \acbr \bar\eta &=\Big(t \acl u \aap^{-1} (\ap \acr u^{-1})\Big)
 \acbr (\epr \acbl u^{-1})\\
&=\Big(\bar t\acl (\ap \acr u^{-1})\Big) \acbr (\epr \acbl u^{-1})\\
&= (\,\bar t \,\acbr \epr)\acbl \Big((\bar t\,  \acl \aap) \acr u^{-1}\Big)\ .
\end{split}
\end{equation*}
Now, from (\ref{pupi}) and the fact that $ (t \acr u )^{-1}=(\bar t\,
\acl \aap) \acr u^{-1} $, we only need to show that
\begin{eqnarray}\label{pupiii}\bar t \,\acbr \epr=\Big(\big(( t \acl
u(s\acr \aa^{-1})) \cdot (s\acl \aa^{-1})\big)\acbr \eta\Big)
\acbl \tau(t \acl u ,s)^{-1}\ .
\end{eqnarray}
From the formula for the composition of the $M$ `action' the right hand
side of (\ref{pupiii}) becomes
$\bar p\, \acbr(\bar s\, \acbr \eta)\,=\,\bar p\, \acbr \epr$,
where ${\bar p}^{'}=t \acl u(s\acr \aa^{-1})$ and
$\bar p = {\bar p}^{'} \acl \tau(\bar s ,\a)
\tau(\langle \bar s\, \acbr \eta\rangle ,\bar s\, \acl \aa)^{-1} $.  We have
used the fact that $\tau(t \acl u ,s)=\tau({\bar p}^{'} \acl (\bar
s\, \acr \aa) ,\bar s\, \acl \aa)$.
Now we just have to prove that $ \bar p = \bar t$.  Because $\tau(\bar s
,\a)^{-1}(\bar s\, \acr \aa)= \tau(\langle \bar s\, \acbr \eta\rangle
, \bar s\, \acl \aa)^{-1} |\bar s\, \acbr \eta|$ and knowing that
$(\bar s\, \acr \aa)=(s\acr \aa^{-1})^{-1}$, we can write $\bar p$ as
follows
\begin{equation*}
\begin{split}
\bar p
&={\bar p}^{'} \acl (\bar s\, \acr \aa)|\bar s\, \acr \eta|^{-1}\\
&=t \acl u(s\acr \aa^{-1})(s\acr \aa^{-1})^{-1}\aap^{-1}\\
&=t \acl u\aap^{-1}= \bar t\ .
\end{split}
\end{equation*}
   Next we show that $|\!|\!|us \dot\acr \chi (\eta) |\!|\!| =us\,|\!|\!|
   \chi (\eta) |\!|\!|\,(us)^{-1}$ where $u\in G$
   and $s\in M$.
\begin{equation*}
\begin{split}
|\!|\!|x \dot\acr \chi (\eta) |\!|\!|&=|\!|\!|\chi \Big( \big( (s\acl\aa^{-1})
\acbr \eta \big) \acbl u^{-1}\Big) |\!|\!|\\
&=\langle  \epr \acbl u^{-1} \rangle^{-1} \, | \epr \acbl u^{-1} |\,
=\,u \ap^{-1} \aap u^{-1}\\
&=u \langle \bar s\,
\acbr \eta \rangle^{-1} |\bar s\,\acbr \eta | u^{-1}\,=\,u ( \bar s\,
\acl \aa ) \ai \aa ( \bar s\,\acl \aa )^{-1} u^{-1}\\
&=u s \ai \aa s^{-1} u^{-1}\ .\quad \square
\end{split}
\end{equation*}

\begin{theorem}\hspace{0.1cm}
The functor $\chi$ is invertible.
\end{theorem}
\textbf{Proof.}  We have already proved
 in the previous proposition that the  $X$-grading
 $|\!|\!|.|\!|\!|$ and the action $\dot \acr$ \, give a representation
 of $D(X)$, so we only need to show that $\chi$ is a 1-1
correspondence,  which we do by giving its inverse $\chi^{-1}$ as
the following:  Let $W$ be a representation of $D(X)$, with $kX$
action $\dot \acr$ and X-grading $|\!|\!|.|\!|\!|$.
  Define a $D$ representation as follows:  $\chi^{-1}(W)$
will be the same as $W$ as a vector space.  There will be $G$ and
$M$ gradings given by the factorization
$$
{|\!|\!|\xi|\!|\!|}^{-1}={|\chi^{-1}(\xi)|}^{-1}
\langle\chi^{-1}(\xi)\rangle , \,\,\,\xi \in W
,\,\langle\chi^{-1}(\xi)\rangle \in M, \,\,|\chi^{-1}(\xi)|\, \in
G.
$$
The action of $s\in M$ and $u \in G$ are given by
$$
s \acbr \chi^{-1}(\xi)=\chi^{-1}\big ((s \acl|\chi^{-1}(\xi)|)
\dot \acr \xi\big)\,\,, \,\, \chi^{-1}(\xi) \acbl
u=\chi^{-1}(u^{-1} \dot \acr\xi)
$$
Checking the rest is left to the reader.\quad$\square$

\begin{Proposition} For $\delta_y\otimes x \in{\cal D}$,
$
\chi\big(\xi\achl(\delta_y\otimes x)\big)=
\delta_{y,\bbb}\, x^{-1}\dacr \chi(\xi).
$
\end{Proposition}
\textbf{Proof.} \,\,\, Starting with the left hand side,
\begin{equation*}
\begin{split}
\chi\big(\xi\achl(\delta_y\otimes x)\big)&=
\chi(\delta_{y,\bbb}\,\xi\achl x )\,=\,\delta_{y,\bbb}\,\chi(\xi\achl x ).
\end{split}
\end{equation*}
Putting $x=us$ for $u\in G$ and $s\in M$,
\begin{equation*}
\begin{split}
\xi\achl x=\xi\achl us=(\xi\acbl u) \achl s=
\big((s^{L}\acl u^{-1}\bbi)\acbr \xi \big)
\acbl(s^{L} \acr u^{-1} )^{-1} \tau(s^{L},s).
\end{split}
\end{equation*}
Now put ${\bar u}=\tau(s^{L},s)^{-1} (s^{L} \acr u^{-1}
)$ and $\bar s =s^{L} \acl u^{-1}$.  Then
\begin{equation*}
\begin{split}
\chi\big(\xi\achl(\delta_y\otimes x)\big)&=\delta_{y,\bbb}\,\chi \Big(
\big( (\bar s\acl\bb^{-1}) \acbr \xi \big) \acbl {\bar u}^{-1}\Big)\\
&=\delta_{y,\bbb}\,\bar u \bar s \dacr \chi(\xi)\\
&=\delta_{y,\bbb}\,\tau(s^{L},s)^{-1}
(s^{L} \acr u^{-1} ) (s^{L} \acl u^{-1}) \dacr \chi(\xi)\\
&=\delta_{y,\bbb}\, s^{-1} {s^{L}}^{-1} s^{L} u^{-1} \, \dacr \chi(\xi)\\
&=\delta_{y,\bbb}\, (us)^{-1} \, \dacr \chi(\xi)\,=\,
\delta_{y,\bbb}\, x^{-1} \, \dacr \chi(\xi).\quad \square
\end{split}
\end{equation*}

\begin{Proposition} \label{dmap} Define a map $\psi :D \longrightarrow D(X)$
by\,
 $
\psi (\delta_y\otimes x)=\delta_{x^{-1}yx}\otimes x^{-1}\ .
  $
    Then $\psi$ satisfies the equation
$
\chi(\xi \achl (\delta_y\otimes x))=\psi (\delta_y\otimes x) \dacr \chi(\xi)\ .
$
\end{Proposition}
\textbf{Proof.} \,\,\,
Use the previous proposition.\quad$\square$

\medskip The reader will recall that $D$ is in general a non-trivially
associated algebra (i.e.\ it is only associative in the category
${\cal D}$ with its non-trivial associator).  Thus, in general, it
can not be isomorphic to  $D(X)$, which really is associative. In
general, $\psi$ can not be an algebra map.

\begin{Proposition} For $a$ and $b$  elements of the algebra $D$ in the
category ${\cal D}$,
\begin{equation*}
\begin{split}
\psi (b)\psi (a)=\psi (ab)\Big(\sum_{ \begin{subarray}\,
y\in Y \end{subarray}}\delta_{y}\otimes \tau(\langle a\rangle,\langle
b\rangle)^{-1}\Big).
\end{split}
\end{equation*}
\end{Proposition}
\textbf{Proof.} \,\,by {\ref{dmap}} we have
\begin{equation*}
\begin{split}
\chi\big((\xi \achl a)\achl b\big)&=\psi (b) \dacr
 \chi(\xi \achl a)\\
&=\psi (b)\dacr \big(\psi (a)\dacr \chi(\xi)\big)\\
&=\psi (b) \psi (a)\dacr \chi(\xi).
\end{split}
\end{equation*}
But also, where $f=\sum_{y}\delta_{y}\otimes\tau(\langle
a\rangle,\langle b\rangle)$,
\begin{equation*}
\begin{split}
\chi\big((\xi \achl a)\achl b\big)&=\chi\Big(\big(\xi \achl
\tilde \tau(\|a \|,\|b \|)\big)\achl ab\Big)\\
&=\,\psi (ab)\dacr \chi\big(\xi \achl \tilde
\tau(\|a \|,\|b \|)\big)
\,=\, \psi (ab)\dacr \chi\big(\xi \achl \tilde
\tau(\langle a\rangle,\langle b\rangle)\big)
\\
&=\psi (ab)\psi (f)\dacr \chi(\xi  ).\quad \square
\end{split}
\end{equation*}

\begin{Def}\,\,\,
 Let \,$ V$ and $W$ be objects of ${\cal D}$.    The map
   $c :\chi(V)\otimes \chi (W)\longrightarrow \chi(V\otimes W)$
 is defined by:
\begin{equation*}
\begin{split}
c \big(\chi(\eta)\otimes \chi (\xi)\big)= \chi\Big(\big((\b \acl
\aai )\acbr \eta \big)\otimes \xi\Big)\ .
\end{split}
\end{equation*}
\end{Def}

\begin{Proposition} The map $c$, defined above, is a $D(X)$
module map, i.e.
\begin{equation*}
\begin{split}
|\!|\!|c \big(\chi(\eta)\otimes \chi (\xi)\big)|\!|\!|&=
|\!|\!|\chi(\eta)\otimes \chi (\xi)|\!|\!|\
,\\
 x \dacr c \big(\chi(\eta)\otimes \chi (\xi)\big)&=
c\Big(x \dacr  \big(\chi(\eta)\otimes \chi (\xi)\big)\Big)\
\qquad\forall x\in X\ .
\end{split}
\end{equation*}
\end{Proposition}
\textbf{Proof.} \,\,\, We will begin with the grading first.
It is known that
$$
|\!|\!|\chi(\eta)\otimes \chi (\xi)|\!|\!|=
|\!|\!|\chi(\eta)|\!|\!||\!|\!| \chi (\xi)|\!|\!|=
\ai \aa \bi \bb.
$$
But on the other hand we know, from the definition of $c$, that
\begin{equation*}
\begin{split}
|\!|\!|c \big(\chi(\eta)\otimes \chi (\xi)\big)|\!|\!|&=
|\!|\!|\chi\big((\b \acl
\aai )\acbr \eta \otimes \xi\big)|\!|\!|\\
&=\big\langle (\b \acl
\aai )\acbr \eta \otimes \xi \big\rangle^{-1} \big|(\b \acl
\aai )\acbr \eta \otimes \xi\big|\\
&=\bi \langle \bar \eta \rangle^{-1} |\bar \eta | \bb \\
&=\bi \langle \bar s \acbr \eta  \rangle^{-1}
 |\bar s \acbr \eta| \bb\\
&=\bi (\bar s \acl \aa) \ai \aa (\bar s \acl \aa)^{-1} \bb\\
&=\ai \aa \bi \bb ,
\end{split}
\end{equation*}
where $\bar s=\b \acl
\aai$ and $\bar \eta=(\b \acl
\aai )\acbr \eta=\bar s \acbr \eta$, which gives the result.

For the $G$ action, we know from the definitions that
\begin{equation*}
\begin{split}
u \dacr \big(\chi(\eta)\otimes \chi (\xi)\big)&=
\chi(\eta\acbl u^{-1})\otimes \chi (\xi \acbl u^{-1}),\\
c \Big( u \dacr \big(\chi(\eta)\otimes \chi (\xi)\big)\Big)
&=\chi \Big( \big( (\langle \xi \acbl u^{-1}  \rangle \acl
|\eta \acbl u^{-1}|^{-1}) \acbr (\eta \acbl u^{-1})\big)
\otimes (\xi \acbl u^{-1})\Big).
\end{split}
\end{equation*}
By using the properties of
the $G$ and $M$ gradings,
 \begin{equation*}
\begin{split}
\langle \xi \acbl u^{-1}  \rangle \acl
|\eta \acbl u^{-1}|^{-1}&= (\b \acl u^{-1})
 \acl u \aai (\langle\eta\rangle \acr u^{-1} )\\
&= \b \acl  \aai (\langle\eta\rangle \acr  u^{-1} )\\
(\langle \xi \acbl u^{-1}  \rangle \acl
|\eta \acbl u^{-1}|^{-1}) \acbr (\eta \acbl u^{-1})&=
\big((\b \acl  \aai)\acl (\langle\eta\rangle \acr
u^{-1} )\big)\acbr (\eta \acbl u^{-1})\\
&=\big((\b \acl  \aai)\acbr \eta \big)\acbl  \Big(\big((\b \acl \aai)\acl\aa\big)\acr u^{-1}\Big)\\
&=\big((\b \acl  \aai)\acbr \eta \big)\acbl   (\b  \acr u^{-1}).
\end{split}
\end{equation*}
Now we can write
 \begin{equation}\label{uupq3}
     c \Big( u \dacr \big(\chi(\eta)\otimes \chi
 (\xi)\big)\Big) =\chi \Big( \big((\b \acl \aai)\acbr \eta \big)\acbl
 (\b \acr u^{-1}) \otimes (\xi \acbl u^{-1})\Big).
\end{equation}
On the other hand,
\begin{equation*}
\begin{split}
u \dacr c \big(\chi(\eta)\otimes \chi (\xi)\big)&=u \dacr \chi\Big(\big((\b \acl
\aai )\acbr \eta \big)\otimes \xi\Big)\\
&=\chi\Big(\Big(\big((\b \acl \aai )\acbr \eta \big)\otimes
\xi\Big)\acbl u^{-1}\Big)\ ,
\end{split}
\end{equation*}
which gives the same as (\ref{uupq3}).

Now we show that  $c$ preserves the $M$ action. For $s\in M$,
\begin{equation*}
\begin{split}
s \dacr  \big(\chi(\eta)\otimes \chi (\xi)\big)&= \chi\big((s \acl
\aai )\acbr \eta \big)\otimes \chi\big((s \acl
\bbi )\acbr \xi \big)\\
c \Big( s \dacr \big(\chi(\eta)\otimes \chi (\xi)\big)\Big)&=
\chi\Big(\big(\big\langle ( s \acl \bbi) \acbr \xi \big\rangle
\acl \big|(s \acl
\aai )\acbr \eta \big|^{-1}\big) \acbr \big((s \acl
\aai )\acbr \eta \big) \\
&\,\,\,\,\,\,\,\,\,\,\,\,\,\,\,\,\,\otimes \big((s \acl
\bbi )\acbr \xi \big) \Big).
\end{split}
\end{equation*}
Using the `action' property for $\acbr$, we get
\begin{equation*}
\begin{split}
\big(\big\langle ( s \acl \bbi) \acbr \xi \big\rangle
\acl \big|(s \acl
\aai )\acbr \eta \big|^{-1} \big) \acbr \big((s \acl
\aai )\acbr \eta \big) =\big( (p^{'} \cdot \bar t) \acbr \eta \big) \acbl
\tau \big(p^{'}\acl (\bar t \acr \aa ),\bar t \acl\aa \big)^{-1}\ ,
\end{split}
\end{equation*}
where $\bar t =s \acl\aai$ and
\begin{equation*}
\begin{split}
 p^{'}&=\big\langle ( s \acl \bbi) \acbr \xi \big\rangle
\acl \big|\bar t \acbr \eta \big|^{-1}\, \tau \big(\langle \bar
 t \acbr \eta \rangle,\bar t \acl\aa \big)\ \tau \big( \bar
 t ,\a \big)^{-1}.
\end{split}
\end{equation*}
But using the connections between the grading and the
actions, we know that $\big|\bar t \acbr \eta \big|^{-1}=
(\bar t \acr \aa )^{-1}\,\tau \big( \bar
 t ,\a \big) \tau \big(\langle \bar
 t \acbr \eta \rangle,\bar t \acl\aa \big)^{-1}$,\, so
\begin{equation*}
\begin{split}
 p^{'}&=\big\langle ( s \acl \bbi) \acbr \xi \big\rangle
\acl (\bar t \acr \aa )^{-1}\\
&=\big\langle ( s \acl \bbi) \acbr \xi \big\rangle
\acl\big((s \acl\aai) \acr \aa \big)^{-1}\\
&=\big\langle ( s \acl \bbi) \acbr \xi \big\rangle
\acl (s \acr\aai)  .
\end{split}
\end{equation*}
Substituting in the equation above gives
\begin{equation*}
\begin{split}
\Big(\big\langle ( s \acl \bbi) \acbr \xi \big\rangle
\acl \big|(s \acl
\aai )\acbr \eta \big|^{-1} \Big) \acbr \big((s \acl
\aai )\acbr \eta \big) \qquad \qquad\qquad \qquad\qquad \qquad \qquad\\
=\Big( \Big(\big(\big\langle ( s \acl \bbi) \acbr \xi \big\rangle \acl
(s \acr\aai)\big) \cdot (s \acl\aai)\Big) \acbr \eta \Big) \acbl \tau
\big(\big\langle ( s \acl \bbi) \acbr \xi \big\rangle,s \big)^{-1}\\
=\Big( \Big(\big(\big\langle ( s \acl \bbi) \acbr \xi \big\rangle
\cdot s\big) \acl\aai\Big) \acbr \eta \Big) \acbl \tau
\big(\big\langle ( s \acl \bbi) \acbr \xi \big\rangle,s
\big)^{-1}\qquad \qquad\quad\, \\
=\Big( \Big(\big( ( s \acl \bbi) \cdot\b \big) \acl\aai\Big) \acbr
\eta \Big) \acbl \tau \big(\big\langle ( s \acl \bbi) \acbr \xi
\big\rangle,s \big)^{-1}.\qquad \qquad\qquad
\end{split}
\end{equation*}
On the other hand, we know that
\begin{equation*}
\begin{split}
s \dacr c \big(\chi(\eta)\otimes \chi (\xi)\big)&=s \dacr
 \chi\Big(\big((\b \acl
\aai )\acbr \eta \big)\otimes \xi\Big)\\
&=s \dacr \chi\big(\bar \eta\otimes \xi\big)
=\chi\big((s \acl|\bar \eta\otimes \xi|^{-1}) \acbr (\bar \eta\otimes \xi)\big),
\end{split}
\end{equation*}
where $\bar \eta=(\b \acl \aai )\acbr \eta$.  Next we calculate
\begin{equation*}
\begin{split}
 |\bar \eta\otimes \xi|&=\tau(\langle \bar \eta
\rangle, \b)^{-1} |\bar \eta| \bb\ ,\\
s \acl|\bar \eta\otimes \xi|^{-1}&=s \acl\, \bbi |\bar \eta|^{-1}\tau(\langle \bar \eta
\rangle, \b)\ .
\end{split}
\end{equation*}
If we put $\,\bar s =s \acl\, \bbi |\bar \eta|^{-1}$, then
\begin{equation*}
\begin{split}
 (s \acl|\bar \eta\otimes \xi|^{-1}) \acbr (\bar \eta\otimes
\xi)&=\big(\bar s \acl\, \tau(\langle \bar \eta
\rangle, \b)\big) \acbr (\bar \eta \otimes \xi)\\
&=(\bar s \acbr \bar \eta)\acbl\,\tau(\bar s \acl |\bar \eta|,
 \b)\,\tau(\langle (\bar s \acl |\bar \eta|)\acbr \xi
\rangle, \bar s \acl|\bar \eta|\bb )^{-1}\\
&\,\,\,\,\,\,\,\,\,\otimes (\bar s \acl |\bar \eta|)\acbr \xi\\
&=(\bar s \acbr \bar \eta)\acbl\,\tau( s \acl \bbi,
 \b)\,\tau(\langle (s \acl \bbi)\acbr \xi
\rangle,  s  )^{-1}\\
&\,\,\,\,\,\,\,\,\,\otimes (s \acl \bbi)\acbr \xi\ .
\end{split}
\end{equation*}
Using the `action' property again,
\begin{equation*}
\begin{split}
\bar s \acbr \bar \eta &=(s \acl\, \bbi |\bar \eta|^{-1})
\acbr \big((\b \acl
\aai )\acbr \eta \big)\\
&=\Big( \big(q' \cdot(\b \acl \aai)\big) \acbr \eta \Big) \acbl\, \tau
\big(q' \acl \big((\b \acl \aai) \acr \aa \big), \b \big)^{-1}\\
&=\Big( \big(q' \cdot(\b \acl \aai)\big) \acbr \eta \Big) \acbl\, \tau
\big(q' \acl (\b \acr \aai )^{-1}, \b \big)^{-1}\ ,
\end{split}
\end{equation*}
where
\begin{equation*}
\begin{split}
q'&=(s \acl\, \bbi |\bar \eta|^{-1}) \acl
\,\tau \big( \big \langle (\b \acl \aai)
\acr \eta \big \rangle, \b \big)\,
\tau \big( \b \acl \aai , \a \big)^{-1}\\
&=(s \acl\, \bbi ) \acl (\b \acr \aai )\ ,
\end{split}
\end{equation*}
as
$$
|\bar \eta|^{-1}= \big( (\b \acl \aai) \acr \aa \big)^{-1}
\tau \big( \b \acl \aai , \a \big)\,
\tau \big( \big \langle (\b \acl \aai)
\acr \eta \big \rangle, \b \big)^{-1}\ .
$$
Hence substituting with the value of $q'$ we get
\begin{equation*}
\begin{split}
\bar s \acbr \bar \eta &=
\Big( \Big(\big((s \acl\, \bbi ) \acl (\b \acr \aai ) \big) \cdot(\b \acl
\aai)\Big) \acbr \eta \Big) \acbl\, \tau \big((s \acl\, \bbi ), \b
\big)^{-1}\\
&=\Big( \Big(\big((s \acl\, \bbi ) \cdot  \b \big) \acl \aai\Big) \acbr
\eta \Big) \acbl\, \tau \big(s \acl\, \bbi , \b \big)^{-1}\ ,
\end{split}
\end{equation*}
giving the required result
\begin{equation*}
\begin{split}
(\bar s \acbr \bar \eta)\acbl\,\tau( s \acl \bbi,
 \b)\,\tau\big(\langle (s \acl \bbi)\acbr \xi
\rangle,  s  \big)^{-1} \qquad \qquad\qquad\qquad\qquad\qquad\qquad\qquad\\
\qquad =\Big( \Big(\big((s \acl\, \bbi ) \cdot \b \big) \acl \aai\Big)
\acbr \eta \Big) \acbl \tau\big(\big\langle (s \acl \bbi)\acbr \xi
\big \rangle, s \big)^{-1}\ .\square
\end{split}
\end{equation*}


\begin{thebibliography}{AAA}

\bibitem{Bak}
{\sc Bakalov B. \& Kirillov A.,}
{\sl Lectures on Tensor Categories and Modular Functors.
\, American Mathematical Society,  $2000$. }

\bibitem{BNT}
{\sc Beggs  E. J.,} {\sl Making non-trivially associated tensor
categories from left coset representatives. Journal of Pure and
Applied Algebra, vol 177 , 5 - 41, $2003$.}

\bibitem{BGM}
{\sc Beggs  E. J., Gould J. D. and Majid S.,}
{\sl  Finite group factorizations and braiding. \, J. Algebra, vol 181 no. 1, 112 - 151, $1996$.}

\bibitem{BM1}
{\sc Beggs  E. J. and Majid S.,}
{\sl Quasitriangular and differential structures on bicrossproduct Hopf algebras.
 \, J. Algebra, vol 219 no. 2, 682 - 727, $1999$.}




\bibitem{Grove}
{\sc Grove L. C.,}
{\sl Groups and Characters.
\, Wiley-Interscience,  $1997$. }

\bibitem{GM}
{\sc Gurevich D. I. and  Majid S.,}
{\sl
Braided groups of  Hopf algebras obtained  by twisting.
\,Pacific  J.\,  Math, vol 162, 27 - 44, $1994$.}



\bibitem{Maj1}
{\sc Majid S.,}
{\sl Physics for algebraists:  Non-commutative and non-cocommutative
 Hopf algebras by bicrossproduct construction.
\, J.  Algebra, vol 130, 17 - 64, $1990$.  From PhD Thesis, Harvard, 1988.}

\bibitem{Maj2}
{\sc Majid S.,}
{\sl The quantum double as quantum mechanics.
\, J.  Geom. Phys., vol 13, 169 - 202, $1994$. }

\bibitem{MajBook}
{\sc Majid S.,}
{\sl Foundations of Quantum Group Theory.
\, Cambridge University Press,  $1995$. }



\bibitem{Tak1}
{\sc Takeuchi  M. ,}
{\sl Matched pairs of groups and bismash products of Hopf algebras.
 \, Commun. Alg., vol. 9,  no. 8, 841-882 , $1981$.}

\bibitem{Tak2}
{\sc Takeuchi  M. ,}
{\sl Modular Categories and  Hopf Algebras.
 \, Journal of Algebra, vol. 243,  no. 2, 631-643, $2001$.}

\bibitem{TWood}
{\sc Thomas A. D. \& Wood G. V.,}
{\sl Group Tables.
\, Shiva Publishing Ltd., $1980$.  }

\bibitem{TWen}
{\sc Turaev  V. \& Wenzl  H.,}
{\sl Semisimple and modular categories from link invariants.
 \, Math. Ann., vol. 309,  no. 3, 411-461, $1997$.}




\end{thebibliography}
\end{document}